\documentclass[11pt]{article}
\usepackage{amsmath, amssymb, amsthm, amsbsy, ifthen, graphicx}

\date{}
\title{\vspace{-1cm}Maximizing the number of $q$-colorings}
\author{
Po-Shen Loh \thanks{Department of Mathematics,
Princeton University, Princeton, NJ 08544. E-mail:
ploh@math.princeton.edu.
Research supported in part by a Fannie and John Hertz Foundation Fellowship, an NSF Graduate
Research Fellowship, and a Princeton Centennial Fellowship.}
\and
Oleg Pikhurko \thanks{Department of Mathematical Sciences, Carnegie Mellon University, Pittsburgh, PA 15123.  E-mail: pikhurko@andrew.cmu.edu.
Research supported in part by NSF grant DMS-0758057.}
\and
Benny Sudakov \thanks{Department of Mathematics, UCLA, Los Angeles, CA 90095. E-mail:
bsudakov@math.ucla.edu.
Research supported in part by NSF CAREER award DMS-0812005, and a USA-Israeli BSF grant.
}
}

\newtheorem{theorem}{Theorem}[section]
\newtheorem{lemma}[theorem]{Lemma}
\newtheorem{fact}{Fact}[section]
\newtheorem{corollary}[theorem]{Corollary}
\newtheorem{inequality}[theorem]{Inequality}
\newtheorem{proposition}[theorem]{Proposition}
\newtheorem{definition}[theorem]{Definition}

\newcommand{\obj}{\text{\sc obj}}
\newcommand{\objS}{\text{\sc obj}^*}

\newcommand{\opt}{\text{\sc opt}}
\newcommand{\optS}{\text{\sc opt}^*}

\newcommand{\ec}{\text{\sc e}}
\newcommand{\vc}{\text{\sc v}}

\newcommand{\feas}{\text{\sc Feas}}
\newcommand{\feasS}{\text{\sc Feas}^*}

\newcommand{\blpha}{{\boldsymbol\alpha}}

\newcommand{\bnu}{{\boldsymbol\nu}}

\begin{document}
\maketitle

\begin{abstract}
  Let $P_G(q)$ denote the number of proper $q$-colorings of a graph
  $G$.  This function, called the \emph{chromatic polynomial}\/ of
  $G$, was introduced by Birkhoff in 1912, who sought to attack the
  famous four-color problem by minimizing $P_G(4)$ over all planar
  graphs $G$.  Since then, motivated by a variety of applications,
  much research was done on minimizing or maximizing $P_G(q)$ over
  various families of graphs.

  In this paper, we study an old problem of Linial and Wilf, to find
  the graphs with $n$ vertices and $m$ edges which maximize the number
  of $q$-colorings.  We provide the first approach which enables one
  to solve this problem for many nontrivial ranges of parameters.
  Using our machinery, we show that for each $q \geq 4$ and
  sufficiently large $m < \kappa_q n^2$ where $\kappa_q \approx 1/(q
  \log q)$, the extremal graphs are complete bipartite graphs minus
  the edges of a star, plus isolated vertices.  Moreover, for $q=3$,
  we establish the structure of optimal graphs for all large $m \leq
  n^2/4$, confirming (in a stronger form) a conjecture of Lazebnik
  from 1989.
\end{abstract}

\section{Introduction}

The fundamental combinatorial problem of graph coloring is as
ancient as the cartographer's task of coloring a map without using the
same color on neighboring regions.  In the context of general graphs,
we say that an assignment of a color to every vertex is a \emph{proper
  coloring}\/ if no two adjacent vertices receive the same color, and
we say that a graph is \emph{$q$-colorable}\/ it has a proper coloring
using only at most $q$ different colors.  

The problem of counting the number $P_G(q)$ of $q$-colorings of a
given graph $G$ has been the focus of much research over the past
century.  Although it is already NP-hard even to determine whether
this number is nonzero, the function $P_G(q)$ itself has very
interesting properties.  $P_G(q)$ was first introduced by Birkhoff
\cite{Birkhoff-1912}, who proved that it is always a polynomial in
$q$.  It is now called the \emph{chromatic polynomial}\/ of $G$.
Although $P_G(q)$ has been studied for its own sake (e.g., Whitney
\cite{Whitney} expressed its coefficients in terms of graph theoretic
parameters), perhaps more interestingly there is a long history of
diverse applications which has led researchers to minimize or maximize
$P_G(q)$ over various families of graphs.  In fact, Birkhoff's
original motivation for investigating the chromatic polynomial was to
use it to attack the famous four-color theorem.  Indeed, one way to
show that every planar graph is 4-colorable is to minimize $P_G(4)$
over all planar $G$, and show that the minimum is nonzero.  In this
direction Birkhoff \cite{Birkhoff-1930} proved the tight lower bound
$P_G(q) \geq q(q-1)(q-2)(q-3)^{n-3}$ for all $n$-vertex planar graphs
$G$ when $q \geq 5$, later conjecturing with Lewis in
\cite{Birkhoff-Lewis-1946} that it extended to $q=4$ as well.

Linial \cite{Linial} arrived at the problem of minimizing the
chromatic polynomial from a completely different motivation.  The
worst-case computational complexity of determining whether a
particular function $f : V(G) \rightarrow \mathbb{R}$ is a proper
coloring (i.e., satisfies $f(x) \neq f(y)$ for every pair of adjacent
vertices $x$ and $y$) is closely related to the number of \emph{acyclic
  orientations}\/ of a graph, which equals $|P_G(-1)|$, obtained by
substituting $q = -1$ into the formal polynomial expression of
$P_G(q)$.
% (See Manber-Tompa \cite{Manber-Tompa} and Stanley \cite{Stanley} for
% details.)
Lower bounding the worst-case complexity therefore corresponds to
minimizing $|P_G(-1)|$ over the family $\mathcal{F}_{n,m}$ of graphs
with $n$ vertices and $m$ edges.  Linial showed that that
surprisingly, for any $n,m$ there is a graph which
\emph{simultaneously}\/ minimizes each $|P_G(q)|$ over
$\mathcal{F}_{n,m}$, for \emph{every}\/ integer $q$.  This graph is
simply a clique $K_k$ with an additional vertex adjacent to $l$
vertices of the $K_k$, plus $n-k-1$ isolated vertices, where $k,l$ are
the unique integers satisfying $m = {k \choose 2} + l$ with $k > l
\geq 0$.  At the end of his paper, Linial posed the problem of
maximizing $P_G(q)$ over all graphs in $\mathcal{F}_{n,m}$.

Around the same time, Wilf arrived at exactly that maximization
problem while analyzing the \emph{backtrack}\/ algorithm for finding a
proper $q$-coloring of a graph (see \cite{Bender-Wilf, Wilf}).
Although this generated much interest in the problem, it was only
solved in sporadic cases.  The special case $q=2$ was completely
solved for all $m,n$, by Lazebnik in \cite{Lazebnik-q23}.  For $q \geq
3$, the only nontrivial pairs $m,n$ for which extremal graphs were known
corresponded to the number of vertices and edges in the Tur\'an graph
$T_r(n)$, which is the complete $r$-partite graph on $n$ vertices with
all parts of size either $\lfloor n/r \rfloor$ or $\lceil n/r \rceil$.
In this vein, Lazebnik \cite{Lazebnik-largeq} proved that $T_r(n)$ is
optimal for very large $q = \Omega(n^6)$, and proved with Pikhurko and Woldar
\cite{LPW} that $T_2(2k)$ is optimal when $q=3$ and asymptotically
optimal when $q=4$.

Outside these isolated cases, very little was known for general $m,n$.
Although many upper and lower bounds for $P_G(q)$ were proved by
various researchers \cite{Byer, Lazebnik-q23, Lazebnik-bounds, Liu},
these bounds were widely separated.  Even the $q=3$ case resisted
solution: twenty years ago, Lazebnik \cite{Lazebnik-q23} conjectured
that when $m \leq n^2/4$, the $n$-vertex graphs with $m$ edges which
maximized the number of 3-colorings were complete bipartite graphs
minus the edges of a star, plus isolated vertices.  Only very
recently, Simonelli \cite{Simonelli} managed to make some progress on
this conjecture, verifying it under the additional very strong
assumption that all optimal graphs are already bipartite.

Perhaps part of the difficulty for general $m,n,q$ stems from the fact
that the maximal graphs are substantially more complicated than the
minimal graphs that Linial found.  For number-theoretic reasons, it is
essentially impossible to explicitly construct maximal graphs for
general $m,n$.  Furthermore, even their coarse structure depends on
the density $\frac{m}{n^2}$.  For example, when $\frac{m}{n^2}$ is
small, the maximal graphs are roughly complete bipartite graphs, but
after $\frac{m}{n^2} > \frac{1}{4}$, the maximal graphs become
tripartite.  At the most extreme density, when $m,n$ correspond to the
Tur\'an graph $T_q(n)$, the unique maximal graph is obviously the
complete $q$-partite graph.  Therefore, in order to tackle the general
case of this problem, one must devise a unified approach that can
handle all of the outcomes.

In this paper, we propose such an approach, developing the machinery
that one might be able to use to determine the maximal graphs in many
nontrivial ranges of $m,n$.  Our methodology can be roughly outlined
as follows.  We show, via Szemer\'edi's Regularity Lemma, that the
asymptotic solution to the problem reduces to a certain
quadratically-constrained linear program in $2^q - 1$ variables.  For any
given $q$, this task can in principle be automated by a computer code
that symbolically solves the optimization problem, although a more
sophisticated approach was required to solve this for all $q$.  Our
solutions to the optimization problem then give us the approximate
structure of the maximal graphs.  Finally, we use various local
arguments, such as the so-called ``stability'' approach introduced by
Simonovits \cite{Simonovits}, to refine their structure into precise
results.

We successfully applied our machinery to solve the Linial-Wilf problem
for many nontrivial ranges of $m,n$, and $q \geq 3$.  In particular,
for $q=3$, our results confirm a stronger form of Lazebnik's
conjecture when $m$ is large.  In addition, for each $q \geq 4$ we
show that for all densities $\frac{m}{n^2}$ up to approximately
$\frac{1}{q \log q}$, the extremal graphs are also complete bipartite
graphs minus a star.  In order to state our results precisely, we need
the following definition.

\begin{definition}
  \label{def:semi-complete}
  Let $a \leq b$ be positive integers.  We say that $G$ is a
  \textbf{semi-complete subgraph of $\boldsymbol{K_{a,b}}$} if the
  number of missing edges $E(K_{a,b}) \setminus E(G)$ is less than
  $a$, and they form a star (i.e., they share a common endpoint $v$
  which we call the \textbf{center}).  If $v$ belongs to the larger
  side of $K_{a,b}$, then we also say that $G$ is \textbf{correctly
    oriented}.
\end{definition}

Define the constant $\kappa_q = \left( \sqrt{\frac{\log (q/(q-1))}{\log
      q}} + \sqrt{\frac{\log q}{\log (q/(q-1))}} \right)^{-2} \approx
\frac{1}{q \log q}$.  All logarithms here and in the rest of the paper
are in base $e \approx 2.718$.  In the following theorems, we write
$o(1)$ to represent a quantity that tends to zero as $m,n \rightarrow
\infty$.

\begin{theorem}
  \label{thm:main:sparse}
  For every fixed integer $q \geq 3$, and any $\kappa < \kappa_q$, the
  following holds for all sufficiently large $m$ with $m \leq \kappa n^2$.  Every
  $n$-vertex graph with $m$ edges which maximizes the number of
  $q$-colorings is a semi-complete subgraph (correctly oriented if $q
  \geq 4$) of some $K_{a,b}$, plus isolated vertices, where $a =
  (1+o(1))\sqrt{m \cdot \log \frac{q}{q-1} / \log q}$ and $b =
  (1+o(1)) \sqrt{m \cdot \log q / \log \frac{q}{q-1}}$.  The
  corresponding number of $q$-colorings is $q^n
  e^{(-c+o(1))\sqrt{m}}$, where $c = 2\sqrt{\log \frac{q}{q-1} \log
    q}$.
\end{theorem}

\noindent \textbf{Remark.}\, The part sizes of the maximal graphs
above all have the ratio roughly $\log q / \log \frac{q}{q-1}$.  The
constant $\kappa_q$ corresponds to the density $m/n^2$ at which the
number of isolated vertices becomes $o(n)$ in the optimal
construction.

\vspace{3mm}

For 3 colors, we can push our argument further, beyond the density
$\kappa_3$.  Now, due to the absence of isolated vertices, a rare
exception occurs, which requires us to include an additional
possibility.  Here, a ``pendant edge'' means that a new vertex is
added, along with a single edge between it and any other vertex in the
graph.  Proposition \ref{prop:pendant-necessary} shows that this
outcome is in fact necessary.

\begin{theorem}
  \label{thm:main:q=3}
  The following holds for all sufficiently large $m \leq n^2/4$.
  Every $n$-vertex graph with $m$ edges and the maximum number of
  3-colorings is either \textbf{(i)} a semi-complete subgraph of some
  $K_{a,b}$, plus isolated vertices if necessary, or \textbf{(ii)} a
  complete bipartite graph $K_{a,b}$ plus a pendant edge.
  Furthermore:
  \begin{itemize}
  \item If $m \leq \kappa_3 n^2$, then $a = (1+o(1))\sqrt{m \cdot
      \frac{\log 3/2}{\log 3}}$ and $b = (1+o(1))\sqrt{m \cdot
      \frac{\log 3}{\log 3/2}}$.  The corresponding number of
    colorings is $3^n e^{-(c+o(1))\sqrt{m}}$, where $c = 2\sqrt{ \log
      \frac{3}{2} \cdot \log 3 }$.
  \item If $\kappa_3 n^2 \leq m \leq \frac{1}{4} n^2$, then $a =
    (1+o(1))\frac{n-\sqrt{n^2-4m}}{2}$ and $b =
    (1+o(1))\frac{n+\sqrt{n^2-4m}}{2}$.  The corresponding number of
    colorings is $2^{b + o(n)}$.
  \end{itemize}
\end{theorem}

We also considered another conjecture of Lazebnik (see, e.g.,
\cite{LPW}), that the Tur\'an graphs $T_r(n)$ are always extremal
when $r \leq q$.  Building upon the techniques in \cite{LPW} that
answered the $r=2, q=3$ case, we confirmed this conjecture for large
$n$ and $r = q-1$.

\begin{theorem}
  \label{thm:exact:turan}
  Fix an integer $q \geq 4$.  For all sufficiently large $n$, the
  Tur\'an graph $T_{q-1}(n)$ has more $q$-colorings than any other
  graph with the same number of vertices and edges.
\end{theorem}

We close by mentioning some related work.  Tomescu \cite{Tomescu-max,
  Tomescu-q3-conn, Tomescu-min, Tomescu-1975, Tomescu-hamiltonian,
  Tomescu-conn-planar, Tomescu-2conn, Tomescu-blocks} and Dohmen
\cite{Dohmen-1, Dohmen-2} considered the problem of maximizing or
minimizing the number of $q$-colorings of $G$ given some other
parameters, such as chromatic number, connectedness, planarity, and
girth.  Wright \cite{Wright} asymptotically determined the total
number of $q$-colored labeled $n$-vertex graphs with $m$ edges, for
the entire range of $m$; this immediately gives an asymptotic
approximation for the \emph{average}\/ value of $P_G(q)$ over all
labeled $n$-vertex graphs with $m$ edges.

Graph coloring is also a special case of a homomorphism problem, and
as we will discuss in our concluding remarks, our approach easily
extends to that more general setting.  Recall that a graph
homomorphism $\phi : G \rightarrow H$ is a map from the vertices of
$G$ to those of $H$, such that adjacent vertices in $G$ are mapped to
adjacent vertices in $H$.  Thus, the number of $q$-colorings of $G$ is
precisely the number of homomorphisms from $G$ to $K_q$.  Another
interesting target graph $H$ is the two-vertex graph consisting of a
single edge, plus a loop at one vertex.  Then, the number of
homomorphisms is precisely the number of independent sets in $G$, and
the problem of estimating that number given some partial information
about $G$ is motivated by various questions in statistical physics and
the theory of partially ordered sets.  Alon \cite{Alon-indep-sets}
studied the maximum number of independent sets that a $k$-regular
graph of order $n$ can have, and Kahn \cite{Kahn-hard-core,
  Kahn-dedekind} considered this problem under the additional
assumption that the $k$-regular graph is bipartite.  Galvin and Tetali
\cite{Galvin-Tetali} generalized the main result from
\cite{Kahn-hard-core} to arbitrary target graphs $H$.

Another direction of related research was initiated by the question of
Erd\H{o}s and Rothschild (see Erd\H{o}s \cite{Erdos-1974, Erdos-1992},
Yuster \cite{Yuster}, Alon, Balogh, Keevash, and Sudakov \cite{ABKS},
Balogh \cite{Balogh}, and others), about the maximum over all
$n$-vertex graphs of the number of $q$-edge-colorings (not necessarily
proper) that do not contain a monochromatic $K_r$-subgraph.  Our
method is somewhat similar to that in \cite{ABKS}, and these two
problems may be more deeply related than just a similarity in their
formulations.

\vspace{3mm}

The rest of this paper is organized as follows.  The next section
contains some definitions, and a formulation of the Szemer\'edi
Regularity Lemma.  In Section \ref{sec:reduction-to-opt}, we prove
Theorems \ref{thm:asymp-number} and \ref{thm:asymp-stability}, which
(asymptotically) reduce the general case of the problem to a
quadratically constrained linear program.  Then, in the next section
we solve the relevant instances of the optimization problem to give
approximate versions of our main theorems.  Sections
\ref{sec:exact:sparse} and \ref{sec:exact:q=3} refine these into the
precise forms of Theorems \ref{thm:main:sparse} and
\ref{thm:main:q=3}.  We prove Theorem \ref{thm:exact:turan} in Section
\ref{sec:exact:turan}.  The final section contains some concluding
remarks and open problems.

\section{Preliminaries}
\label{sec:preliminaries}

The following (standard) asymptotic notation will be utilized
extensively.  For two functions $f(n)$ and $g(n)$, we write $f(n) =
o(g(n))$ if $\lim_{n \rightarrow \infty} f(n)/g(n) = 0$, and $f(n) =
O(g(n))$ or $g(n) = \Omega(f(n))$ if there exists a constant $M$ such
that $|f(n)| \leq M|g(n)|$ for all sufficiently large $n$.  We also
write $f(n) = \Theta(g(n))$ if both $f(n) = O(g(n))$ and $f(n) =
\Omega(g(n))$ are satisfied.

We will use $[q]$ to denote the set $\{1, 2, \ldots, q\}$, and
$2^{[q]}$ to denote the collection of all of its subsets.  As
mentioned in the introduction, the \emph{Tur\'an graph} $T_q(n)$ is
the complete $r$-partite graph on $n$ vertices with all parts of size
either $\lfloor n/r \rfloor$ or $\lceil n/r \rceil$.

Given two graphs with the same number of vertices, their \emph{edit
  distance}\/ is the minimum number of edges that need to be added or
deleted from one graph to make it isomorphic to the other.  We say
that two graphs are \emph{$d$-close}\/ if their edit distance is at
most $d$.

%\subsection{Szemer\'edi's Regularity Lemma}
%\label{sec:tools-regularity}

The rest of this section is devoted to formulating
%one of the most powerful tools in modern combinatorics, 
the celebrated Szemer\'edi Regularity Lemma.  This theorem roughly
states that every graph, no matter how large, can be approximated by
an object of bounded complexity, which corresponds to a union of a
bounded number of random-looking graphs.  To measure the randomness of
edge distribution, we use the following definition.  Let the edge
density $d(A, B)$ be the fraction $\frac{e(A, B)}{|A| |B|}$, where
$e(A, B)$ is the number of edges between $A$ and $B$.

\begin{definition}
  A pair $(X, Y)$ of disjoint subsets of a graph is
  \textbf{$\epsilon$-regular} if every pair of subsets $X' \subset X$
  and $Y' \subset Y$ with $|X'| \geq \epsilon |X|$ and $|Y'| \geq
  \epsilon |Y|$ has $|d(X', Y') - d(X, Y)| < \epsilon$.
\end{definition}

In this paper, we use the following convenient form of the Regularity
Lemma, which is essentially Theorem IV.5.$29'$ in the textbook
\cite{B-modern-graph-theory}.

\begin{theorem}
  \label{thm:regularity-lemma} For every $\epsilon > 0$, there is a
  natural number $M' = M'(\epsilon)$ such that \textbf{every} graph $G
  = (V, E)$ has a partition $V = \bigcup_{i=1}^M V_i$ with the
  following properties.  The sizes of the vertex clusters $V_i$ are as
  equal as possible (differing by at most 1), their number is between
  $1/\epsilon \leq M \leq M'$, and all but at most $\epsilon M^2$ of
  the pairs $(V_i, V_j)$ are $\epsilon$-regular.
\end{theorem}

\section{Reduction to an optimization problem}
\label{sec:reduction-to-opt}

In this section, we show that the solution of the following
quadratically constrained linear\footnote{Observe that the logarithms
  are merely constant multipliers for the variables $\alpha_A$.}
program answers our main problem asymptotically.

\vspace{3mm}

\noindent \textbf{Optimization Problem 1.}\, Fix an integer $q \geq 2$
and a real parameter $\gamma$.  Consider the following objective and
constraint functions:
\begin{displaymath}
  \obj(\blpha) := \sum_{A \neq \emptyset} \alpha_A \log |A|\,;
  \quad\quad\quad
  \vc(\blpha) := \sum_{A \neq \emptyset} \alpha_A,
  \quad
  \ec(\blpha) := \sum_{A \cap B = \emptyset} \alpha_A \alpha_B.
\end{displaymath}
The vector $\blpha$ has $2^q - 1$ coordinates $\alpha_A \in
\mathbb{R}$ indexed by the nonempty subsets $A \subset [q]$, and the
sum in $\ec(\blpha)$ runs over \emph{unordered}\/ pairs of disjoint
nonempty sets $\{A,B\}$.  Let $\feas(\gamma)$ be the \emph{feasible set}\/ of
vectors defined by the constraints $\blpha \geq 0$, $\vc(\blpha) = 1$,
and $\ec(\blpha) \geq \gamma$.  We seek to maximize $\obj(\blpha)$
over the set $\feas(\gamma)$, and we define $\opt(\gamma)$ to be this
maximum value, which exists by compactness.  We will write that the
vector $\blpha$ \emph{solves}\/ $\opt(\gamma)$ when both $\blpha \in
\feas(\gamma)$ and $\obj(\blpha) = \opt(\gamma)$.

\vspace{3mm}

\noindent \textbf{Note.}\, In the remainder of this paper, we will
write $\sum_A$ instead of $\sum_{A \neq \emptyset}$ because it is
clear from the definition of $\blpha$ that the empty set is excluded.

\vspace{3mm}

%Next, we show how a vector $\blpha$ with large $\obj(\blpha)$ produces
%a graph with many $q$-colorings.

%\vspace{2mm}

\noindent \textbf{Construction 1: $\boldsymbol{G_\alpha(n)}$.}\, Let
$n$ and $m$ be the desired numbers of vertices and edges, and let
$\blpha \in \feas(m/n^2)$ be a feasible vector.  Consider the
following $n$-vertex graph, which we call $G_\blpha(n)$.  Partition
the vertices into (possibly empty) clusters $V_A$ such that each
$|V_A|$ differs from $n \alpha_A$ by less than 1.  For every pair of
clusters $(V_A, V_B)$ which is indexed by disjoint subsets, place a
complete bipartite graph between the clusters.

\vspace{3mm}

Observe that any coloring that for each cluster $V_A$ uses only colors
from $A$ is a proper coloring.  Therefore, if all $n \alpha_A$
happened to be integers, then $G_\blpha(n)$ would have at least
$\prod_A |A|^{n\alpha_A} = e^{\obj(\blpha)n}$ colorings, and also
precisely $\ec(\blpha) n^2$ edges.  But we cannot simply apply
Construction 1 to the $\blpha$ that solves $\opt(m/n^2)$, because it
may happen that $G_\blpha(n)$ has fewer than $m$ edges if the entries
of $\blpha$ are not integer multiples of $1/n$.  Fortunately, the
shortfall cannot be substantial:

\begin{proposition}
  \label{prop:construction-asymp-edges}
  The number of edges in any $G_\blpha(n)$ differs from $\ec(\blpha)
  n^2$ by less than $2^q n$.  Also, for any other vector $\bnu$, the
  edit-distance between $G_\blpha(n)$ and $G_\bnu(n)$ is at most $\|
  \blpha - \bnu \|_1 n^2 + 2^{q+1} n$, where $\| \cdot \|_1$ is the
  $L^1$-norm.
\end{proposition}

The proof is elementary and routine, so we will defer it to Section
\ref{sec:pf-asymp-number-ii} so as not to interrupt this exposition.
To recover from the $O(n)$ edge deficit, we extend the construction
in the following way.

\vspace{3mm}

\noindent \textbf{Construction 2: $\boldsymbol{G_\blpha'(n)}$.}\, Let
$n$ and $m$ be the desired numbers of vertices and edges, and let
$\blpha \in \feas(m/n^2)$ be a feasible vector.  If $G_\blpha(n)$ from
Construction 1 already has at least $m$ edges, then set $G_\blpha'(n)
= G_\blpha(n)$.

Otherwise, $G_\blpha(n)$ is short by, say, $k$ edges, and $k = O(n)$
by Proposition \ref{prop:construction-asymp-edges}.  Let $V_A$ be its
largest cluster whose index $A$ is not a singleton.  Suppose first
that $|V_A| \geq 2\lceil \sqrt{k} \rceil$.  So far $V_A$ does not span
any edges, so we can add $k$ edges to $G_\blpha(n)$ by selecting two
disjoint subsets $U_1, U_2 \subset V_A$ of size $\lceil \sqrt{k}
\rceil$, and putting a $k$-edge bipartite graph between them.  Call
the result $G_\blpha'(n)$.

The last case is $|V_A| < 2\lceil \sqrt{k} \rceil$.  We will later
show that this only arises when the maximum number of colorings is
only $2^{o(n)}$, and this is already achieved by the Tur\'an graph
$T_q(n)$.  So, to clean up the statements of our theorems, we just
define $G_\blpha'(n) = T_q(n)$ here.

\subsection{Structure of asymptotic argument}
\label{sec:dense-structure-argument}

We are now ready to state our theorem, which shows that solutions to
Optimization Problem 1 produce graphs which asymptotically maximize
the number of $q$-colorings.

\begin{theorem}
  \label{thm:asymp-number}
  For any $\epsilon > 0$, the following holds for any sufficiently
  large $n$, and any $m$ less than or equal to the number of edges in
  the Tur\'an graph $T_q(n)$.
  \begin{description}

  \item[(i)] Every $n$-vertex graph with $m$ edges has fewer than
    $e^{(\opt(m/n^2) + \epsilon) n}$ proper $q$-colorings.

  \item[(ii)] Any $\blpha$ which solves $\opt(m/n^2)$ yields a graph
    $G_\blpha'(n)$ via Construction 2 which has at least $m$ edges and
    more than $e^{(\opt(m/n^2) - \epsilon) n}$ proper $q$-colorings.

  \end{description}

\end{theorem}

\noindent \textbf{Remark.}\, The number of colorings can only increase
when edges are deleted, so one may take an arbitrary $m$-edge subgraph
of $G_\blpha'(n)$ if one requires a graph with exactly $m$ edges.

\vspace{3mm}

The key ingredient in the proof of Theorem \ref{thm:asymp-number} is
Szemer\'edi's Regularity Lemma.  Part (ii) is routine, and full
details are given in Section \ref{sec:pf-asymp-number-ii}.  On the
other hand, the argument for part (i) is more involved, so we
highlight its structure here so that the reader does not get lost in
the details.  The proof breaks into the following claims.

\begin{description}
\item[Claim 1.]  For any $\delta > 0$, there exists $n_0$ such that
  the following holds for any graph $G = (V, E)$ with $n > n_0$
  vertices and $m$ edges.  The Regularity Lemma gives a special
  partition of the vertex set into sets $V_1$, \ldots, $V_M$ of almost
  equal size, where $M$ is upper bounded by a constant depending only
  on $\delta$.  Then, we may delete at most $\delta n^2$ edges of $G$
  in such a way that the resulting graph $G'$ has the following
  properties.
  \begin{description}
  \item[(i)] Each $G'[V_i]$ spans no edges.
  \item[(ii)] If $G'$ has any edges at all between two parts $V_i$
    and $V_j$, then in fact it has an edge between every pair of
    subsets $U \subset V_i$, $W \subset V_j$ with $|U| \geq \delta
    |V_i|$ and $W \geq \delta |V_j|$.
  \end{description}
  Note that since $G'$ is a subgraph of $G$, the number of
  $q$-colorings can only increase.

\item[Claim 2.] Let $\mathcal{C}_1$ be the set of colorings of $G'$.
  Then, if we keep only those colorings $\mathcal{C}_2 \subset
  \mathcal{C}_1$ with the property that in each $V_i$, any color is
  used either zero times or at least $\delta |V_i|$ times, we will
  still have $|\mathcal{C}_2| \geq e^{-c_\delta n} |\mathcal{C}_1|$.
  Here, $c_\delta$ is a constant which tends to zero with $\delta$.
  Now each coloring in $\mathcal{C}_2$ has the special property that
  whenever the same color appears on two parts $V_i$ and $V_j$, then
  there cannot be any edges between those entire parts.

\item[Claim 3.] By looking at which colors appear on each part $V_i$,
  we may associate each coloring with a map $[M] \rightarrow
  2^{[q]}$.  Let $\phi : [M] \rightarrow 2^{[q]}$ be a map which is
  associated with the maximum number of colorings in $\mathcal{C}_2$.
  Then, if we keep only those colorings $\mathcal{C}_3 \subset
  \mathcal{C}_2$ which give $\phi$, we still have $|\mathcal{C}_3|
  \geq 2^{-qM} |\mathcal{C}_2|$.

\item[Claim 4.] For every nonempty $A \subset [q]$, let $V_A$ be the
  union of those parts $V_i$ for which $\phi(i) = A$.  (These are the
  parts that in all colorings in $\mathcal{C}_3$ are colored using
  exactly colors from $A$.)  Define the vector $\blpha$ by setting
  each $\alpha_A = |V_A|/n$.  Then $G' \subset G_\blpha(n)$, and since
  $G'$ only differs from our original $G$ by at most $\delta n^2$
  edges, we also have $\blpha \in \feas(m/n^2 - \delta)$.  Thus:
  \begin{displaymath}
    |\mathcal{C}_3| 
    \ \leq \ \prod_A |A|^{|V_A|} 
    \ = \ e^{\obj(\blpha) n}
    \ \leq \ e^{\opt(m/n^2 - \delta) n} \, .
  \end{displaymath}

\item[Claim 5.] The function $\opt$ is uniformly continuous.  Thus,
  for an appropriate (sufficiently small) choice of $\delta > 0$, we
  have for all sufficiently large $n$ that
  \begin{displaymath}
    P_G(q)
    \ \leq \ P_{G'}(q)
    \ \leq \ e^{c_\delta n} \cdot 2^{qM} \cdot e^{\opt(m/n^2 - \delta) n}
    \ < \ e^{(\opt(m/n^2) + \epsilon) n} \, ,
  \end{displaymath}
  as desired.  (Recall that $P_G(q)$ is the number of $q$-colorings of $G$.)

\end{description}

By combining these five claims with an elementary analysis argument,
we also obtain a stability result, which roughly states that if a
graph has ``close'' to the optimal number of colorings, then it must
resemble a graph from Construction 1.  A stability result is very
useful, because the approximate structure later allows us to apply
combinatorial arguments to refine our asymptotic results into exact
results.  We quantify this in terms of the edit-distance, which we
defined in Section \ref{sec:preliminaries}.  Recall that we say that two
graphs are $d$-close when their edit distance is at most $d$.  We
prove the following theorem in Section \ref{sec:asymp-stability}.

% \begin{theorem}
%   \label{thm:asymp-stability}
%   For any $\epsilon > 0$, there exists $\delta > 0$ such that the
%   following holds for all sufficiently large $n$.  Let $G$ be an
%   $n$-vertex graph with $m$ edges and at least $e^{(\opt(m/n^2) -
%     \delta) n}$ proper $q$-colorings.  Then $G$ is $\epsilon
%   n^2$-close to some $G_\blpha(n)$ from Construction 1, for an
%   $\blpha$ which solves $\opt(\gamma)$ for some $|\gamma - m/n^2| \leq
%   \delta$.  
% \end{theorem}

\begin{theorem}
  \label{thm:asymp-stability}
  For any $\epsilon, \kappa > 0$, the following holds for all
  sufficiently large $n$.  Let $G$ be an $n$-vertex, graph with $m
  \leq \kappa n^2$ edges, which maximizes the number of $q$-colorings.
  Then $G$ is $\epsilon n^2$-close to some $G_\blpha(n)$ from
  Construction 1, for an $\blpha$ which solves $\opt(\gamma)$ for some
  $|\gamma - m/n^2| \leq \epsilon$ with $\gamma \leq \kappa$.
\end{theorem}

\noindent \textbf{Remark.}\, This theorem is only useful if the
resulting $\gamma$ falls within the range of densities for which the
solution of $\opt$ is known.  The technical parameter $\kappa$ is used
to keep $\gamma$ within this range.

\subsection{Finer resolution in the sparse case}
\label{sec:sparse-summary}

The Regularity Lemma is nontrivial only for graphs with positive edge
density (i.e., quadratic number of edges).  This typically presents a
serious and often insurmountable obstacle when trying to extend
Regularity-based results to situations involving sparse graphs.
Although much work has been done to develop sparse variants of the
Regularity Lemma, the resulting analogues are weaker and much more
difficult to apply.

Let us illustrate the issue by attempting to apply Theorem
\ref{thm:asymp-number} when $m = o(n^2)$.  Then, we find that the
maximum number of $q$-colorings of any $n$-vertex graph with $m$ edges
is $e^{cn + o(n)}$, where $c = \opt(0) = \log q$ is a constant entirely
determined by $q$.  Note that the final asymptotic is independent of
$m$, even if $m$ grows extremely slowly compared to $n^2$.  This is
because the key parameter was the density $m/n^2$, which already
vanished once $m = o(n^2)$.  Thus, the interesting question in the
sparse case is to distinguish between sparse graphs and very sparse
graphs, by looking inside the $o(n)$ error term in the exponent.

We are able to circumvent these difficulties by making the following
key observation which allows us to pass to a dense subgraph.  As it
turns out, every sparse graph which maximizes the number of
$q$-colorings has a nice structure: most of the vertices are isolated,
and \emph{all}\/ of the edges are contained in a subgraph which is
dense, but not too dense.  Section \ref{sec:sparse-proofs} contains
the following lemma's short proof, which basically boils down to a
comparison against the smallest Tur\'an graph with at least $m$ edges.

\begin{lemma}
  \label{lem:sparse-has-core}
  Fix an integer $q \geq 2$ and a threshold $\kappa > 0$.  Given any
  positive integer $m$, there exists an $n_0 = \Theta(\sqrt{m})$ with
  $m/n_0^2 \leq \kappa$ such that the following holds for any $n \geq
  n_0$.  In every $n$-vertex graph $G$ with $m$ edges, which maximizes
  the number of $q$-colorings, there is a set of $n_0$ vertices which
  spans all of the edges.
\end{lemma}

The fact that our graph is sparse becomes a benefit rather than a
drawback, because it allows us to limit the edge density from above by
any fixed threshold.  This is useful, because we can
completely solve the optimization problem for all densities below
$\kappa_q = \left( \sqrt{\frac{\log q/(q-1)}{\log q}} +
  \sqrt{\frac{\log q}{\log q/(q-1)}} \right)^{-2}$.  We will prove the
following proposition in Section \ref{sec:solve-opt-sparse}.

\begin{proposition}
  \label{prop:solve-opt-sparse}
  Fix an integer $q \geq 3$.  For any $0 \leq \gamma \leq
  \kappa_q$, the \textbf{unique} solution (up to a permutation of the
  ground set $[q]$) to $\opt(\gamma)$ has the following form.
  \begin{equation}
    \label{eq:sparse-opt-soln}
    \alpha_{\{1\}} = \sqrt{\gamma \cdot \log \frac{q}{q-1} \, / \, \log q},
    \quad\quad
    \alpha_{\{2, \ldots q\}} = \frac{\gamma}{\alpha_{\{1\}}},
    \quad\quad
    \alpha_{[q]} = 1 - \alpha_{\{1\}} - \alpha_{\{2, \ldots q\}},
  \end{equation}
  with all other $\alpha_A = 0$.  This gives $\opt(\gamma) = \log q -
  2\sqrt{\gamma \cdot \log \frac{q}{q-1} \cdot \log q}$.
\end{proposition}

Since we have the complete solution of the relevant instance of the
optimization problem, we can give explicit bounds when we transfer our
asymptotic results from the previous section to the sparse case.  We
can also explicitly describe the graph that approximates any optimal
graph, as follows.  Let $t_1$ and $t_2$ be real numbers that satisfy
$t_1/t_2 = \log \frac{q}{q-1} / \log q$ and $t_1 t_2 = m$.  Take a
complete bipartite graph between two vertex clusters $V_1$ and $V_2$
with sizes $|V_i| = \lceil t_i \rceil$, and add enough isolated
vertices to make the total number of vertices exactly $n$.  Call the
result $G_{n,m}$.

\begin{proposition}
  \label{prop:asymp-sparse}
  Fix an integer $q \geq 3$.  The following hold for all sufficiently
  large $m \leq \kappa_q n^2$.
  \begin{description}
  \item[(i)] The maximum number of $q$-colorings of an $n$-vertex
    graph with $m$ edges is $q^n e^{(-c + o(1)) \sqrt{m}}$, where $c =
    2\sqrt{\log \frac{q}{q-1} \log q}$.  Here, the $o(1)$ term tends
    to zero as $m \rightarrow \infty$.
  \item[(ii)] For any $\epsilon > 0$, as long as $m$ is sufficiently
    large, every $n$-vertex graph $G$ with $m$ edges, which maximizes
    the number of $q$-colorings, is $\epsilon m$-close to the graph
    $G_{n, m}$ which we described above.
  \end{description}
\end{proposition}

We prove this proposition in Section \ref{sec:sparse-proofs}.  Note
that part (i) is precisely the final claim of Theorem
\ref{thm:main:sparse}.

\subsection{Proof of Theorem \ref{thm:asymp-number}, part (i)}

This section contains the proofs of the claims in Section
\ref{sec:dense-structure-argument}, except for Claim 3, which is
obvious.  Together, these establish part (i) of Theorem
\ref{thm:asymp-number}, which gives the asymptotic upper bound for the
number of $q$-colorings of an $n$-vertex graph with $m$ edges.

\vspace{3mm}

\noindent \textbf{Proof of Claim 1.}\, Apply Szemer\'edi's Regularity
Lemma (Theorem \ref{thm:regularity-lemma}) with parameter $\epsilon =
\delta/3$ to partition of $V$ into nearly-equal parts $V_1$, \ldots,
$V_M$.  Then, all but $\epsilon M^2$ of the pairs $(V_i, V_j)$ are
$\epsilon$-regular, and $M \geq 1/\epsilon$.  Importantly, $M$ is also
upper bounded by a constant independent of $n$.  We clean up the graph
in a way typical of many applications of the Regularity Lemma.  Delete
all edges in each induced subgraph $G[V_i]$, all edges between pairs
$(V_i, V_j)$ which are not $\epsilon$-regular, and all edges between
pairs $(V_i, V_j)$ whose edge density is at most $\epsilon$. Since all
$|V_i| = (1+o(1))n/M$, the number of deleted edges is at most
\begin{displaymath}
(1+o(1)) \left[ 
M {n/M \choose 2} 
+ 
\epsilon M^2 (n/M)^2
+
\epsilon {n \choose 2}
\right]
\ \leq \
(1+o(1)) [
\epsilon n^2 / 2
+
\epsilon n^2
+ 
\epsilon n^2 / 2
],
\end{displaymath}
which is indeed less than $\delta n^2$ when $n$ is sufficiently large.

It remains to show property (ii).  The only edges remaining in $G'$
are those between $\epsilon$-regular pairs $(V_i, V_j)$ with
edge-density greater than $\epsilon$.  By definition of
$\epsilon$-regularity (and since $\delta > \epsilon)$, 
the edge density between
every pair of
sets $|U| \geq \delta |V_i|$, $|W| \geq \delta |V_j|$ must be
positive.  In particular, there must be at least one
edge, which establishes property (ii).  \hfill $\Box$

\vspace{3mm}

\noindent \textbf{Proof of Claim 2.}\, We aim to establish
$|\mathcal{C}_2| \geq e^{-c_\delta n} |\mathcal{C}_1|$, with $c_\delta
= q \delta \log \frac{e^2}{\delta}$.  It is a simple calculus exercise
to verify that $c_\delta \rightarrow 0$ as $\delta \rightarrow 0$.
Let us show that we can obtain any coloring $\psi \in \mathcal{C}_1$
by starting with an appropriate coloring $\psi' \in \mathcal{C}_2$,
and changing only a few color choices.  Since we may assume $\delta <
\frac{1}{q}$, every part $V_i$ has some color $c_i^*$ which appears on
at least $\delta$-fraction of its vertices.  Now consider each $V_i$.
For every color $c$ which appears less than $\delta |V_i|$ times in
$V_i$, use color $c_i^*$ to re-color all vertices of $V_i$ that had
color $c$ under $\psi$.  Now all colors appear either 0 or at least
$\delta |V_i|$ times, so once we verify that the coloring is still
proper, we will have our desired $\psi' \in \mathcal{C}_2$.  But the
only way to make a monochromatic edge is to have two distinct parts
$V_i$, $V_j$, with $c_i^* = c_j^*$, joined by at least one edge.  Then
part (ii) of Claim 1 implies that there is also some edge between the
$\delta |V_i|$ vertices in $V_i$ originally colored $c_i^*$ under
$\psi$, and the $\delta |V_j|$ vertices in $V_j$ originally colored
$c_j^*$.  This contradicts the fact that $\psi$ was a proper coloring.

Reversing the process, it is clear that $\psi$ can be recovered by
taking $\psi' \in \mathcal{C}_2$ and changing the colors of at most
$\delta |V_i|$ vertices for every color $c \in [q]$ and every $1 \leq i \leq M$.
Note that for each $c \in [q]$, we recolor a subset of $G$ of total
size at most $\sum_i \delta |V_i| = \delta n$.  Using the bounds ${n
  \choose r} \leq (en/r)^r$ and $(1+x) \leq e^x$, we see that the
total number of distinct ways in which we can modify any given $\psi' \in \mathcal{C}_2$ is at most
\begin{displaymath}
  \left[ \sum_{r = 0}^{\delta n} {n \choose r} \right]^q
  \ \leq \ \left[ (1 + \delta n) {n \choose \delta n} \right]^q
  \ \leq \ \left[ e^{\delta n} \left(\frac{en}{\delta n}\right)^{\delta n} \right]^q
  \ = \ e^{c_\delta n},
\end{displaymath}
which provides the desired upper bound on $|\mathcal{C}_1| /
|\mathcal{C}_2|$.

The final part of this claim is a simple consequence of property (ii)
of Claim 1.  Indeed, suppose that some coloring in $\mathcal{C}_2$
assigns the same color $c$ to some vertices $U_i \subset V_i$ and $U_j
\subset V_j$.  Since this is a proper coloring, there cannot be any
edges between $U_i$ and $U_j$.  Yet $|U_i| \geq \delta |V_i|$ and
$|U_j| \geq \delta |V_j|$ by definition of $\mathcal{C}_2$. Therefore, by
property (ii) of Claim 1, there are no edges at all between
$V_i$ and $V_j$, as claimed.  \hfill $\Box$

\vspace{3mm}

\noindent \textbf{Proof of Claim 4.}\, Recall that $G_\blpha(n)$ was
obtained in Construction 1 by putting a complete bipartite graph
between every pair ($V_A, V_B$) indexed by disjoint subsets.  The last
part of Claim 2 implies that $G'$ has no edges at all between parts
$V_i$ and $V_j$ which receive overlapping color sets under
$\mathcal{C}_3$.  Furthermore, each $G'[V_i]$ is empty by part (i) of
Claim 1.  So, $G'$ has no edges in each $V_A$, and also has no edges
between any $V_A$ and $V_B$ that are indexed by overlapping sets.
Hence $G'$ is indeed a subgraph of $G_\blpha(n)$.

Furthermore, $G_\blpha(n)$ has at least $m - \delta n^2$ edges,
because $G'$ differs from $G$ by at most $\delta n^2$ edges.  Yet all
$n \alpha_A$ are integers by construction, so $G_\blpha(n)$ has
precisely $\ec(\blpha) n^2$ edges.  Therefore, $\blpha \in \feas(m/n^2
- \delta)$, as claimed.  The final inequality in Claim 4 follows from
the fact that $\mathcal{C}_3$ only uses colors from $A$ to color each
$V_A$, and the definitions of $\alpha_A = |V_A|/n$ and $\obj(\blpha) =
\sum_A \alpha_A \log |A|$.  \hfill $\Box$

\vspace{3mm}

\noindent \textbf{Proof of Claim 5.}\, The only nontrivial part of
this claim is the continuity of $\opt$ on its domain, which is the set
of $\gamma$ for which $\feas(\gamma) \neq \emptyset$.  This is easily
recognized as the interval $\big(-\infty, \frac{q-1}{2q}\big]$, where
the upper endpoint, which corresponds to the $q$-partite Tur\'an
graph, equals $\ec(\blpha)$ for the vector $\blpha$ with $\alpha_A =
1/q$ for all singletons $A$.  Note that the constraint $\blpha \geq 0$
already guarantees that $\ec(\blpha) \geq 0$, so $\opt$ is constant on
$(-\infty, 0]$.

Fix an $\epsilon > 0$.  Since $\opt$ is monotonically decreasing by
definition, and constant on $(-\infty, 0]$, it suffices to show that
any $0 \leq \gamma < \gamma' \leq \frac{q-1}{2q}$ with $|\gamma' -
\gamma| < \epsilon^2$ has $\opt(\gamma') > \opt(\gamma) - 2^{q+1}
\epsilon \log q$.  Select any $\blpha$ which solves $\opt(\gamma)$.
We will adjust $\blpha$ to find an $\blpha' \in \feas(\gamma')$ with
$\obj(\blpha') > \obj(\blpha) - 2^{q+1} \epsilon \log q$, using
essentially the same perturbation as in Construction 2.

If there is an $\alpha_A \geq 2\epsilon$ with $|A| \geq 2$, shift
$\epsilon$ of $\alpha_A$'s value\footnote{Formally, $\alpha_A$ falls
  by $2\epsilon$, and each of $\alpha_{\{i\}}$ and $\alpha_{\{j\}}$
  increase by $\epsilon$.} to each of $\alpha_{\{i\}}$ and
$\alpha_{\{j\}}$ for distinct $i, j \in A$.  This clearly keeps
$\vc(\blpha)$ invariant, and it increases $\ec(\blpha)$ by at least
$\epsilon^2$ because $\alpha_{\{i\}} \alpha_{\{j\}}$ is a summand of
$\ec(\blpha)$.  Yet it only reduces $\obj(\blpha)$ by at most
$2\epsilon \log |A| \leq 2 \epsilon \log q$, so $\obj(\blpha') \geq
\obj(\blpha) - 2 \epsilon \log q$, finishing this case.

On the other hand, if all non-singletons $A$ have $\alpha_A <
2\epsilon$, then $\obj(\blpha)$ is already less than $2^q \cdot 2
\epsilon \log q$.  Since $\opt$ is always nonnegative, we trivially
have $\opt(\gamma') \geq 0 > \opt(\gamma) - 2^{q+1} \epsilon \log q$,
as desired.  \hfill $\Box$

\subsection{Proof of Theorem \ref{thm:asymp-number}, part (ii)}
\label{sec:pf-asymp-number-ii}

In this section, we establish the asymptotic tightness of our upper
bound, by showing that Construction 2 produces graphs that
asymptotically maximize the number of $q$-colorings.  We will need
Proposition \ref{prop:construction-asymp-edges}, so we prove it first.

\vspace{3mm}

\noindent \textbf{Proof of Proposition
  \ref{prop:construction-asymp-edges}.}\, Define the variables $n_A =
n \alpha_A$ (not necessarily integers), and call the expressions
$\sum_A n_A$ and $\sum_{A \cap B = \emptyset} n_A n_B$ the numbers of
\emph{fractional vertices}\/ and \emph{fractional edges},
respectively.  Initially, there are exactly $n$ fractional vertices
and $\ec(\blpha) n^2$ fractional edges.

Recall that the construction rounds each $n_A$ either up or down to
the next integer.  Let us perform these individual roundings
sequentially, finishing all of the downward roundings before the upward
roundings.  This ensures that the number of fractional vertices is
kept $\leq n$ throughout the process.  But each iteration changes the
number of fractional edges by at most $\sum_A n_A \leq n$, and there
are at most $2^q$ iterations, so our final number of edges is indeed
within $2^q n$ of $m$.

The second part of the proposition is proved similarly.  We can apply
the same iterative process to change each part size from $\alpha_A n$
to $\nu_A n$, in such a way that all downward adjustments are
performed first.  When updating the coordinate indexed by $A \subset
[q]$, we affect at most $(|\alpha_A n - \nu_A n| + 2)n$ edges,
where the extra 2 comes from the fact that the part sizes were rounded
off.  Therefore, after the $\leq 2^q$ total iterations, the total
number of edges we edit is indeed at most $\| \blpha - \bnu \|_1 n^2 +
2^{q+1} n$.  \hfill $\Box$

\vspace{3mm}

\noindent \textbf{Proof of Theorem \ref{thm:asymp-number}(ii).}\, Let
$n$ and $m$ be given, with $m$ less than the number of edges in the
Tur\'an graph $T_q(n)$.  Suppose we have a vector $\blpha \in
\feas(m/n^2)$ which achieves the maximum $\obj(\blpha) = \opt(m/n^2)$.
Construction 2 produces a graph $G_\blpha'(n)$ with $n$ vertices and
at least $m$ edges, which we will show has more than
$e^{(\opt(m/n^2)-\epsilon)n}$ proper $q$-colorings, as long as $n$ is
sufficiently large.

If $G_\blpha(n)$ already has at least $m$ edges, then we defined
$G_\blpha'(n) = G_\blpha(n)$, which has at least $\prod_A |A|^{\lfloor
  n\alpha_A \rfloor} \geq \prod_A |A|^{n \alpha_A - 1} =
e^{\obj(\blpha)n} / \prod_A |A| = e^{\obj(\blpha)n - O(1)}$ colorings,
because all colorings that use only colors from $A$ for each $V_A$ are
proper.

Otherwise, $G_\blpha(n)$ is short by, say, $k$ edges, which is $\leq
2^q n$ by Proposition \ref{prop:construction-asymp-edges}.  If the
largest $|V_A|$ indexed by a non-singleton is at least $2 \lceil
\sqrt{k} \rceil$, our construction places a $k$-edge bipartite graph
between $U_1, U_2 \subset V_A$.  Let $c_1$ and $c_2$ be two distinct
colors in $A$.  Even if we force every vertex in each $U_i$ to take
the color $c_i$, we only lose at most a factor of $q^{2 \lceil
  \sqrt{k} \rceil} = e^{o(n)}$ compared to the bound in the previous
paragraph.  This is because each of the $2 \lceil \sqrt{k} \rceil$
vertices in $U_1 \cup U_2$ had its number of color choices reduced
from $|A| \leq q$ to 1.  So, $G_\blpha'(n)$ still has at least
$e^{\obj(\blpha)n - o(n)}$ colorings.

The final case is when all parts $V_A$ indexed by non-singletons are
smaller than $2 \lceil \sqrt{k} \rceil$.  Here, the construction
simply defines $G_\blpha'(n)$ to be the Tur\'an graph $T_q(n)$.  Since
$\log |A| = 0$ for singletons $A$, the upper bound on $|V_A|$ implies
that $\obj(\blpha) \leq 2^q \cdot \frac{2 \lceil \sqrt{k} \rceil}{n}
\cdot \log q$.  This is less than $\epsilon$ for sufficiently large
$n$, because we had $k \leq 2^q n$.  Then, $e^{(\opt(m/n^2) -
  \epsilon)n} < 1$, which is of course less than the number of
$q$-colorings of the Tur\'an graph $T_q(n)$.  This completes our proof.
\hfill $\Box$

\subsection{Proof of Theorem \ref{thm:asymp-stability}}
\label{sec:asymp-stability}

In this section, we prove that any $n$-vertex graph with $m$ edges,
which maximizes the number of $q$-colorings, is in fact close (in
edit-distance) to a graph $G_\blpha(n)$ from Construction 1.  In fact,
we prove something slightly stronger: if a graph has ``close'' to the
maximum number of $q$-colorings, then it must be ``close'' (in
edit-distance) to an asymptotically optimal graph from Construction 1.

\begin{lemma}
  \label{lem:asymp-stability}
  For any $\epsilon, \kappa > 0$, there exists $\delta > 0$ such that
  the following holds for all sufficiently large $n$.  Let $G$ be an
  $n$-vertex graph with $m \leq \kappa n^2$ edges and at least
  $e^{(\opt(m/n^2) - \delta) n}$ proper $q$-colorings.  Then $G$ is
  $\epsilon n^2$-close to some $G_\blpha(n)$ from Construction 1, for
  an $\blpha$ which solves $\opt(\gamma)$ for some $|\gamma - m/n^2|
  \leq \epsilon$ with $\gamma \leq \kappa$.
\end{lemma}

Note that this lemma immediately implies Theorem
\ref{thm:asymp-stability}, because Theorem \ref{thm:asymp-number}
established that the maximum number of colorings of an $n$-vertex
graph with $m$ edges was $e^{(\opt(m/n^2) +o(1)) n}$.  Its proof is an
elementary analysis exercise in compactness, which only requires the
continuity of $\obj$, $\opt$, $\vc$, and $\ec$, the fact that $\blpha$
and the edge densities $m/n^2$ reside in compact spaces, and the
following consequence of Claims 1--4 of Section
\ref{sec:dense-structure-argument} (whose simple proof we omit):

\begin{corollary}
  \label{cor:stability-to-suboptimal}
  For every $\delta > 0$, the following holds for all sufficiently
  large $n$.  Every $q$-colorable, $n$-vertex graph $G$ with $m$ edges
  is $\delta n^2$-close to a subgraph of some $G_\blpha(n)$ with
  $\blpha \in \feas(m/n^2 - \delta)$.  Also, $G$ has at most
  $e^{(\obj(\blpha) + \delta) n}$ proper $q$-colorings.
\end{corollary}

% \noindent \textbf{Proof.}\, By stringing together Claims 1--4, we see
% that for any fixed $\delta > 0$, the following holds for all
% sufficiently large $n$.  For every $q$-colorable graph $G$ with $n$
% vertices and $m$ edges we can find an $\blpha \in \feas(m/n^2 -
% \delta)$ such that $\#\{\text{colorings of $G$}\} \leq e^{c_\delta n}
% \cdot 2^{qM} \cdot e^{\obj(\blpha)n}$, where $M$ was bounded for every
% fixed $\delta$, and the constant $c_\delta \rightarrow 0$ when $\delta
% \rightarrow 0$.  Therefore, given any $\epsilon > 0$, it is possible
% to choose a $\delta > 0$ such that $\obj(\blpha) > \log
% \#\{\text{colorings of $G$}\} - \epsilon$ whenever $n$ is
% sufficiently large. \hfill $\Box$

\vspace{3mm}

\noindent \textbf{Proof of Lemma \ref{lem:asymp-stability}.}\, We
proceed by contradiction.  Then, there is some fixed $\epsilon > 0$, a
sequence $\delta_i \rightarrow 0$, and a sequence of graphs $G_i$ with
the following properties.
\begin{description}
\item[(i)] $G_i$ has at least as many vertices as required to apply
  Corollary \ref{cor:stability-to-suboptimal} with parameter
  $\delta_i$.
\item[(ii)] $G_i$ has at least $e^{(\opt(m_i/n_i^2) - \delta_i) n_i}$
  colorings, where $n_i$ and $m_i$ are its numbers of vertices and
  edges, and $m_i \leq \kappa n_i^2$.
\item[(iii)] $G_i$ is at least $\epsilon n_i^2$-far from
  $G_\blpha(n_i)$ for every $\blpha$ that solves $\opt(\gamma)$ with
  $|\gamma - m_i/n_i^2| \leq \epsilon$.
\end{description}
Applying Corollary \ref{cor:stability-to-suboptimal} to each $G_i$
with parameter $\delta_i$, we find vectors $\blpha_i
\in \feas(m_i/n_i^2 - \delta_i)$ such that $G_i$ is $\delta_i
n_i^2$-close to some subgraph $G_i'$ of $G_{\blpha_i}(n_i)$, and each
$G_i$ has at most $e^{(\obj(\blpha_i)+\delta_i)n_i}$ proper
$q$-colorings.  Combining this with property (ii) above, we find that
each $\obj(\blpha_i) \geq \opt(m_i/n_i^2) - 2\delta_i$.  The densities
$m_i/n_i^2$ and the vectors $\blpha_i$ live in bounded (hence compact)
spaces.  So, by passing to a subsequence, we may assume that
$m_i/n_i^2 \rightarrow \gamma \leq \kappa$ and $\blpha_i \rightarrow
\blpha$ for some limit points $\gamma$ and $\blpha$.

Observe that by continuity, both $\blpha \in \feas(\gamma)$ and $\obj(\blpha) \geq \opt(\gamma)$.
Therefore $\blpha$ solves $\opt(\gamma)$, i.e., $\obj(\blpha)=\opt(\gamma)$.
Furthermore, although \emph{a priori}\/ we only knew
that $\ec(\blpha) \geq \gamma$, maximality implies that in fact
$\ec(\blpha) = \gamma$. Indeed, if not then one could shift more mass to
$\alpha_{[q]}$ to increase $\obj(\blpha)$ while staying within the
feasible set. This would contradict that $\obj(\blpha)=\opt(\gamma)$.

We finish by showing that eventually $G_i$ is $\epsilon n_i^2$-close
to $G_\blpha(n_i)$, contradicting (iii).  To do this, we show that all
three of the edit-distances between $G_i \leftrightarrow G_i'
\leftrightarrow G_{\blpha_i}(n_i) \leftrightarrow G_\blpha(n_i)$ are
$o(n_i^2)$.  The closeness of the first pair follows by construction
since $\delta_i \rightarrow 0$, and the closeness of the last pair
follows from Proposition \ref{prop:construction-asymp-edges}
because $\blpha_i \rightarrow \blpha$.

For the central pair, recall that $G_i'$ is actually contained in
$G_{\blpha_i}(n_i)$, so we only need to compare their numbers of
edges.  In fact, since we already established $o(n_i^2)$-closeness of
the first and last pairs, it suffices to show that the difference
between the number of edges in $G_i$ and $G_\blpha(n_i)$ is
$o(n_i^2)$.  Recall from above that $\ec(\blpha) = \gamma$, and
therefore by Proposition \ref{prop:construction-asymp-edges},
$G_\blpha(n_i)$ has $\ec(\blpha)n_i^2 + o(n_i^2) = (\gamma + o(1))
n_i^2$ edges.  Yet $G_i$ also has $(\gamma+o(1)) n_i^2$ edges, because
$m_i/n_i^2 \rightarrow \gamma$.  This completes the proof. \hfill
$\Box$

\subsection{Proofs for the sparse case}
\label{sec:sparse-proofs}

In this section, we prove the statements which refine our results in
the case when the graph is sparse, i.e., $m = o(n^2)$.  We begin with
the lemma which shows that every sparse graph with the maximum number
of colorings has a dense core which spans all of the edges.

\vspace{3mm}

\noindent \textbf{Proof of Lemma \ref{lem:sparse-has-core}.}\, Let
$n_1$ be the number of non-isolated vertices in $G$, and let $r$ be
the number of connected components in the subgraph induced by the
non-isolated vertices.  Since all such vertices there have degree at
least 1, we have $r \leq n_1/2$.

Any connected graph on $t$ vertices has at most $q (q-1)^{t-1}$ proper
$q$-colorings, because we may iteratively color the vertices along a
depth-first-search tree rooted at an arbitrary vertex; when we visit
any vertex other than the root, there will only be at most $q-1$
colors left to choose from.  So, $G$ has at most $q^{n-n_1} \cdot q^r
\cdot (q-1)^{n_1-r}$ colorings, where the first factor comes from the
fact that isolated vertices have a free choice over all $q$ colors.
Using $r \leq n_1/2$, this bound is at most $q^{n-n_1/2}
(q-1)^{n_1/2}$.

But since $G$ is optimal, it must have at least as many colorings as
the Tur\'an graph $T_q(n_2)$ plus $n-n_2$ isolated vertices, where
$n_2 = \Theta(\sqrt{m})$ is the minimum number of vertices in a
$q$-partite Tur\'an graph with at least $m$ edges.  The isolated
vertices already give the latter graph at least $q^{n-n_2}$ colorings,
so we must have $q^{n-n_2} \leq q^{n-n_1/2} (q-1)^{n_1/2}$, which
implies that
\begin{equation}
  \label{eq:n_1}
  n_1 \leq n_2 \cdot (2 \log q) / \left(\log \frac{q}{q-1}\right) .
\end{equation}
The expression on the right hand side is $\Theta(n_2) =
\Theta(\sqrt{m})$, so if we define the integer $n_0$ to be the maximum
of right hand side in \eqref{eq:n_1} and $\sqrt{m/\kappa}$ (rounding up
to the next integer if necessary) then we indeed have $n_1 \leq n_0 =
\Theta(n_2) = \Theta(\sqrt{m})$.  \hfill $\Box$

% \vspace{3mm}
% 
% \noindent \textbf{Remark.}\, Observe that for any fixed $q$, the ratio
% $m/n_0^2$ tends to a limit as $m \rightarrow \infty$.  This is
% because the right hand side of equation \eqref{eq:n_1} is a constant
% multiple of $n_2$, which in turn was the minimum number of vertices in
% a $q$-partite Tur\'an graph with at least $m$ edges, i.e., $n_2 =
% (1+o(1)) \sqrt{m \cdot \frac{2q}{q-1} }$.

\vspace{3mm}

Next, we prove the first part of Proposition \ref{prop:asymp-sparse},
which claims that the maximum number of $q$-colorings of an $n$-vertex
graph with $m \leq \kappa_q n^2$ edges is asymptotically $q^n e^{(-c +
  o(1)) \sqrt{m}}$, where $\kappa_q = \left( \sqrt{\frac{\log
      q/(q-1)}{\log q}} + \sqrt{\frac{\log q}{\log q/(q-1)}}
\right)^{-2}$ and $c = 2\sqrt{\log \frac{q}{q-1} \log q}$.

\vspace{3mm}

\noindent \textbf{Proof of Proposition \ref{prop:asymp-sparse}(i).}\,
Let $G$ be an $n$-vertex graph with $m$ edges, which maximizes the
number of $q$-colorings.  Let $n_0$ be the integer obtained by
applying Lemma \ref{lem:sparse-has-core} with threshold $\kappa_q$.
If $n \geq n_0$, the lemma gives a dense $n_0$-vertex subgraph $G'
\subset G$ which contains all of the edges.  Otherwise, set $G' = G$.
In either case, we obtain a graph $G'$ whose number of vertices $n'$
is $\Theta(\sqrt{m})$, and $m/(n')^2 \leq \kappa_q$.

Since the vertices in $G \setminus G'$ (if any) are isolated, the
number of $q$-colorings of $G$ is precisely $q^{n-n'}$ times the
number of $q$-colorings of $G'$.  Therefore, $G'$ must also have the
maximum number of $q$-colorings over all $n'$-vertex graphs with $m$
edges.  Applying Theorem \ref{thm:asymp-number} to $G'$, we find that
$G'$ has $e^{(\opt(m/(n')^2) + o(1))n'}$ colorings.  Proposition
\ref{prop:solve-opt-sparse} gives us the precise answer
$\opt(m/(n')^2) = \log q - 2\sqrt{\frac{m}{(n')^2} \cdot \log
  \frac{q}{q-1} \cdot \log q}$, so substituting that in gives us that
the number of $q$-colorings of $G$ is:
\begin{displaymath}
  q^{n-n'} \cdot e^{(\opt(m/(n')^2) + o(1))n'} 
  \ = \
  q^{n-n'} \cdot q^{n'} e^{(-c+o(1))\sqrt{m}}
  \ = \
  q^n e^{(-c+o(1))\sqrt{m}},
\end{displaymath}
where $c$ is indeed the same constant as claimed in the statement of
this proposition.  \hfill $\Box$

\vspace{3mm}

We finish this section by proving the stability result which shows
that any optimal sparse graph is $\epsilon m$-close (in edit-distance)
to the graph $G_{n, m}$ defined in Section \ref{sec:sparse-summary}.

\vspace{3mm}

\noindent \textbf{Proof of Proposition \ref{prop:asymp-sparse}(ii).}\,
Let $G$ be an $n$-vertex graph with $m$ edges, which maximizes the
number of $q$-colorings.  We will actually show the equivalent
statement that $G$ is $O((\epsilon + \sqrt{\epsilon})m)$-close to
$G_{n, m}$.

As in the proof of part (i) above, we find a dense $n'$-vertex
subgraph $G' \subset G$ that spans all of the edges, which itself must
maximize the number of $q$-colorings.  Using the same parameters as
above, we have $n' = \Theta(\sqrt{m})$ and $m \leq \kappa_q (n')^2$.
By Theorem \ref{thm:asymp-stability}, $G'$ must be $\epsilon
(n')^2$-close to a graph $G_\blpha(n')$ from Construction 1, for some
$\blpha$ that solves $\opt(\gamma)$ with $\gamma \leq \kappa_q$.
Since $n' = \Theta(\sqrt{m})$, the graphs are $O(\epsilon m)$-close.
The $\gamma$ is within the range in which Proposition
\ref{prop:solve-opt-sparse} solved Optimization Problem 1, so
$G_\blpha(n')$ is a complete bipartite graph plus isolated vertices,
which indeed resembles $G_{n, m}$.

Moreover, the ratio between the sizes of the sides of the complete
bipartite graph in $G_\blpha(n')$ is correct, because it tends to the
constant $\log \frac{q}{q-1} / \log q$ regardless of the value of
$\gamma$.  Also, their product, which equals the number of edges in
$G_\blpha(n')$, is within $O(\epsilon m)$ of $m$ because
$G_\blpha(n')$ is $O(\epsilon m)$-close to the $m$-edge graph $G'$.
Therefore, each of the sides of the complete bipartite graph in
$G_\blpha(n')$ differs in size from its corresponding side in $G_{n,
  m}$ by at most $O(\sqrt{\epsilon m})$.  Since each side of the
bipartite graph in $G_{n, m}$ has size $\Theta(\sqrt{m})$, we can
transform $G_\blpha(n')$ into $G_{n, m}$ by adding isolated vertices
and editing at most $O(\sqrt{\epsilon} \cdot m)$ edges.  Yet by
construction of $\blpha$, the graphs $G'$ and $G_\blpha(n')$ were
$O(\epsilon m)$-close, modulo isolated vertices.  Therefore, $G$ and
$G_{n, m}$ are indeed $O((\epsilon + \sqrt{\epsilon})m)$-close, as
claimed.  \hfill $\Box$

\section{Solving the optimization problem}
\label{sec:solve-opt}

In this section, we solve the optimization problem for low densities,
for all values of $q$.  We also solve it for all densities in the case
when $q=3$.

\subsection{Sparse case}

\label{sec:solve-opt-sparse}

The key observation is that when the edge density is low, we can
reduce the optimization problem to one with no edge density parameter
and no vertex constraint.  This turns out to be substantially easier
to solve.

\vspace{3mm}

\noindent \textbf{Optimization Problem 2.}\, Fix an integer $q$, and
consider the following objective and constraint functions:
\begin{displaymath}
  \objS(\blpha) := \sum_A \alpha_A \log \frac{|A|}{q}\,;
  \quad\quad\quad
  \ec(\blpha) := \sum_{A \cap B = \emptyset} \alpha_A \alpha_B.
\end{displaymath}
The vector $\blpha$ has $2^q - 2$ coordinates $\alpha_A \in
\mathbb{R}$ indexed by the nonempty \textbf{proper} subsets $A \subset
[q]$, and the sum in $\ec(\blpha)$ runs over unordered pairs of
disjoint sets $\{A,B\}$.  Let $\feasS$ be the feasible set of vectors
defined by the constraints $\blpha \geq 0$ and $\ec(\blpha) \geq 1$.
We seek to maximize $\objS(\blpha)$ over the set $\feasS$, and we
define $\optS$ to be this maximum value, which we will show to exist
in Section \ref{sec:sparse:observ}.  We write that the vector $\blpha$
\emph{solves}\/ $\optS$ when both $\blpha \in \feasS$ and
$\objS(\blpha) = \optS$.

\begin{proposition}
  \label{prop:solve-opt-2}
  For any given $q \geq 3$, the \textbf{unique} solution (up to a permutation
  of the base set $[q]$) to Optimization Problem 2 is the vector
  $\blpha^*$ with
  \begin{displaymath}
    \alpha_{\{1\}}^* = \sqrt{\log \frac{q}{q-1} \, / \, \log q},
    \quad\quad\quad
    \alpha_{\{2, \ldots q\}}^* = \frac{1}{\alpha_{\{1\}}^*},
    \quad\quad\quad
    \text{and all other $\alpha_A^* = 0$.}
  \end{displaymath}
  This gives $\objS(\blpha^*) = -2\sqrt{\log \frac{q}{q-1} \log q}$.
\end{proposition}

Let us show how Proposition \ref{prop:solve-opt-2} implies Proposition
\ref{prop:solve-opt-sparse}, which gave the solution to Optimization
Problem 1 for sufficiently low edge densities $\gamma$.

\vspace{3mm}

\noindent \textbf{Proof of Proposition \ref{prop:solve-opt-sparse}.}\,
Let $\blpha^*$ be the unique maximizer for Optimization Problem 2, and
consider any number $t \geq \vc(\blpha^*)$.  Then $\blpha^*$ is still
the unique maximizer of $\objS(\blpha)$ when $\blpha$ is required to
satisfy the vacuous condition $\vc(\blpha) \leq t$ as well.  Let
$\overline{\blpha}$ be the vector obtained by dividing every entry of
$\blpha^*$ by $t$, and adding a new entry $\overline{\alpha}_{[q]}$ so
that $\vc(\overline{\blpha}) = 1$.

Then, $\overline{\blpha}$ is the unique maximizer of $\objS(\blpha)$
when $\blpha$ is constrained by $\vc(\blpha) = 1$ and $\ec(\blpha)
\geq t^{-2}$.  But when $\vc(\blpha) = 1$ is one of the constraints,
then $\objS(\blpha) = \obj(\blpha) - \log q$, so this implies that
$\overline{\blpha}$ is the unique solution to $\opt(t^{-2})$.  Using
the substitution $\gamma = t^{-2}$, we see that $\overline{\blpha}$ is
precisely the vector described in \eqref{eq:sparse-opt-soln}.  Since
$t \geq \vc(\blpha^*)$ was arbitrary, we conclude that this holds for
all $\gamma$ below $\vc(\blpha^*)^{-2} = \left( \sqrt{\frac{\log
      q/(q-1)}{\log q}} + \sqrt{\frac{\log q}{\log q/(q-1)}}
\right)^{-2} = \kappa_q$. \hfill $\Box$

\subsubsection{Observations for Optimization Problem 2}

\label{sec:sparse:observ}

We begin by showing that $\objS$ attains its maximum on the feasible
set $\feasS$.  Since $\feasS$ is clearly nonempty, there is some
finite $c \in \mathbb{R}$ for which $\optS \geq c$.  In the formula
for $\objS$, all coefficients $\log \frac{|A|}{q}$ of the $\alpha_A$
are negative, so we only need to consider the compact region bounded
by $0 \leq \alpha_A \leq c/\log \frac{|A|}{q}$ for each $A$.
Therefore, by compactness, $\objS$ indeed attains its maximum on
$\feasS$.

Now that we know the maximum is attained, we can use perturbation
arguments to determine its location.  The following definition will be
convenient for our analysis.

\begin{definition}
  Let the \textbf{support} of a vector $\blpha$ be the collection of
  $A$ for which $\alpha_A \neq 0$.
\end{definition}

The following lemma will allow us to reduce to the case of considering
optimal vectors whose supports are a partition of $[q]$.

\begin{lemma}
  \label{lem:sparse:support=partition}
  One of the vectors $\blpha$ which solves $\optS$ has support that is
  a partition\footnote{A collection of disjoint sets whose union
    is $[q]$.} of $[q]$.  Furthermore, if the only partitions that
  support optimal vectors consist of a singleton plus a $(q-1)$-set,
  then in fact every vector which solves $\optS$ is supported by such
  a partition.
\end{lemma}

\noindent \textbf{Proof.}\, We begin with the first statement.  Let
$\blpha$ be a vector which solves $\optS$, and suppose that its
support contains two intersecting sets $A$ and $B$.  We will perturb
$\alpha_A$ and $\alpha_B$ while keeping all other $\alpha$'s fixed.
Since $A$ and $B$ intersect, the polynomial $\ec(\blpha)$ has no
products $\alpha_A \alpha_B$, i.e., it is of the form $x \alpha_A + y
\alpha_B + z$, for some constants $x, y, z \geq 0$.

Furthermore, $x \neq 0$, or else we could reduce $\alpha_A$ to zero
without affecting $\ec(\blpha)$, but this would strictly increase
$\objS(\blpha)$ because all coefficients $\log \frac{|A|}{q}$ in
$\objS$ are negative.  Similarly, $y \neq 0$.  Therefore, we may
perturb $\alpha_A$ by $+ty$ and $\alpha_B$ by $-tx$, while keeping
$\ec(\blpha)$ fixed.  Since we may use both positive and negative $t$
and $\objS$ itself is linear in $\alpha_A$ and $\alpha_B$, optimality
implies that $\objS$ does not depend on $t$.  Hence we may choose a
$t$ which drives one of $\alpha_A$ or $\alpha_B$ to zero (we are free to pick which one), and $\objS$
will remain unchanged.

Repeating this process, we eventually obtain a vector $\blpha$ which
is supported by disjoint sets.  Their union must be the entire $[q]$,
because otherwise we could simply grow one of the sets in the support
by adding the unused elements of $[q]$.  This would not affect
$\ec(\blpha)$, but it would strictly increase $\objS$.

It remains to prove the second part of our lemma.  Let $\blpha$ be an
optimal vector, and apply the above reduction process to simplify its
support.  At the end, we will have a vector supported by $|A| = 1$ and
$|B| = q-1$, by assumption.  Each iteration of the reduction removes
exactly one set from the support, so the second to last stage will have
some $\blpha'$ supported by three distinct sets, two of which are the
final $A$ and $B$, and the third which we call $C$.

In the reduction, when we consider two overlapping sets, we are free
to select which one is removed.  Therefore, we could choose to keep
the third set $C$ and remove one of $A$ and $B$, and then continue
reducing until the support is disjoint, while keeping $\objS$
unchanged.  Yet no matter what $C$ was, it is impossible for this
alternative reduction route to terminate in a partition of $[q]$,
contradicting the above observation that any reduction must terminate
in a partition.  \hfill $\Box$

\begin{definition}
  \label{def:sparse:IJ}
  Let $\blpha$ be a fixed vector whose support is a partition of
  $[q]$.  For each $A \subset [q]$, define the expressions:
  \begin{displaymath}
    I_A \ = \ \alpha_A \sum_{B \neq A} \alpha_B
    \quad\quad\quad
    J_A \ = \ \frac{1}{\objS(\blpha)} \cdot \alpha_A \log \frac{|A|}{q}.
  \end{displaymath}
\end{definition}

% Note that these are normalized so that $\sum_A I_A = 2$ and $\sum_A
% J_A = 1$.

\begin{lemma}
  \label{lem:sparse:IJ}
  Let $\blpha$ be a vector which solves $\optS$, whose support is a
  partition of $[q]$.  Then:
  \begin{description}
  \item[(i)] For every $A \subset [q]$, we have $I_A = 2J_A$.  In
    particular, for each $A$ in the support, $I_A / \alpha_A = 2J_A /
    \alpha_A$.
  \item[(ii)] Suppose $A$ and $B$ are both in the support, and $|A| =
    |B|$.  Then $\alpha_A = \alpha_B$ as well.
  \end{description}
\end{lemma}

\noindent \textbf{Proof.}\, We begin with part (i).  Fix any $A
\subset [q]$.  Consider the following operation for small $\epsilon >
0$.  First, replace $\alpha_A$ by $(1+\epsilon)\alpha_A$.
Observe that $I_A = \alpha_A \sum_{B: B \cap A
  = \emptyset} \alpha_B$ because the support of $\blpha$ is a
partition of $[q]$.  Therefore we increase $\ec(\blpha) = \sum_{A \cap B=\emptyset} \alpha_A \alpha_B$ by
$\epsilon I_A$. Next,
multiply all $\alpha$'s (including the one we just increased) by $(1 +
\epsilon I_A)^{-1/2}$. Then $\ec(\blpha)$ is still at least 1 and 
our perturbed vector is in $\feasS$. Its new objective equals $\objS(\blpha) \cdot \frac{1 + \epsilon
  J_A}{\sqrt{1 + \epsilon I_A}}$.  Since $\blpha$ maximized the
objective (which is always negative), we must have $\frac{1 + \epsilon
  J_A}{\sqrt{1 + \epsilon I_A}} \geq 1$.  Rearranging, this implies
that $I_A \leq 2J_A + \epsilon J_A^2$.  Sending $\epsilon \rightarrow
0$, we see that $I_A \leq 2J_A$.  The opposite inequality follows from
considering the replacement of $\alpha_A$ by $(1-\epsilon)\alpha_A$,
and then multiplying $\alpha$'s by $(1 - \epsilon I_A)^{-1/2}$.  This
establishes part (i).

% Part (ii) is obvious because $\alpha_A \neq 0$ for $A$ in the support.

For part (ii), let $S = \sum_C \alpha_C$.  Since the support of
$\blpha$ is a partition of $[q]$, $S - \alpha_A = I_A/\alpha_A$.  By
part (i), this equals $2J_A/\alpha_A = \log \frac{|A|}{q} /
\objS(\blpha)$, which is determined by the cardinality of $A$.
Therefore, $S - \alpha_A = S - \alpha_B$, which implies (ii).  \hfill
$\Box$

\subsubsection{Solution to Optimization Problem 2 for $\boldsymbol{q < 9}$}

In its original form, Optimization Problem 2 involves exponentially
many variables, but Lemma \ref{lem:sparse:support=partition}
dramatically reduces their number by allowing us to consider only
supports that are partitions of $[q]$.  Therefore, we need to make one
computation per partition of $[q]$, which can actually be done
\emph{symbolically}\/ (hence exactly) by \emph{Mathematica}.  The
running time of \emph{Mathematica}'s symbolic maximization is
double-exponential in the number of variables, so it was particularly
helpful to reduce the number of variables.  The entire computation for
$q \in \{3, \ldots, 8\}$ took less than an hour, and the complete
\emph{Mathematica}\/ program and output appear in Appendix
\ref{sec:mathematica}.

Let us illustrate this process by showing what needs to be done for
the partition $7 = 2+2+3$.  This corresponds to maximizing $\alpha_A
\log \frac{2}{7} + \alpha_B \log \frac{2}{7} + \alpha_C \log
\frac{3}{7}$ subject to the constraints $\alpha_A \alpha_B + \alpha_B
\alpha_C + \alpha_C \alpha_A \geq 1$ and $\blpha \geq 0$.  By Lemma
\ref{lem:sparse:IJ}(ii), we may assume $\alpha_A = \alpha_B$, so it
suffices to maximize $2 x \log \frac{2}{7} + y \log \frac{3}{7}$
subject to $x^2 + 2xy \geq 1$ and $x, y \geq 0$.  This is achieved by
\emph{Mathematica}'s \texttt{Maximize} function:
\begin{verbatim}
Maximize[{2 x Log[2/7] + y Log[3/7], x^2 + 2 x y >= 1 && x >= 0 && y >= 0}, {x, y}]
\end{verbatim}
\emph{Mathematica}\/ answers that the maximum value is
$-\sqrt{-\big(\log\frac{7}{3}\big)^2 + 4 \log \frac{7}{3} \log
  \frac{7}{2}} \approx -1.9$, which is indeed less than the claimed
value $-2 \sqrt{\log \frac{7}{7-1} \log 7} \approx -1.1$.

We performed one such computation per partition of each $q \in \{3,
\ldots, 8\}$.  In every case except for the partition $q = 1 + (q-1)$,
the maximum indeed fell short of the claimed value.  That final
partition is completely solved analytically (i.e., including the
uniqueness result) by Lemma \ref{lem:sparse:2sets=>done} in the next
section.  This completes the analysis for all $q < 9$.

\subsubsection{Solution to Optimization Problem 2 for $\boldsymbol{q \geq 9}$}

\label{sec:solve-opt-sparse-q>=9}

We begin by ruling out several extreme partitions that our general
argument below will not handle.  As one may expect, each of these
special cases has a fairly pedestrian proof, so we postpone the
proofs of the following two lemmas to the appendix.

\begin{lemma}
  \label{lem:sparse:2sets=>done}
  Fix any integer $q \geq 3$, and let $\blpha$ be a vector which
  solves $\optS$.  If the support of $\blpha$ is a partition of $[q]$
  into exactly two sets, then (up to permutation of the ground set
  $[q]$) $\blpha$ must be equal to the claimed unique optimal vector
  $\blpha^*$ in Proposition \ref{prop:solve-opt-2}.
\end{lemma}

\begin{lemma}
  \label{lem:sparse:extreme}
  Fix any integer $q \geq 4$, and let $\blpha$ be a vector which
  solves $\optS$, whose support is a partition of $[q]$.  Then that
  partition cannot have any of the following forms:
  \begin{description}
  \item[(i)] all singletons;  % actually holds for q>=3
  \item[(ii)] all singletons, except for one 2-set;
  \item[(iii)] have a $(q-2)$-set as one of the parts. % actually also for q>=3
  \end{description}
\end{lemma}

The heart of the solution to the optimization problem is the following
general case, which we will prove momentarily.

\begin{lemma}
  \label{lem:sparse:<q-2=>done}
  Fix any integer $q \geq 9$, and let $\blpha$ be a vector which
  solves $\optS$, whose support is a partition of $[q]$.  Then that
  partition must have a set of size at least $q-2$.
\end{lemma}

These collected results show that $\optS$ has the unique solution that
we claimed at the beginning of this section.

\vspace{3mm}

\noindent \textbf{Proof of Proposition \ref{prop:solve-opt-2} for
  $\boldsymbol{q \geq 9}$.}\, Let $\blpha$ be a vector which solves
$\optS$.  By Lemma \ref{lem:sparse:support=partition}, we may assume
that its support is a partition of $[q]$.  It cannot be a single set
(of cardinality $q$), because then $\ec(\blpha) = 0$, and by Lemmas
\ref{lem:sparse:extreme}(iii) and \ref{lem:sparse:<q-2=>done}, the
support cannot contain a set of size $\leq q-2$.  

Thus, the support must contain a set of size $q-1$, and since it is a
partition, the only other set is a singleton.  Then Lemma
\ref{lem:sparse:2sets=>done} gives us that $\blpha$ equals the claimed
unique optimal vector $\blpha^*$, up to a permutation of the ground
set $[q]$.  This completes the proof.  \hfill $\Box$

\vspace{3mm}

In the remainder of this section, we prove the general case (Lemma
\ref{lem:sparse:<q-2=>done}).  The following definition and fact are
convenient, but the proof is a routine calculus exercise, so we
postpone it to the appendix.

\begin{lemma}
  \label{lem:sparse:Fq}
  Define the function $F_q(x) = \log \frac{q}{q-x} \cdot \log
  \frac{q}{x}$.
  \begin{description}
  \item[(i)] For $q > 0$, $F_q(x)$ strictly increases on $0 < x <
    q/2$ and strictly decreases on $q/2 < x < q$.
  \item[(ii)] For $q \geq 9$, we have the inequality $F_q(3) > 2F_q(1)
    \cdot \frac{q-3}{q-2}$.
  \end{description}
\end{lemma}

\vspace{3mm}

\noindent \textbf{Proof of Lemma \ref{lem:sparse:<q-2=>done}.}\,
Assume for the sake of contradiction that all sets in the support of
the optimal $\blpha$ have size at most $q-3$.  In terms of the
expressions $I$ and $J$ from Definition \ref{def:sparse:IJ}, we have
the following equality, where the sums should be interpreted as only
over sets in the support of $\blpha$:
\begin{displaymath}
  \frac{2 \log \frac{|A|}{q} }{\objS(\blpha)} 
  \ \ = \ \ \frac{2 J_A}{\alpha_A}
  \ \ = \ \ \frac{I_A}{\alpha_A}
  \ \ = \ \ \sum_{B \neq A} \alpha_B
  \ \ = \ \ \sum_{B \neq A} \frac{J_B \cdot \objS(\blpha)}{\log \frac{|B|}{q}}.
\end{displaymath}
(The second equality is Lemma \ref{lem:sparse:IJ}(i), and the other
three equalities come from the definitions of $I$ and $J$.)  Note
that the above logarithms are always negative.  It is cleaner to work
with positive quantities, so we rewrite the above equality in the
equivalent form:
\begin{displaymath}
  \frac{2 \log \frac{q}{|A|} }{\objS(\blpha)} = \sum_{B \neq A} \frac{J_B \cdot \objS(\blpha)}{\log \frac{q}{|B|}}.
\end{displaymath}
Since every $B$ in the above sum is disjoint from $A$ and we assumed
all sets in the support have size at most $q-3$, we have that every
$B$ above has size $|B|\leq q - \max\{|A|, 3\}$.  This gives the upper
bound:
\begin{eqnarray*}
  \frac{2 \log \frac{q}{|A|} }{\objS(\blpha)} 
  &\leq& \sum_{B \neq A} \frac{J_B \cdot \objS(\blpha)}{\log \frac{q}{q-\max\{|A|, 3\}}} \\
  \frac{ 2 \cdot \log \frac{q}{|A|} \cdot \log \frac{q}{q-\max\{|A|, 3\}} }{\objS(\blpha)^2} &\leq& \sum_{B \neq A} J_B.
\end{eqnarray*}
Since $|A| \leq \max\{|A|, 3\}$, the left hand side is at least
$2F_q(\max\{|A|, 3\}) / \objS(\blpha)^2$.  Also, $F_q(x)$ is symmetric
about $x=q/2$ and we assumed that $3 \leq q/2$ and $|A| \leq q-3$, so
Lemma \ref{lem:sparse:Fq}(i) implies that this is in turn $\geq
2F_q(3) / \objS(\blpha)^2$.  Lemma \ref{lem:sparse:Fq}(ii) bounds this
in terms of $F_q(1)$, which ultimately gives us the following
bound for $\sum_{B \neq A} J_B$:
\begin{equation}
  \label{ineq:sparse:Jsum}
  \frac{q-3}{q-2} 
  \ \ \leq \ \ \frac{ \objS(\blpha^*)^2 }{ \objS(\blpha)^2 } \cdot \frac{q-3}{q-2}
  \ \ = \ \ \frac{ 4 F_q(1) }{ \objS(\blpha)^2 } \cdot \frac{q-3}{q-2}
  \ \ < \ \ \frac{2F_q(3)}{\objS(\blpha)^2}
  \ \ \leq \ \ \sum_{B \neq A} J_B.
\end{equation}
Here, $\blpha^*$ is the claimed optimal vector in Proposition
\ref{prop:solve-opt-2}, and we recognize $4F_q(1) =
\objS(\blpha^*)^2$.  The first inequality follows from the maximality
of $\blpha$, and its direction is reversed because $\objS$ is always
negative.

Let $t$ be the number of sets in the support of $\blpha$.  Summing
\eqref{ineq:sparse:Jsum} over all sets $A$ in the support:
\begin{displaymath}
  t \cdot \frac{q-3}{q-2} \ \ < \ \ \sum_A \sum_{B \neq A} J_B \ \ = \ \ \sum_B J_B (t-1).
\end{displaymath}
Yet $\sum_B J_B = 1$ by definition, so this implies $\frac{t}{t-1} <
\frac{q-2}{q-3}$, which forces $t > q-2$.  Then, the support must be
all singletons, except possibly for a single 2-set.  This contradicts
Lemma \ref{lem:sparse:extreme}, and completes our proof.  \hfill
$\Box$

\subsection{Solving the optimization problem for 3 colors}

In this section, we provide the complete analytic solution to
Optimization Problem 1, for the entire range of the edge density
parameter $\gamma$ when the number of colors $q$ is exactly 3.  To
simplify notation, we will write $\alpha_{12}$ instead of
$\alpha_{\{1, 2\}}$, etc.

\begin{proposition}
  \label{prop:solve-opt}
  Define the constant $c = \left( \sqrt{\frac{\log 3/2}{\log 3}} +
    \sqrt{\frac{\log 3}{\log 3/2}} \right)^{-2} \approx 0.1969$.
  Then, the \textbf{unique} solution (up to a permutation of the index
  set $\{1,2,3\}$) of Optimization Problem 1 with edge density
  parameter $\gamma$ is the vector $\blpha$ defined as follows.
  (All unspecified $\alpha_A$ below are zero.)
  \begin{description}
  \item[(i)] If $0 \leq \gamma \leq c$, then $\alpha_3 = \sqrt{\gamma
      \cdot \frac{\log 3/2}{\log 3}}$, $\alpha_{12} =
    \frac{\gamma}{\alpha_3}$, and $\alpha_{123} = 1-\alpha_{12} -
    \alpha_3$.  This gives $\opt(\gamma) = \log 3 - 2\sqrt{\gamma
      \cdot \log 3 \cdot \log \frac{3}{2}}$.

  \item[(ii)] If $c \leq \gamma \leq \frac{1}{4}$, then $\alpha_{12} =
    \frac{1+\sqrt{1-4\gamma}}{2}$ and $\alpha_3 = 1-\alpha_{12}$,
    which gives $\opt(\gamma) = \frac{1+\sqrt{1-4\gamma}}{2} \cdot
    \log 2$.

  \item[(iii)] If $\frac{1}{4} \leq \gamma \leq \frac{1}{3}$, then
    $\alpha_{12} = \frac{1 - \sqrt{12\gamma - 3}}{2}$, $\alpha_1 =
    \alpha_2 = \frac{1-2\alpha_{12}}{3}$, and $\alpha_3 =
    \frac{1+\alpha_{12}}{3}$, which gives $\opt(\gamma) = \frac{1 -
      \sqrt{12\gamma - 3}}{2} \cdot \log 2$.
  \end{description}
\end{proposition}

This covers the entire range of admissible $\gamma$, because $\gamma =
1/3$ corresponds to the density of the Tur\'an graph $T_3(n)$, which
is the densest 3-colorable graph.

\subsubsection{Outline of solution}

The strategy of the solution is as follows.  Suppose we have some
$\blpha$ that solves $\opt(\gamma)$.  Since we may permute the index
set, we may assume without loss of generality that $\alpha_1 \leq
\alpha_2 \leq \alpha_3$.  We then use perturbation arguments to
pinpoint the location of $\blpha$.  Although the problem initially
looks cumbersome (there are 7 nontrivially-related variables), the
solution cleanly follows from 6 short steps.

\begin{description}
\item[Step 1.] By \emph{shifting mass}\footnote{Adjusting the values
    of the $\alpha_A$ while conserving their sum $\sum_A \alpha_A =
    \vc(\blpha)$.} between the $\alpha_A$ with $|A| = 2$, we deduce
  that $\alpha_{23}$ and $\alpha_{13}$ are both zero.

\item[Step 2.] By smoothing together $\alpha_1$ and $\alpha_2$, we
  deduce that $\alpha_1 = \alpha_2$.

\item[Step 3.] By shifting mass between the variables $\alpha_A$ with
  $|A| = 1$, we reduce to one of the following two situations.  Either
  $\alpha_1 = \alpha_2 = 0$, or $0 < \alpha_1 = \alpha_2 = \alpha_3 -
  \alpha_{12}$.

\item[Step 4.] We solve the first case resulting from Step 3, which is
  vastly simpler than the original problem.  We find that the solution
  corresponds to outcomes (i) and (ii) of Proposition
  \ref{prop:solve-opt}.

\item[Step 5.] It remains to consider the second case resulting from
  Step 3.  By taking mass away from both $\alpha_{123}$ and
  $\alpha_1$, and giving it to $\alpha_{12}$, we conclude that
  $\alpha_{123} = 0$.

\item[Step 6.] We are left with the situation where the only nonzero
  variables are $\alpha_1$, $\alpha_2$, $\alpha_3$, and $\alpha_{12}$,
  and they are related by the equation $\alpha_1 = \alpha_2 = \alpha_3
  - \alpha_{12}$.  Again, this is vastly simpler than the original
  problem, and we find that its solution corresponds to outcome (iii)
  of Proposition \ref{prop:solve-opt}.
\end{description}

\subsubsection{Details of solution}

% The perturbations that we perform in this section all consist of
% \emph{shifting mass}\/ between the coordinates of the vector $\blpha$.
% This means that we adjust the values of the coordinates while
% conserving the total sum $\sum_A \alpha_A = \vc(\blpha)$.  
We begin by recording a simple result that we will use repeatedly in
the solution.

\begin{lemma}
  \label{lem:shift-equal-index}
  Let $\blpha$ be a vector that solves $\opt(\gamma)$.  Then
  $\ec(\blpha) = \gamma$.  Furthermore, if $\blpha'$ is obtained from
  $\blpha$ by shifting mass from some $\alpha_A$ to another $\alpha_B$
  with $|A| = |B|$, then $\ec(\blpha') \leq \ec(\blpha)$.
\end{lemma}

\noindent \textbf{Proof.}\, Suppose for contradiction that
$\ec(\blpha) > \gamma$.  The slack in the edge constraint lets us
shift some more mass to $\alpha_{123}$ while keeping $\ec(\blpha) \geq
\gamma$.  But in the definition of $\obj$, the coefficient ($\log 3$)
of $\alpha_{123}$ is the largest, so this shift strictly increases
$\obj$, contradicting maximality of $\blpha$.

For the second claim, observe that $\obj$ is invariant under the shift
since $|A| = |B|$.  Now suppose for contradiction that $\ec(\blpha') >
\ec(\blpha)$.  Then, as above, we could shift more mass to
$\alpha_{123}$, which would strictly increase $\obj$, again
contradicting the maximality of $\blpha$.  \hfill $\Box$

\vspace{3mm}

\noindent \textbf{Step 1.}\, Consider shifting mass among
$\{\alpha_{12}, \alpha_{23}, \alpha_{13}\}$.  If we hold all other
$\alpha_A$ constant, then $\ec(\blpha) = \alpha_1 \alpha_{23} +
\alpha_2 \alpha_{13} + \alpha_3 \alpha_{12} + \text{constant}$, which
is linear in the three variables of interest.

Let us postpone the uniqueness claim for a moment.  Since we ordered
$\alpha_1 \leq \alpha_2 \leq \alpha_3$, shifting all of the mass from
$\{\alpha_{13}, \alpha_{23}\}$ to $\alpha_{12}$ will either strictly
grow $\ec(\blpha)$ if $\alpha_2 < \alpha_3$, or keep $\ec(\blpha)$
unchanged.  Also, $\obj(\blpha)$ will be invariant.  Therefore, if we
are only looking for an upper bound for $\opt(\gamma)$, we may perform
this shift, and reduce to the case when $\alpha_{13} = 0 =
\alpha_{23}$ without loss of generality.

We return to the topic of uniqueness.  The next five steps of this
solution will deduce that, conditioned on $\alpha_{13} = 0 =
\alpha_{23}$, the unique optimal $\blpha$ always has either $\alpha_2
< \alpha_3$ or $\alpha_{12} = \alpha_{13} = \alpha_{23} = 0$.  We
claim that this implies that our initial shift of mass to
$\alpha_{12}$ \emph{never happened}.  Indeed, in the case with
$\alpha_2 < \alpha_3$, the previous paragraph shows that an initial
shift would have strictly increased $\ec(\blpha)$, violating Lemma
\ref{lem:shift-equal-index}.  And in the case with $\alpha_{12} =
\alpha_{13} = \alpha_{23} = 0$, there was not even any mass at all to
shift.  Therefore, this will imply the full uniqueness result.

\vspace{3mm}

\noindent \textbf{Step 2.}\, Consider shifting mass between $\alpha_1$
and $\alpha_2$ until they become equal.  If we hold all other
$\alpha_A$ constant, then $\ec(\blpha) = \alpha_1 \alpha_2 + (\alpha_1
+ \alpha_2) \alpha_3 + \text{constant}$.  This ``smoothing'' operation
strictly increases the first term, while keeping the other terms
invariant.  But Lemma \ref{lem:shift-equal-index} prohibits
$\ec(\blpha)$ from increasing, so we conclude that we must have had
$\alpha_1 = \alpha_2$.

\vspace{3mm}

\noindent \textbf{Step 3.}\, Consider shifting mass among $\{\alpha_1,
\alpha_2, \alpha_3\}$.  That is, fix $S = \alpha_1 + \alpha_2 +
\alpha_3$, and vary $t = \alpha_3$ in the range $0 \leq t \leq S$.  By
Step 2, $\alpha_1 = \alpha_2 = \frac{S-t}{2}$.  Step 1 gave $\alpha_{13} =
\alpha_{23} = 0$, so we have:
\begin{eqnarray*}
  \ec(\blpha) 
  \ \ = \ \ \alpha_1 \alpha_2 + \alpha_1 \alpha_3 + \alpha_2 \alpha_3 + \alpha_{12} \alpha_3 
  &=& \frac{(S-t)^2}{4} + 2 \cdot \frac{S-t}{2} \cdot t + \alpha_{12} t \\
  &=& -\frac{3}{4} t^2 + \left(\frac{S}{2} + \alpha_{12}\right) t + \frac{S^2}{4}.
\end{eqnarray*}
By Lemma \ref{lem:shift-equal-index}, $\alpha_3 = t$ must maximize
this downward-opening parabola in the range $0 \leq t \leq S$.  Recall
that quadratics $f(x) = ax^2 + bx + c$ reach their extreme value at $x
= -\frac{b}{2a}$, which corresponds to $t = -\big(\frac{S}{2} +
\alpha_{12}\big)/\big(2 \cdot \big(-\frac{3}{4}\big)\big) = \frac{S +
  2\alpha_{12}}{3}$ above.  Thus, if $\frac{S + 2\alpha_{12}}{3} < S$,
then we must have $\alpha_3 = \frac{S + 2\alpha_{12}}{3} =
\frac{\alpha_1 + \alpha_2 + \alpha_3 + 2\alpha_{12}}{3}$.  Step 2 gave
us $\alpha_1 = \alpha_2$, which forces $0 < \alpha_1 = \alpha_2 =
\alpha_3 - \alpha_{12}$.  This is the second claimed outcome of this
step.

On the other hand, if $\frac{S + 2\alpha_{12}}{3} \geq S$, then the
quadratic is strictly increasing on the interval $0 \leq t \leq S$.
Therefore, we must have $\alpha_3 = S$, forcing $\alpha_1 = \alpha_2 =
0$.  This is the first claimed outcome of this step.

\vspace{3mm}

\noindent \textbf{Step 4.}\, In this case, only $\alpha_3$,
$\alpha_{12}$, and $\alpha_{123}$ are nonzero.  Then the edge
constraint is simply $\ec(\blpha) = \alpha_3 \alpha_{12} = \gamma$
(Lemma \ref{lem:shift-equal-index} forces equality).  Note that since
$\alpha_3 + \alpha_{12} \leq \vc(\blpha) = 1$, their product $\alpha_3
\alpha_{12}$ is always at most $1/4$, \textbf{so we can only be in
  this case when} $\boldsymbol{\gamma \leq 1/4}$.

Now let $x = \alpha_3$ and $y = \alpha_{12}$.  The vertex constraint
forces $\alpha_{123} = 1-x-y$, so we are left with the routine problem
of maximizing $\obj = y \log 2 + (1-x-y) \log 3 = \log 3 - x \log 3 -
y \log \frac{3}{2}$ subject to the constraints
\begin{displaymath}
  x,y \geq 0,
  \quad \quad
  x+y \leq 1,
  \quad \quad
  xy = \gamma.
\end{displaymath}
These constraints specify a segment of a hyperbola (a convex function)
in the first quadrant of the $xy$-plane, and the objective is linear
in $x$ and $y$.  Therefore, by convexity, the maximum would be at the
global maximum of $\obj$ on the entire first quadrant branch of the
hyperbola, unless that fell outside the segment, in which case it must
be at an endpoint, forcing $x+y=1$.

The maximum over the entire branch of $xy = \gamma$ follows easily
from the inequality of arithmetic and geometric means: $\obj \leq \log
3 - 2\sqrt{x\log 3 \cdot y \log \frac{3}{2}} = \log 3 - 2\sqrt{\gamma
  \cdot \log 3 \cdot \log \frac{3}{2}}$, with equality when $x \log 3
= y \log \frac{3}{2}$.  Using $xy = \gamma$ to solve for $x$ and $y$,
we see that the unique global maximum is at $x = \sqrt{\gamma \cdot
  \frac{\log 3/2}{\log 3}}$ and $y = \sqrt{\gamma \cdot \frac{\log
    3}{\log 3/2}}$.  This lies on our segment (satisfies $x+y \leq
1$) precisely when $\gamma$ is below the constant $c \approx 0.1969$
in Proposition \ref{prop:solve-opt}, and these values of $\alpha_3 =
x$ and $\alpha_{12} = y$ indeed match those claimed in that regime.

On the other hand, when $\gamma > c$, we are outside the segment, so
by the above we must have $x+y = 1$, and we may substitute $x = 1-y$.
We are left with the single-variable maximization of $\obj = y \log 2$
subject to $0 \leq y \leq 1$ and $(1-y)y = \gamma$.  By the quadratic
formula, this is at $\alpha_{12} = y = \frac{1+\sqrt{1-4\gamma}}{2}
\leq 1$, which produces $\alpha_3 = x = 1-y = 1-\alpha_{12}$.  This
indeed matches outcome (ii) of our proposition.

\vspace{3mm}

\noindent \textbf{Step 5.}\, The remaining case is $0 < \alpha_1 =
\alpha_2 = \alpha_3 - \alpha_{12}$, and we will show that this forces
$\alpha_{123} = 0$.  Indeed, suppose for the sake of contradiction
that $\alpha_{123} > 0$.  Shift mass to $\alpha_{12}$ by taking
$\epsilon$ from $\alpha_{123}$ and $\epsilon' = \epsilon \alpha_3 /
\alpha_2$ from $\alpha_1$.  Since many $\alpha_A$ are zero,
$\ec(\blpha) = \alpha_1(\alpha_2 + \alpha_3) + \alpha_2 \alpha_3 +
\alpha_{12} \alpha_3 $.  Our perturbation decreases the first term by
$\epsilon' (\alpha_2 + \alpha_3)$, increases the third term by
$(\epsilon + \epsilon')\alpha_3$, and does not change the second term,
so our choice of $\epsilon'$ keeps $\ec(\blpha)$ invariant.

On the other hand, $\obj$ increases by $(\epsilon + \epsilon') \log 2
-\epsilon \log 3$.  Since we know $\alpha_2 = \alpha_3 - \alpha_{12}$,
in particular we always have $\alpha_3 \geq \alpha_2$, which implies
that $\epsilon' \geq \epsilon$ because we assume $\alpha_2,\alpha_3 > 0$.  Hence the increase in $\obj$ is
$(\epsilon + \epsilon') \log 2 - \epsilon \log 3
  \geq
  (\epsilon + \epsilon) \log 2 - \epsilon \log 3
  > 0$,
contradicting the maximality of $\blpha$.  Therefore, we must have had
$\alpha_{123} = 0$.

\vspace{3mm}

\noindent \textbf{Step 6.}\, Now only $\alpha_1$, $\alpha_2$,
$\alpha_3$, and $\alpha_{12}$ remain.  Let $t = \alpha_3$ and $r =
\alpha_{12}$.  Step 3 gives $\alpha_1 = \alpha_2 = \alpha_3 -
\alpha_{12} = t-r$.  We use the vertex constraint to eliminate $t$: $1
= \vc(\blpha) = 2(t-r) + t + r$, so $t = \frac{1+r}{3}$.  Substituting
this for $t$, we are left with $\alpha_1 = \alpha_2 = \frac{1-2r}{3}$
and $\alpha_3 = \frac{1+r}{3}$.  Since we need all $\alpha_A \geq 0$,
the range for $r$ is $0 \leq r \leq 1/2$.  

The above expressions give $\ec(\blpha) = \left(\frac{1 -
    2r}{3}\right)^2 + 2\left(\frac{1-2r}{3}\right)
\left(\frac{1+r}{3}\right) + \left(\frac{1+r}{3}\right) r = \frac{r^2
  - r + 1}{3}$, and Lemma \ref{lem:shift-equal-index} forces
$\ec(\blpha) = \gamma$.  The quadratic formula gives the roots $r =
\frac{1 \pm \sqrt{12\gamma - 3}}{2}$.  These are only real when
$12\gamma - 3 \geq 0$, so \textbf{this case only occurs when}
$\boldsymbol{\gamma \geq 1/4}$.  Furthermore, the only root within the
interval $0 \leq r \leq 1/2$ is $r = \frac{1 - \sqrt{12\gamma -
    3}}{2}$.  Plugging this value of $r$ into the expressions for the
$\alpha_A$, we indeed obtain outcome (iii) of Proposition
\ref{prop:solve-opt}.

\vspace{3mm}

\noindent \textbf{Conclusion.}\, The only steps which proposed
possible maxima were Steps 4 and 6.  Conveniently, Step 4 also
required that $\gamma \leq 1/4$, while Step 6 required $\gamma \geq
1/4$ (both deductions are bolded above), so we do not need to compare
them except at $\gamma = 1/4$, which is trivial.  Finally, note that
all extremal outcomes indeed have $\alpha_2 < \alpha_3$, except at
$\gamma = 1/3$, in which case $\alpha_{12} = \alpha_{13} = \alpha_{23}
= 0$.  This justifies the uniqueness argument that we used at the end
of Step 1, and completes our proof of Proposition
\ref{prop:solve-opt}. \hfill $\Box$

\section{Exact result for sparse graphs}
\label{sec:exact:sparse}

In this section, we determine the precise structure of the sparse
graphs that maximize the number of colorings, completing the proof of
Theorem \ref{thm:main:sparse}.  Proposition
\ref{prop:asymp-sparse}(ii) showed that in this regime, the optimal
graphs were close, in edit distance, to complete bipartite graphs.  As
a warm-up for the arguments that will follow in this section, let us
begin by showing that the semi-complete subgraphs of Definition
\ref{def:semi-complete} are optimal among bipartite graphs.  We will
use this in the final stage of our proof of the exact result.

\begin{lemma}
  \label{lem:semi-complete}
  Let $q \geq 3$ and $r < a \leq b$ be positive integers.  Among all
  subgraphs of $K_{a,b}$ with $r$ missing edges, the ones which
  maximize the number of $q$-colorings are precisely:
  \begin{description}
  \item[(i)] both the correctly and incorrectly oriented semi-complete
    subgraphs, when $q=3$, and
  \item[(ii)] the correctly oriented semi-complete subgraph, when $q
    \geq 4$ and $\frac{b}{a} \geq \log q / \log \frac{q-2}{q-3}$ and
    $a$ is sufficiently large (i.e., $a > N_q$, where $N_q$ depends
    only on $q$).
  \end{description}
\end{lemma}

\noindent \textbf{Remark.}\, The above result is not as clean when
more than 3 colors are used, but is sufficient for our purposes.  In
the sparse case, we encounter only highly unbalanced bipartite graphs,
all of which have part size ratio approximately $\log q / \log
\frac{q}{q-1}$.  Apparently out of sheer coincidence (and good
fortune), this is just barely enough to satisfy the additional
condition of the lemma.  Nevertheless, it would be nice to remove that
condition.

\vspace{3mm}

\noindent \textbf{Proof of Lemma \ref{lem:semi-complete}(ii).}\, Let
$A \cup B$ be the vertex partition of $K_{a,b}$, with $|A| = a$ and
$|B| = b$.  Let $F^*$ be the correctly oriented semi-complete subgraph
of $K_{a,b}$ with exactly $r$ missing edges.  Let $F$ be another
non-isomorphic subgraph of $K_{a,b}$ with the same number of edges.
We will show that $F$ has fewer colorings.  Since $F$ and $F^*$ are
both bipartite, they share every coloring that uses disjoint sets of
colors on the sides of the bipartition.  Discrepancies arise when the
same color appears on both sides.  Note, however, that whenever this
occurs, every edge between same-colored vertices must be missing from
the graph.  This set of forced missing edges,\footnote{In this lemma,
  \emph{missing edges}\/ refer only to those missing from the
  bipartite $K_{a,b}$, not the entire $K_{a+b}$.}  which we call the
coloring's \emph{footprint}, is always a union of vertex-disjoint
complete bipartite graphs, one per color that appears on both sides.
For each subset $H$ of the missing edges of $F$, let $n_H$ be the
number of colorings of $F$ with footprint $H$.  Then, $\sum n_H$ is
exactly the number of colorings of $F$.  To give each $n_H$ a
counterpart from $F^*$, fix an arbitrary bijection $\phi$ between the
missing edges of $F$ and $F^*$, and let $n_H^*$ be the number of
colorings of $F^*$ with footprint $\phi(H)$.  Since $F^*$ has $\sum
n_H^*$ colorings, it suffices to show that $n_H \leq n_H^*$ for all
$H$, with strict inequality for at least one $H$.

Clearly, when $H$ is empty, or a star centered in $B$, then $n_H =
n_H^*$.  We observed that all footprints are unions $\Gamma_1 \cup
\cdots \cup \Gamma_k$ of vertex-disjoint complete bipartite graphs, so
all $H$ not of that form automatically have $n_H = 0 \leq n_H^*$.  It
remains to consider $H$ that have this form, but are not stars
centered in $B$.  Colorings with this footprint are monochromatic on
each $\Gamma_i$, and there are ${q \choose k} k!$ ways to choose a
distinct color for each $\Gamma_i$.  The remaining $q-k$ colors are
partitioned into two sets, one for $A \setminus V(H)$ and one for $B
\setminus V(H)$.  Crucially, $|B \setminus V(H)| \leq b-2$ because $H$ is
not a star centered in $B$.  Thus,
\begin{eqnarray*}
  n_H &\leq& \left[ {q \choose k} k! \right] \cdot 
  \sum_{i=1}^{q-k-1} {q-k \choose i} i^{|A \setminus V(H)|}
  (q-k-i)^{|B \setminus V(H)|}  \\
  &\leq&
  q^k \cdot \sum_{i=1}^{q-k-1} {q-k \choose i} i^a (q-k-i)^{b-2}.
\end{eqnarray*}
To see that the sum is dominated by the $i=1$ term, note that since we
assumed that $\frac{b}{a} \geq \log q / \log \frac{q-2}{q-3}$, for
sufficiently large $a$ we have 
\begin{displaymath}
  \frac{b-2}{a} \geq \log (q-1) / \log
  \frac{q-2}{q-3} \geq \log (q-k) / \log \frac{q-k-1}{q-k-2},
\end{displaymath}
so we may apply Inequality \ref{ineq:partition-colors}(ii) from the Appendix.  This
gives $n_H \leq q^k \cdot 1.1 (q-k) (q-k-1)^{b-2}$.  Next, we claim
that this bound is greatest when $k$ is smallest.  Indeed, when $k$
increases by one, $q^k$ increases by the factor $q$, but
$(q-k-1)^{b-2}$ decreases by a factor of at least $\big(
\frac{q-2}{q-3} \big)^{b-2} \gg q$ for large $b$.  Hence we have $n_H
\leq 1.1 q (q-1) (q-2)^{b-2}$.

On the other hand, $\phi(H)$ is always a star centered in $B$, so we
can easily construct $q(q-1)(q-2)^{b-1}$ colorings of $F^*$.  Indeed,
choose one color for the vertices of the graph $\phi(H)$, a different
color for the remainder of $A \setminus \phi(H)$, and allow each
vertex left in $B \setminus \phi(H)$ to take any of the other $q-2$
colors.  Since $\phi(H)$ intersects $B$ in exactly one vertex, $n_H^*
\geq q(q-1)(q-2)^{b-1}$, as claimed.  But $q-2 \geq 2$, so we have the
desired strict inequality $n_H^* \geq 2 q(q-1) (q-2)^{b-2} > n_H$ for
all remaining $H$. \hfill $\Box$

\vspace{3mm}

Part (i) is a consequence of the following more precise result, which
we will also need later.

\begin{lemma}
  \label{lem:subgraph-bipartite:q=3}
  Let $F$ be a subgraph of the complete bipartite graph $K_{a,b}$ with
  vertex partition $A \cup B$, and $r < \max\{a, b\}$ missing
  edges.  Suppose $F$ has $x \in A$ and $y \in B$ with $x$ complete to
  $B$ and $y$ complete to $A$.  Then its number of 3-colorings is
  precisely $3 \cdot 2^a + 3 \cdot 2^b - 6 + 6s$, where $s$ is
  the number of nonempty subsets of missing edges which form complete
  bipartite graphs.  This is at most $3 \cdot 2^a + 3 \cdot 2^b +
  6 \cdot (2^r-2)$, with equality exactly when the missing edges form a star.
\end{lemma}

\noindent \textbf{Proof.}\, As in the proof of Lemma
\ref{lem:semi-complete}(ii), let $n_H$ be the number of 3-colorings of
$F$ with footprint $H$.  The key observation is that for every
nonempty $H$, $n_H = 6$ when $H$ is a complete bipartite graph, and
$n_H = 0$ otherwise.  Indeed, if $H$ is not a complete bipartite
graph, then it cannot be a footprint of a 3-coloring, so $n_H = 0$.
Otherwise, there are 3 ways to choose a color for the vertices of $H$,
and then by definition of footprint, the remaining two colors must be
split between $A \setminus H$ and $B \setminus H$.  Both of these sets
are nonempty, because $A \setminus H$ must contain the given vertex
$x$ and $B \setminus H$ must contain $y$, so the only way to split the
two colors is to use one on all of $A \setminus H$ and the other on
all of $B \setminus H$.  There are 2 ways to decide how to do this.
So, $n_H = 3 \cdot 2 = 6$, as claimed, and this produces the $6s$ in
the formula.

The rest of the formula follows from $n_\emptyset = 3 \cdot 2^a +
3 \cdot 2^b - 6$.  Indeed, the terms correspond to the colorings
that use a single color (for which there are three choices) on $B$ and
allow the other two on $A$, those that use one on $A$ and allow the
others on $B$, and those that use only one on each of $A$ and $B$
(hence were double-counted).  The final claim in the statement comes
from the fact that stars are the only $r$-edge graphs which have all
$2^r-1$ of their nonempty subgraphs complete bipartite. \hfill $\Box$

\vspace{3mm}

\noindent \textbf{Proof of Lemma \ref{lem:semi-complete}(i).}\, Since
the number of missing edges $r$ is less than both $|A|$ and $|B|$, the
vertices $x$ and $y$ of Lemma \ref{lem:subgraph-bipartite:q=3} must
exist.  Therefore, its equality condition implies that the optimal
subgraphs are indeed semi-complete. \hfill $\Box$

\subsection{Structure of proof}

We will use several small constants with relative order of magnitude
$\epsilon_1 \ll \epsilon_2 \ll \epsilon_3$, related by $\epsilon_1 =
\epsilon_2^2 = \epsilon_3^3$.  We do not send them to zero; rather, we
show that there is an eventual choice of the $\epsilon_i$, determined
by $q$ and $\kappa$, that makes our argument work.  So, to avoid
confusion, the $O$, $\Theta$, and $o$ notation that we employ in this
proof will only mask constants depending on $q,\kappa$ alone.  For
example, we will write $X = O(\epsilon_2 Y)$ when there is a constant
$C_{q,\kappa}$ such that $X \leq C_{q,\kappa} \epsilon_2 Y$ for
sufficiently large $m$ and $n$.  Occasionally, we will use phrases
like ``almost all colorings have property $P$'' when
$(1-o(1))$-fraction of all colorings have that property.

\vspace{3mm}

\noindent \textbf{Proof of Theorem \ref{thm:main:sparse}.}\, Let $G =
(V, E)$ be an optimal graph with $n$ vertices and $m \leq \kappa n^2$
edges.  We begin with a convenient technical modification: if $G$ has
an isolated edge $xy$, replace it with an edge between $x$ and another
non-isolated vertex of minimal degree.  Do this only once, even if $G$
had multiple isolated edges.  The number of colorings stays the same
because both graphs share the same partial colorings of $V \setminus
\{x\}$, and each of those has exactly $q-1$ extensions (in each graph)
to the degree-1 vertex $x$.

This adjustment will not compromise the uniqueness claim, because it
cannot create one of the optimal graphs listed in Theorem
\ref{thm:main:sparse}.  Indeed, if it did, then the degree-1 vertex
$x$ would now have to be the center of the missing star of the
semi-complete subgraph $H \subset K_{a,b}$.  But we made $x$
adjacent to a vertex of minimal degree, so $x$ must be on the smaller
side of $H$'s bipartition.  Then the number of $K_{a,b}$-edges missing
from the semi-complete $H$ is precisely $b-d(x) = b-1$.  This exceeds
$a$ for all optimal graphs listed in Theorem \ref{thm:main:sparse},
but our definition of semi-completeness required that the number of
missing edges was strictly less than the size of the smaller part.
This contradiction shows that we may assume without loss of generality
that if $G$ has an isolated edge $uv$, then it also contains a
degree-1 vertex $x \not \in \{u, v\}$.

Define $u_1 = \sqrt{m \cdot \log \frac{q}{q-1} / \log q}$ and $u_2 =
\sqrt{m \cdot \log q / \log \frac{q}{q-1}}$, and note that
$\frac{u_1}{u_2} = \log \frac{q}{q-1} / \log q$ and $u_1 u_2 = m$.
So, Proposition \ref{prop:asymp-sparse}(ii) gives disjoint subsets
$U_1, U_2 \subset V$ of size $|U_i| = \lceil u_i \rceil$, such that by
editing at most $\epsilon_1 m$ edges, we can transform $G$ into the
complete bipartite graph between $U_1$ and $U_2$, with all other
vertices isolated.  Call that graph $G^*$.

Let $(V_1, V_2)$ be a max-cut partition of the \textbf{non-isolated}
vertices of $G$, such that $V_1$ contains at least as many vertices of
$U_1$ as $V_2$ does.  We would like to show that this partition is
very close to $(U_1, U_2)$, so we keep track of the $U_i$ by defining
$U_i' = U_i \cap V_i$ and $U_i'' = U_i \cap V_{3-i}$ for each $i \in
\{1, 2\}$.  To help us recognize vertices that are ``mostly correct,''
let $X_i \subset U_i'$ be the vertices that are adjacent to all but at
most $\epsilon_2\sqrt{m}$ vertices of $U_{3-i}'$.

The following series of claims will complete the proof of Theorem
\ref{thm:main:sparse}, since Proposition \ref{prop:asymp-sparse}(i)
already determined the asymptotic maximum number of colorings.

\begin{description}

\item[Claim 1.] For each $i$, $|U_i'|$ is within $O(\epsilon_1
  \sqrt{m})$ of $u_i$, $|X_i|$ is within $O(\epsilon_2 \sqrt{m})$ of
  $u_i$, and $|U_i''| \leq O(\epsilon_1 \sqrt{m})$.

\item[Claim 2.] Almost all colorings of $G$ are \emph{$(X_1,
    X_2)$-regular}, which means that they only use one color on $X_1$,
  and avoid that color on $X_2$.

\item[Claim 3.] At most one non-isolated vertex $v_0$ has degree $\leq
  2 \epsilon_3 \sqrt{m}$.  We use this to show that each $|V_i|$ is
  within $O(\epsilon_2 \sqrt{m})$ of $u_i$.  Let $V_0 = \{v_0\}$ if it
  exists; otherwise, let $V_0 = \emptyset$.  Let $V_i^* = V_i
  \setminus V_0$.

\item[Claim 4.]  Almost all colorings are \emph{$(V_1^*, V_2^*)$-regular},
  i.e., use one color for $V_1^*$, and avoid it on $V_2^*$.

\item[Claim 5.] Each $V_i^*$ is an independent set, and $v_0$ (if it
  exists) has neighbors in only one of the $V_i^*$.  Hence $G$ is a
  bipartite graph plus isolated vertices.

\item[Claim 6.] $G$ is a semi-complete subgraph of $K_{|V_1|, |V_2|}$
  plus isolated vertices, correctly oriented if $q \geq 4$.

\end{description}

\subsection{Details of proof}
\label{sec:exact:sparse:details}

\noindent \textbf{Proof of Claim 1.}\, We know that by editing at most
$\epsilon_1 m$ edges, $G$ can be transformed into $G^*$, the complete
bipartite graph between $(U_1, U_2)$, plus isolated vertices.  Since
$|U_i| = \lceil u_i \rceil = \Theta(\sqrt{m})$, all vertices in the
$U_i$ have degree $\Theta(\sqrt{m})$ in $G^*$.  So, the number of
$U_i$-vertices that are isolated in $G$ is at most $\frac{\epsilon_1
  m}{\Theta(\sqrt{m})} = O(\epsilon_1 \sqrt{m})$, implying in
particular that the number of $U_1$-vertices in $V_1 \cup V_2$ is at
least $|U_1| - O(\epsilon_1 \sqrt{m}) \geq \frac{2}{3} u_1$.  (Recall
that $(V_1, V_2)$ is a max-cut partition of the \emph{non-isolated}
vertices of $G$.)  Since more $U_1$-vertices are in $V_1$ than in
$V_2$, and $U_1' = U_1 \cap V_1$, we have $|U_1'| \geq \frac{1}{3} u_1
= \Theta(\sqrt{m})$.

Also, $G^*$ has at least $m$ edges crossing between $(U_1, U_2)$, so
$G$ has at least $m-\epsilon_1 m$ edges crossing between $(U_1, U_2)$,
and at least that many between its max-cut $(V_1, V_2)$.  As $G$ has
only $m$ edges, this shows that each $G[V_i]$ spans at most
$\epsilon_1 m$ edges.  But the sets $U_1', U_2'' \subset V_1$ are
complete to each other in $G^*$, so among the $\leq \epsilon_1 m$
edges of $G[V_1]$, at least $|U_1'| |U_2''| - \epsilon_1 m$ of them
must go between $U_1'$ and $U_2''$.  Combining this with the above
result that $|U_1'| \geq \Theta(\sqrt{m})$, we obtain the desired
bound $|U_2''| \leq O(\epsilon_1 \sqrt{m})$.

Then $U_2'$, the set of $U_2$-vertices in $V_2$, has size at least
$u_2 - O(\epsilon_1 \sqrt{m}) \geq \Theta(\sqrt{m})$, because only
$O(\epsilon_1 \sqrt{m})$ of the $U_2$-vertices are isolated and
$|U_2''| \leq O(\epsilon_1 \sqrt{m})$ of them are in $V_1$.  Repeating
the previous paragraph's argument with respect to $U_2'$ and $U_1''$,
we find that $|U_1''| \leq O(\epsilon_1 \sqrt{m})$, which then implies
that $|U_1'| \geq u_1 - O(\epsilon_1 \sqrt{m})$.

It remains to control $X_i$, which we recall to be the vertices of
$U_i'$ which had at most $\epsilon_2\sqrt{m}$ non-neighbors in
$U_{3-i}'$.  The $U_i'$ are complete to each other in $G^*$, so each
vertex not in $X_i$ contributes at least $\epsilon_2\sqrt{m}$ to the
total edit distance of $\leq \epsilon_1 m$.  We set $\epsilon_2^2 =
\epsilon_1$, so this implies that all but at most $\epsilon_2
\sqrt{m}$ vertices of $U_i'$ belong to $X_i$.  Since $|U_i'|$ is
within $O(\epsilon_1 \sqrt{m})$ of $u_i$, this gives the desired
result.  \hfill $\Box$

\vspace{3mm}

\noindent \textbf{Proof of Claim 2.}\, We bound the number of
colorings that are not $(X_1, X_2)$-regular.  For each partition $[q]
= C_0 \cup C_1 \cup C_2 \cup C_3$, we count the colorings which use
the colors $C_1$ in $X_1$ but not $X_2$, use $C_2$ in $X_2$ but not
$X_1$, use $C_3$ in both $X_1$ and $X_2$, and do not use $C_0$ in
either $X_1$ or $X_2$.  Then we sum over all \emph{irregular}\/
partitions, which are all partitions with $|C_1| \geq 2$ or 
$|C_3| \geq 1$.  It suffices to
show that the result is of smaller order than the total number of
colorings of $G$.

For any given partition with $|C_i| = c_i$, we claim that the
corresponding number of colorings is at most $(|X_1| |X_2|)^{c_3}
\cdot c_1^{|X_1| - q \epsilon_2 \sqrt{m}} \cdot c_2^{|X_2| - q
  \epsilon_2 \sqrt{m}} \cdot q^{n - 2c_3 - (|X_1| - q \epsilon_2
  \sqrt{m}) - (|X_2| - q \epsilon_2 \sqrt{m})}$.  The first factor
comes from choosing $c_3$ pairs of vertices $x_i \in X_1$, $y_i \in
X_2$ on which to use each color of $C_3$.  Then, every vertex in the
common neighborhood of $\{y_i\}$ must avoid $C_3$ in order to produce
a proper coloring.  By definition of $X_2$, the number of vertices of
$U_1'$ that are not in this common neighborhood is at most $|C_3|
\epsilon_2 \sqrt{m} \leq q \epsilon_2 \sqrt{m}$.  Thus all but at most
$q \epsilon_2 \sqrt{m}$ vertices of $X_1 \subset U_1'$ are adjacent to
every $\{y_i\}$, and therefore restricted to colors in $C_1$.  This
produces the second factor in our bound, and the third factor is
obtained analogously.  Of course every vertex has at most $q$ color
choices, and we use that trivial bound for all remaining vertices,
producing our final factor.  Using that each $|X_i|$ is within
$O(\epsilon_2 \sqrt{m})$ of $u_i = \Theta(\sqrt{m})$, we find that the
sum $\Sigma_1$ of this bound over all $\leq 4^q$ irregular partitions
is:
\begin{eqnarray*}
  \Sigma_1 &=& \sum_{\text{irregular}} (|X_1| |X_2|)^{c_3} \cdot c_1^{|X_1|
    - q \epsilon_2 \sqrt{m}} \cdot c_2^{|X_2| - q \epsilon_2 \sqrt{m}}
  \cdot q^{n - 2c_3 - (|X_1| - q \epsilon_2 \sqrt{m}) - (|X_2| - q \epsilon_2
    \sqrt{m})} \\
  &\leq& e^{O(\epsilon_2 \sqrt{m})} \sum_{\text{irregular}} (\Theta(\sqrt{m}) \cdot \Theta(\sqrt{m}))^{c_3} \cdot c_1^{u_1} \cdot c_2^{u_2}
  \cdot q^{n - u_1 - u_2} \\
  &\leq& e^{O(\epsilon_2 \sqrt{m})} \cdot 4^q \cdot O(m^q) \cdot 
  \max_{c_1 \geq 2 \text{ or } c_3 \geq 1} \left\{ c_1^{u_1} c_2^{u_2} \right\} \cdot
  q^{n - u_1 - u_2}.
\end{eqnarray*}
For any irregular partition with $c_1 + c_2 < q$, it is clear that
$c_1^{u_1} c_2^{u_2}$ increases when $C_1$ is replaced by $C_1 \cup
C_0 \cup C_3$, and $C_0$ and $C_3$ are reduced to $\emptyset$.  It is
also clear that this procedure gives another irregular partition, but
this time with $c_1 + c_2 = q$.  
Yet $\frac{u_2}{u_1} = \log q /
\log \frac{q}{q-1} \geq \log q / \log \frac{q-1}{q-2}$, so we may apply
Inequality \ref{ineq:partition-colors}(i), which gives
\begin{displaymath}
  \max_{c_1 \geq 2 \text{ or } c_3 \geq 1} c_1^{u_1} c_2^{u_2}
  \ \ = \ \ 
  2^{u_1} (q-2)^{u_2}
  \ \ \leq \ \ 
  1.5^{-u_1} \cdot 1^{u_1} (q-1)^{u_2}
  \ \ = \ \ 
  e^{-\Theta(\sqrt{m})} \cdot (q-1)^{u_2}.
\end{displaymath}
Thus for small $\epsilon_2$, we have $\Sigma_1 \leq
e^{-\Theta(\sqrt{m})} \cdot (q-1)^{u_2} \cdot q^{n - u_1 - u_2}$.

On the other hand, Proposition \ref{prop:asymp-sparse}(i) shows that
the optimal graph has at least $\Sigma_0 := q^n e^{(-c-\epsilon_1)
  \sqrt{m}}$ colorings, where $c = 2\sqrt{\log \frac{q}{q-1} \log
  q}$. Since $u_1 = \sqrt{m \cdot \log \frac{q}{q-1} / \log q}$ and
$u_2 = \sqrt{m \cdot \log q / \log \frac{q}{q-1}}$,  routine
algebra shows that $\Sigma_0$ is precisely $e^{-\epsilon_1 \sqrt{m}}
(q-1)^{u_2} q^{n-u_1-u_2}$.  Therefore, for small $\epsilon_1$ we
have $\Sigma_1 / \Sigma_0 \leq e^{-\Theta(\sqrt{m})} = o(1)$, i.e.,
almost all colorings of $G$ are $(X_1, X_2)$-regular.  \hfill $\Box$

\vspace{3mm}

Before proving the next claim, it is convenient to establish the
following lemma, which should be understood in the context of Claim 3.

\begin{lemma}
  \label{lem:exact:sparse:sum-degrees}
  Let $x,y$ be a pair of non-isolated vertices of $G$, such that $xy$
  is not an isolated edge.  Then $d(x) + d(y) \geq |X_1| - 1$.
\end{lemma}

\noindent \textbf{Proof.}\, Suppose for contradiction that there is
such a pair $x,y$ with $d(x) + d(y) \leq |X_1| - 2$.  Let $G'$ be the
graph obtained by deleting the $\leq |X_1|-2$ edges incident to $x$ or
$y$, and adding back as many edges between $x$ and $X_1 \setminus
\{x,y\}$.  In $G'$, any $(X_1 \setminus \{x,y\}, X_2 \setminus
\{x,y\})$-regular partial coloring\footnote{A proper coloring of the
  vertices $V \setminus \{x,y\}$, which uses only one color on $X_1
  \setminus \{x,y\}$, and avoids that color on $X_2 \setminus
  \{x,y\}$.} of $V \setminus \{x,y\}$ has exactly $q-1$ extensions to
$x$ since only one color appears on $N_{G'}(x) \subset X_1 \setminus
\{x,y\}$, and then exactly $q$ further extensions to the
newly-isolated vertex $y$.  
On the other hand, since $x$ and $y$ both have degree at least 1 and do not form an isolated edge, one of them, 
say $x$, has a  neighbor in the rest of the graph. Therefore, in $G$ 
the same partial coloring has at most $q-1$ extensions to the vertex $x$, 
and then at most $q-1$ further
extensions to the non-isolated vertex $y$.  Yet by Claim 2, almost all
colorings of $G$ arise in this way, so for sufficiently large $m$, $G$
has fewer colorings than $G'$, contradiction.  \hfill $\Box$

\vspace{3mm}

\noindent \textbf{Proof of Claim 3.}\, Recall that our initial
technical adjustment allows us to assume that if $G$ contains an
isolated edge $uv$, then it also contains a degree-1 vertex $x \not
\in \{u,v\}$.  This would give $d(x) + d(u) = 2 \ll |X_1| - 1$,
contradicting Lemma \ref{lem:exact:sparse:sum-degrees} because $xu$
cannot be an isolated edge.  Hence $G$ in fact has no isolated edges.
But then the same lemma implies that at most one vertex $v_0$ has
degree $\leq 2 \epsilon_3 \sqrt{m}$, since $|X_1| = \Theta(\sqrt{m})$
by Claim 1.

It remains to show that each $|V_i|$ is within $O(\epsilon_2
\sqrt{m})$ of $u_i$.  Recall that $U_1'$ and $U_2''$ are the the
$U_1$- and $U_2$-vertices that are in $V_1$.  All other vertices of
$V_1$ are isolated in the graph $G^*$ which is within edit-distance
$\epsilon_1 m$ of $G$.  So by the previous paragraph, each of them
(except $v_0$ if it exists) has degree at least $2 \epsilon_3
\sqrt{m}$, and thus contributes at least $2 \epsilon_3 \sqrt{m}$ to
the edit distance between $G$ and $G^*$.  Therefore, there are at most
$1 + \frac{\epsilon_1 m}{2 \epsilon_3 \sqrt{m}} \ll \epsilon_2
\sqrt{m}$ of them, where we used $\epsilon_3^3 = \epsilon_2^2 =
\epsilon_1$.  Claim 1 controls $|U_i'|$ and $|U_i''|$, so we indeed
find that $|V_1|$ is within $O(\epsilon_2 \sqrt{m})$ of $u_1$.  The
analogous result for $V_2$ follows by a similar argument.  \hfill
$\Box$

\vspace{3mm}

\noindent \textbf{Proof of Claim 4.}\, Since almost all colorings are
$(X_1, X_2)$-regular, it suffices to prove this claim only for those
colorings.  So, we bound the $(X_1, X_2)$-regular colorings that
\textbf{(i)} use a common color on both $V_2^*$ and $V_1^*$, or
\textbf{(ii)} use at most $q-2$ colors on $V_2^*$.  Note that every
$(X_1, X_2)$-regular coloring which avoids both (i) and (ii) must use
exactly $q-1$ colors on $V_2^*$ and only the remaining color on
$V_1^*$, and so is automatically $(V_1^*, V_2^*)$-regular.  It therefore
suffices to show that these two types of colorings constitute
$o(1)$-fraction of all colorings.  The key observation is that every
$v \in V_2^*$ has a neighbor in $X_1$.  Indeed, $(V_1, V_2)$ is a
max-cut, so at least half of the $\geq 2 \epsilon_3 \sqrt{m}$
neighbors of $v$ must be in $V_1$.  These cannot all avoid $X_1$,
because Claims 1 and 3 show that only $O(\epsilon_2 \sqrt{m})$
vertices of $V_1$ are outside $X_1$, and $\epsilon_2 \ll \epsilon_3$.

To bound the number of colorings of type (i) above, first choose a
color $c_1$ for all $X_1$.  By the key observation, $c_1$ cannot
appear on $V_2^*$, so the shared color $c_2$ must be different.  Hence
we have $q-1$ choices for $c_2$, and must pick a pair of vertices $x
\in V_1^* \setminus X_1$ and $y \in V_2^*$ to use it on.  The $\geq
\epsilon_3 \sqrt{m}$ neighbors of $x$ in $V_2^*$ must avoid $c_2$ as
well as $c_1$, so they each have at most $q-2$ color choices.  Every
other vertex of $V_2^*$ must still avoid $c_1$, so we use the bound of
$\leq q-1$ color choices there.  Using the trivial bound $\leq q$ for
all other vertices, and the fact that $|X_i|$ and $|V_i^*|$ are within
$O(\epsilon_2 \sqrt{m})$ of $u_i = \Theta(\sqrt{m})$, we find that the
number of type-(i) colorings is at most:
\begin{eqnarray*}
  \Sigma_2 &:=& q \cdot (q-1) \cdot |V_1^* \setminus X_1| |V_2^*| \cdot
  (q-2)^{\epsilon_3 \sqrt{m}} \cdot (q-1)^{|V_2^*| - \epsilon_3 \sqrt{m}} \cdot
  q^{n-|X_1|-|V_2^*|-1} \\
  &\leq& O(m) \cdot \left( \frac{q-2}{q-1} \right)^{\epsilon_3 \sqrt{m}} \cdot
  (q-1)^{|V_2^*|} \cdot q^{n-|X_1|-|V_2^*|-1} \\
  &\leq& e^{O(\epsilon_2 \sqrt{m})} \cdot \left( \frac{q-2}{q-1} \right)^{\epsilon_3 \sqrt{m}} \cdot
  (q-1)^{u_2} \cdot
  q^{n-u_1-u_2}.
\end{eqnarray*}
On the other hand, we showed at the end of the proof of Claim 2 that
$G$ had at least $\Sigma_0 = e^{-\epsilon_1 \sqrt{m}} (q-1)^{u_2}
q^{n-u_1-u_2}$ colorings.  Since $\epsilon_1 \ll \epsilon_2 \ll
\epsilon_3$, we have $\Sigma_2 / \Sigma_0 \leq e^{-\Theta(\epsilon_3
  \sqrt{m})} = o(1)$, as desired.

The number of type-(ii) colorings is easily bounded by $\Sigma_3 := q
\cdot (q-1) \cdot (q-2)^{|V_2^*|} \cdot q^{n-|X_1|-|V_2^*|}$.  The
four factors correspond to choosing a color for $X_1$, choosing
another color to avoid on $V_2^*$, coloring $V_2^*$, and coloring all
remaining vertices.  Using that $|X_i|$ and $|V_i^*|$ are within
$O(\epsilon_2 \sqrt{m})$ of $u_i$, we obtain $\Sigma_3 \leq
e^{O(\epsilon_2 \sqrt{m})} (q-2)^{u_2} q^{n-u_1-u_2}$, so $\Sigma_3 /
\Sigma_0 \leq e^{O(\epsilon_2 \sqrt{m})}
\big(\frac{q-2}{q-1}\big)^{u_2}$.  Since $u_2 = \Theta(\sqrt{m})$, for
small enough $\epsilon_2$ we indeed have $\Sigma_3 / \Sigma_0 \leq
e^{-\Theta(\sqrt{m})} = o(1)$, as desired.  \hfill $\Box$

\vspace{3mm}

\noindent \textbf{Proof of Claim 5.}\, Almost all colorings are
$(V_1^*, V_2^*)$-regular, so $G[V_1^*]$ spans no edges.  We turn our
attention to $V_2^*$, and start by showing that all degrees within
$G[V_2^*]$ are at most $\epsilon_3 \sqrt{m}$.  Indeed, suppose for
contradiction that some $x \in V_2^*$ has at least $\epsilon_3
\sqrt{m}$ neighbors in $V_2^*$.  Then the number of $(V_1^*,
V_2^*)$-regular colorings is at most $\Sigma_4 := q \cdot (q-1) \cdot
(q-2)^{\epsilon_3 \sqrt{m}} \cdot (q-1)^{|V_2^*| - \epsilon_3
  \sqrt{m}} \cdot q^{n-|V_1^*|-|V_2^*|}$.  Here, the factors
correspond to choosing a color $c_1$ for $|V_1^*|$, choosing a color
$c_2$ for $x$, coloring $V_2^* \cap N(x)$ without $c_1$ or $c_2$,
coloring the rest of $V_2^*$ without $c_1$, and coloring the remaining
vertices.  Using that each $|V_i^*|$ is within $O(\epsilon_2
\sqrt{m})$ of $u_i$, we find that
\begin{eqnarray*}
  \Sigma_4 &\leq& e^{O(\epsilon_2 \sqrt{m})} \cdot q \cdot (q-1) \cdot
  (q-2)^{\epsilon_3 \sqrt{m}} \cdot (q-1)^{u_2 - \epsilon_3
    \sqrt{m}} \cdot q^{n - u_1 - u_2} \\
  &\leq& e^{O(\epsilon_2 \sqrt{m})} \cdot \left(
  \frac{q-2}{q-1} \right)^{\epsilon_3 \sqrt{m}} \cdot (q-1)^{u_2}
  q^{n-u_1-u_2}.
\end{eqnarray*}
Yet we showed at the end of the proof of Claim 2 that $G$ had at least
$\Sigma_0 = e^{-\epsilon_1 \sqrt{m}} (q-1)^{u_2} q^{n-u_1-u_2}$
colorings, so using $\epsilon_1 \ll \epsilon_2 \ll \epsilon_3$, we
obtain $\Sigma_4 / \Sigma_0 \leq e^{-\Theta(\epsilon_3 \sqrt{m})}$.
This contradicts the fact that $\Sigma_4$ includes almost all
colorings.  Therefore, all degrees within $G[V_2^*]$ are indeed at
most $\epsilon_3 \sqrt{m}$.

% now show that all internal degrees are zero

We now use this intermediate bound to show that all such degrees are
in fact zero.  Suppose for contradiction that some $x \in V_2^*$ has
neighbors within $V_2^*$.  Let $G'$ be the graph obtained by deleting
all edges between $x$ and $V_2^*$ and all edges incident to $v_0$ (if
it exists), and adding back as many edges between $V_1^*$ and some
formerly isolated vertex $z$.\footnote{Isolated vertices exist because
  Claim 3 shows that each $|V_i|$ is within $O(\epsilon_2 \sqrt{m})$
  of $u_i$, so the number of non-isolated vertices is $|V_1 \cup V_2|
  \leq u_1 + u_2 + O(\epsilon_2 \sqrt{m})$.  This is strictly below
  $n$ for small $\epsilon_2$, because $u_1 + u_2 = \sqrt{m/\kappa_q}$,
  and we assumed that $m \leq \kappa n^2$ with $\kappa < \kappa_q$.}
This is possible because $d(v_0) \leq 2\epsilon_3 \sqrt{m}$ and $x$
has at most $\epsilon_3 \sqrt{m}$ neighbors within $V_2^*$, while
$|V_1^*| = \Theta(\sqrt{m})$.  Observe that any $(V_1^*, V_2^*
\setminus \{x\})$-regular partial coloring of $V \setminus
\{x,z,v_0\}$ has exactly $(q-1)^2 q^{|V_0|}$ extensions to all of
$G'$, because $x$ and $z$ only need to avoid the single color which
appears on $V_1^*$, and $v_0$ is now isolated, if it exists.  On the
other hand, we claim that the same partial coloring has at most
$(q-2)q(q-1)^{|V_0|}$ extensions in $G$.  Indeed, there are at most
$q-2$ extensions to $x$ because $x$ must avoid the color of $V_1^*$ as
well as some (different) color which appears on its neighbor in
$V_2^*$.  Then, there are $q$ ways to color the isolated vertex $z$,
and finally at most $q-1$ further extensions to the non-isolated
vertex $v_0$ if it exists.  Yet by Claim 2, almost all colorings of
$G$ arise in this way, so for sufficiently large $m$, $G$ has fewer
colorings than $G'$.  This is impossible, so $V_2^*$ must indeed be an
independent set.

% now show $v_0$ has neighbors only on one side

It remains to show that $v_0$, if it exists, has neighbors in only one
$V_i^*$.  Suppose for contradiction that $v_0$ is adjacent to both
$V_1^*$ and $V_2^*$, and consider the graph $G'$ obtained by deleting all edges
incident to $v_0$, and replacing them with edges to $V_1^*$ only.
This is possible because $d(v_0) \leq 2\epsilon_3 \sqrt{m}$ and
$|V_1^*| = \Theta(\sqrt{m})$.  Any partial $(V_1^*, V_2^*)$-regular
coloring of $G \setminus \{v_0\}$ has at most $q-2$ extensions to
$v_0$, because $v_0$'s neighbors in $V_2^*$ are colored differently
from its neighbors in $V_1^*$.  Yet the same partial coloring has
exactly $q-1$ extensions with respect to $G'$, since it uses the same
color on all of $v_0$'s neighbors (now in $V_1^*$).  So, for
sufficiently large $m$, $G'$ has more colorings than $G$, giving the
required contradiction.  \hfill $\Box$

\vspace{3mm}

\noindent \textbf{Proof of Claim 6.}\, First, consider the case when
$V_0$ is empty.  Then all non-isolated vertices are already in the
bipartite graph $(V_1^*, V_2^*)$.  If that subgraph is less than
$|V_1^*|$ edges away from being complete bipartite, then Lemma
\ref{lem:semi-complete} already implies\footnote{$V_1^*$ is the
  smaller side of the bipartite graph $(V_1^*, V_2^*)$ because Claim 3
  shows that $|V_1^*|$ is within $O(\epsilon_2 \sqrt{m})$ of $u_1 =
  \sqrt{m \cdot \log \frac{q}{q-1} / \log q}$ and $|V_2^*|$ is within
  $O(\epsilon_2 \sqrt{m})$ of $u_2 = \sqrt{m \cdot \log q / \log
    \frac{q}{q-1}}$.}  that $G[V_1^* \cup V_2^*]$ is semi-complete
(and correctly oriented if $q \geq 4$), so we are done.  On the other
hand, if that subgraph has at least $|V_1^*|$ missing edges, then we
can construct an $n$-vertex graph $G'$ with at least $m$ edges by
taking $K_{|V_1^*|, |V_2^*| - 1}$ and adding enough isolated vertices.
Then, $G'$ has at least $q(q-1)^{|V_2^*|-1}q^{n-|V_1^*|-|V_2^*|+1}$
colorings because there are $q$ choices of a single color for the
$|V_1^*|$-side, $q-1$ color choices for each vertex on the other side,
and $q$ choices for each remaining (isolated) vertex.  However, the
same counting shows that $G$ has exactly $q(q-1)^{|V_2^*|}
q^{n-|V_1^*|-|V_2^*|}$ colorings that are $(V_1^*, V_2^*)$-regular,
which includes almost all colorings by Claim 4.  Hence for
sufficiently large $m$, $G'$ has more colorings, and this
contradiction completes the case when $V_0$ is empty.

Now suppose the vertex $v_0$ with degree $\leq 2\epsilon_3 \sqrt{m}$
exists.  By counting $(V_1^*, V_2^*)$-regular colorings, we find that
$G$ has at most $\Sigma_5 := (1+o(1))
q(q-1)^{|V_2^*|}(q-1)q^{n-|V_1^*|-|V_2^*|-1}$ colorings.  Here, the
factors correspond to choosing a color for $V_1^*$, coloring $V_2^*$,
coloring the non-isolated vertex $v_0$ which must avoid a neighbor's
color, and coloring the remaining vertices.  Observe that if there
were at least $d(v_0)$ edges missing between $V_1^*$ and $V_2^*$, then
we could isolate $v_0$ by deleting its edges and adding back as many
between $V_1^*$ and $V_2^*$.  The resulting graph would have at least
$q(q-1)^{|V_2^*|} q^{n-|V_1^*|-|V_2^*|}$ colorings, where the factors
correspond to choosing a color for $V_1^*$, coloring $V_2^*$, and
coloring the remaining (isolated) vertices.  For sufficiently large
$m$, this exceeds the number of colorings of $G$, which is impossible.
Therefore, less than $d(v_0)$ edges are missing between $(V_1^*,
V_2^*)$.

By Claim 5, $v_0$ has neighbors in only one $V_i^*$.  If it is
$V_1^*$, we must have $V_1 = V_1^*$ and $V_2 = V_2^* \cup \{v_0\}$
because $(V_1, V_2)$ is a max-cut.  The previous paragraph then
implies that less than $|V_1|$ edges are missing between $(V_1, V_2)$,
so Lemma \ref{lem:semi-complete} shows that $G$ is indeed
semi-complete on its non-isolated vertices (and correctly oriented if
$q \geq 4$).

The only remaining case is when $v_0$ has neighbors only in $V_2^*$,
which we will show is impossible.  This time, the max-cut gives $V_1 =
V_1^* \cup \{v_0\}$ and $V_2 = V_2^*$.  Since $d(v_0) \leq 2\epsilon_3
\sqrt{m}$, there are at least $|V_2| - 2\epsilon_3 \sqrt{m}$ missing
edges between $(V_1, V_2)$.  So, if we let $t = \big\lfloor
\frac{|V_2| - 2\epsilon_3 \sqrt{m}}{|V_1|} \big\rfloor = \big\lfloor
\frac{u_2}{u_1} - O(\epsilon_3) \big\rfloor = \big\lfloor \log q /
\log \frac{q}{q-1} - O(\epsilon_3) \big\rfloor$, we can construct an
$n$-vertex graph $G'$ with at least $m$ edges by taking $K_{|V_1|,
  |V_2|-t}$ and adding enough isolated vertices.  This graph has at
least $\Sigma_6 := q(q-1)^{|V_2|-t} q^{n-|V_1|-|V_2|+t}$ colorings, by
the same counting as earlier in this proof.  Let us compare this with
the number of colorings $\Sigma_5$ of $G$, which we calculated above.
Since $|V_1^*| = |V_1| - 1$ and $|V_2^*| = |V_2|$, we have $\Sigma_6 /
\Sigma_5 \geq (1-o(1)) \big( \frac{q}{q-1} \big)^t \cdot
\frac{1}{q-1}$.

Crucially, $\log q / \log \frac{q}{q-1}$ is always irrational, because
any positive integral solution to $q^x = \big(\frac{q}{q-1}\big)^y$
would require $q$ and $q-1$ to have a nontrivial common factor.  So,
by choosing our $\epsilon$'s sufficiently small in advance (based only
on $q$), we may ensure that $t \geq \log q / \log \frac{q}{q-1} - 1 +
c_q$ for some small positive constant $c_q$.  Since
$\big(\frac{q}{q-1}\big)^{\log q / \log \frac{q}{q-1} - 1} \cdot
\frac{1}{q-1} = 1$, this gives $\Sigma_6 / \Sigma_5 \geq (1-o(1))
\big( \frac{q}{q-1} \big)^{c_q}$, which exceeds 1 for large $m$,
leaving $G'$ with more colorings than $G$.  This contradiction
finishes our last case, and our entire proof.  \hfill $\Box$

\section{Exact result for 3 colors}
\label{sec:exact:q=3}

Our arguments can be pushed further when only three colors are used.
In this section, we complete the proof of Theorem \ref{thm:main:q=3},
determining the precise structure of the graphs that maximize the
number of 3-colorings, for edge densities up to $m \leq \frac{1}{4}
n^2$ (i.e., up to the density of the complete bipartite graph).  The
structure of this proof closely resembles that of the previous
section, so parts that are essentially the same are rewritten briefly.

We would, however, like to draw attention to a new piece of notation.
Recall that, as defined in the previous section, a coloring is $(X,
Y)$-regular if it uses only one color on $X$ and the other $q-1$ on
$Y$.  This time, we will also need a symmetric version of this
concept, which we denote with square brackets.  We will say that a
coloring is \emph{$[X, Y]$-regular} if one of $X$ or $Y$ is monochromatic,
and the other avoids that color entirely.  Note that this is
equivalent to having no colors shared between $X$ and $Y$, because
there are only 3 colors altogether.

\vspace{3mm}

\noindent \textbf{Proof of Theorem \ref{thm:main:q=3}.}\, Theorem
\ref{thm:main:sparse} already established our result for densities up
to $m \leq \kappa n^2$ for some constant $\kappa$, so we may assume
that $m = \Theta(n^2)$.  Routine algebra verifies that Proposition
\ref{prop:solve-opt} and Theorem \ref{thm:asymp-number} establish the
claimed numbers of colorings in this theorem.  This leaves us to
concentrate on the optimal graph structure.  We use several constants
$\epsilon_1 \ll \epsilon_2 \ll \epsilon_3$, related by $\epsilon_1 =
\epsilon_2^2 = \epsilon_3^3$, and show that there is an eventual
choice that makes our argument work.  To avoid confusion, our $O$,
$\Theta$, and $o$ notation will only mask absolute constants.

Let $G = (V, E)$ be an optimal graph whose density $m/n^2$ is between
$\kappa$ and $1/4$.  Let $u_1 = \alpha_3 n$ and $u_2 = \alpha_{12} n$,
where the $\alpha$'s are determined by Proposition
\ref{prop:solve-opt} with density parameter $\gamma = m/n^2$.  Note
that since $\kappa \leq \gamma \leq \frac{1}{4}$, each $u_i =
\Theta(n)$.  Theorem \ref{thm:asymp-stability} gives disjoint subsets
$U_1, U_2 \subset V$ with $|U_i| \in \{\lfloor u_i \rfloor, \lceil u_i
\rceil\}$, such that by editing at most $\epsilon_1 n^2$ edges, we can
transform $G$ into the complete bipartite graph between $U_1$ and
$U_2$, plus isolated vertices.  Call that graph $G^*$.

Let $(V_1, V_2)$ be a max-cut partition of the \textbf{non-isolated}
vertices of $G$, such that $V_1$ contains at least as many vertices of
$U_1$ as $V_2$ does.  Define $U_i' = U_i \cap V_i$ and $U_i'' = U_i
\cap V_{3-i}$, and let $X_i \subset U_i'$ be the vertices that are
adjacent to all but at most $\epsilon_2 n$ vertices of $U_{3-i}'$.
The following series of claims will complete the proof of Theorem
\ref{thm:main:q=3}.

\begin{description}

\item[Claim 1.] For each $i$, $|U_i'|$ is within $O(\epsilon_1 n)$ of
  $u_i$, $|X_i|$ is within $O(\epsilon_2 n)$ of $u_i$, and
  $|U_i''| \leq O(\epsilon_1 n)$.

\item[Claim 2.] Almost all colorings of $G$ are \emph{$[X_1,
    X_2]$-regular}, meaning that one $X_i$ is monochromatic, and the
  other $X_{3-i}$ avoids that color entirely.

\item[Claim 3.] All nonzero degrees are at least $2 \epsilon_3 n$,
  except possibly for either (i) only one isolated edge $w_1 w_2$, or
  (ii) only one non-isolated vertex $v_0$.  We use this to show that
  each $|V_i|$ is within $O(\epsilon_2 n)$ of $u_i$.  Let $V_0 =
  \{w_1, w_2\}$ if exception (i) occurs, let $V_0 = \{v_0\}$ if (ii)
  occurs, and let $V_0 = \emptyset$ otherwise.  Let $V_i^* = V_i
  \setminus V_0$.

\item[Claim 4.]  Almost all colorings are $[V_1^*, V_2^*]$-regular.

\item[Claim 5.] Each $V_i^*$ is an independent set, and $v_0$ (if it
  exists) has neighbors in only one of the $V_i^*$.  Hence $G$ is a
  bipartite graph plus isolated vertices.

\item[Claim 6.] $G$ is either a semi-complete subgraph of $K_{|V_1|,
    |V_2|}$ plus isolated vertices, or a complete bipartite subgraph
  $K_{|V_1^*|, |V_2^*|}$ plus a pendant edge to $v_0$.

\end{description}

\subsection{Supporting claims}

\noindent \textbf{Proof of Claim 1.}\, The sets $|U_i| = \Theta(n)$
are complete to each other in $G^*$, so all $U_i$-vertices have degree
$\Theta(n)$ in $G^*$.  As $G$ is at most $\epsilon_1 n^2$ edges away
from $G^*$, the number of $U_i$-vertices that are isolated in $G$ is
at most $\frac{\epsilon_1 n^2}{\Theta(n)} = O(\epsilon_1 n)$.  Since
$V_1$ received more non-isolated $U_1$-vertices than $V_2$ did, we
must have $|U_1'| \geq \frac{1}{3} u_1 = \Theta(n)$.  By Proposition
\ref{prop:construction-asymp-edges}, $G^*$ has at least $m - O(n)$
edges, all of which cross between $(U_1, U_2)$.  So $G$ has at least
$m - O(n)-\epsilon_1 n^2$ edges there, and at least that many between
its max-cut $(V_1, V_2)$.  As $G$ has only $m$ edges, this shows that
each $G[V_i]$ spans $O(\epsilon_1 n^2)$ edges.  But the sets $U_1',
U_2'' \subset V_1$ are complete to each other in $G^*$, so $|U_1'|
|U_2''| - \epsilon_1 n^2 \leq e(G[V_i]) \leq O(\epsilon_1 n^2)$.
Using $|U_1'| \geq \Theta(n)$, we indeed obtain $|U_2''| \leq
O(\epsilon_1 n)$.

Then $|U_2'| \geq u_2 - O(\epsilon_1 n) \geq \Theta(n)$, because only
$O(\epsilon_1 n)$ of the $U_2$-vertices are isolated and $|U_2''| \leq
O(\epsilon_1 n)$ of them are in $V_1$.  So, repeating the above with
respect to $U_2'$ and $U_1''$ instead of $U_1'$ and $U_2''$, we find
that $|U_1''| \leq O(\epsilon_1 n)$, which then implies that $|U_1'|
\geq u_1 - O(\epsilon_1 n)$.

To control $X_i$, observe that since the $U_i'$ are complete to each
other in $G^*$, each vertex not in $X_i$ contributes at least
$\epsilon_2 n$ to the total edit distance of $\leq \epsilon_1 n^2$
between $G$ and $G^*$.  We set $\epsilon_2^2 = \epsilon_1$, so all but
at most $\epsilon_2 n$ vertices of $U_i'$ belong to $X_i$.  Since
$|U_i'|$ is within $O(\epsilon_1 n)$ of $u_i$, this gives the desired
result.  \hfill $\Box$

\vspace{3mm}

\noindent \textbf{Proof of Claim 2.}\, For each partition $\{1, 2, 3\}
= C_0 \cup C_1 \cup C_2 \cup C_3$, we count the colorings which use
the colors $C_1$ in $X_1$ but not $X_2$, use $C_2$ in $X_2$ but not
$X_1$, use $C_3$ in both $X_1$ and $X_2$, and do not use $C_0$ in
either $X_1$ or $X_2$.  Then we sum over all \emph{irregular}\/
partitions, which are all partitions with $|C_3| \geq 1$.  Note that a
coloring is $[X_1, X_2]$-regular if and only if it does not use any
color on both $X_i$, so this sum will include all other colorings.

For any given partition with $|C_i| = c_i$, we have that the corresponding number
of colorings is at most $(|X_1| |X_2|)^{c_3} \cdot c_1^{|X_1| - 3
  \epsilon_2 n} \cdot c_2^{|X_2| - 3 \epsilon_2 n} \cdot 3^{n - 2c_3 -
  (|X_1| - 3 \epsilon_2 n) - (|X_2| - 3 \epsilon_2 n)}$, by the
calculation in Claim 2 of Section \ref{sec:exact:sparse:details} with
$q$ replaced by 3 and $\sqrt{m}$ replaced by $n$.  Using that each
$|X_i|$ is within $O(\epsilon_2 n)$ of $u_i = \Theta(n)$ and all
irregular colorings have $|C_3| \geq 1 \Rightarrow c_1+c_2 \leq 2$, we
find that the sum $\Sigma_1$ of this bound over all $\leq 4^3$
irregular partitions is:
\begin{eqnarray*}
  \Sigma_1 &=& \sum_{\text{irregular}} (|X_1| |X_2|)^{c_3} \cdot c_1^{|X_1|
    - 3 \epsilon_2 n} \cdot c_2^{|X_2| - 3 \epsilon_2 n}
  \cdot 3^{n - 2c_3 - (|X_1| - 3 \epsilon_2 n) - (|X_2| - 3 \epsilon_2 n)} \\
  &\leq& e^{O(\epsilon_2 n)} \sum_{\text{irregular}} (\Theta(n) \cdot \Theta(n))^{c_3} \cdot c_1^{u_1} \cdot c_2^{u_2}
  \cdot 3^{n - u_1 - u_2} \\
  &\leq& e^{O(\epsilon_2 n)} \cdot 4^3 \cdot O(n^6) \cdot 
  \max_{c_1 + c_2 \leq 2} \left\{ c_1^{u_1} c_2^{u_2} \right\} \cdot
  3^{n - u_1 - u_2}
  \ \ = \ \ e^{O(\epsilon_2 n)} \cdot 3^{n - u_1 - u_2}.
\end{eqnarray*}
On the other hand, Proposition \ref{prop:solve-opt}, Theorem
\ref{thm:asymp-number}, and routine algebra show that just as in the
sparse case, the optimal graph has at least $\Sigma_0 :=
e^{-\epsilon_1 n} \cdot 2^{u_2} \cdot 3^{n-u_1-u_2}$ colorings.  Using
$u_2 = \Theta(n)$, we find that $\Sigma_1 / \Sigma_0 \leq
e^{-\Theta(n)} = o(1)$, i.e., almost all colorings of $G$ are $[X_1,
X_2]$-regular.  \hfill $\Box$

\vspace{3mm}

Before proving the next claim, it is convenient to establish the
following lemma, which should be understood in the context of Claim 3.

\begin{lemma}
  \label{lem:exact:q=3:sum-degrees}
  Let $x,y$ be a pair of non-isolated vertices of $G$, such that $xy$
  is not an isolated edge.  Then $d(x) + d(y) \geq \min\{|X_1|,
  |X_2|\} - 1$.
\end{lemma}

\noindent \textbf{Proof.}\, Suppose for contradiction that there is
such a pair $x,y$ with $d(x) + d(y) \leq \min\{|X_1|, |X_2|\} - 2$.
Also suppose that among the $[X_1 \setminus \{x,y\}, X_2 \setminus
\{x,y\}]$-regular partial colorings of $V \setminus \{x,y\}$,
at least half of them have $X_1 \setminus \{x,y\}$ monochromatic.
(The case when at least half have $X_2 \setminus \{x,y\}$
monochromatic follows by a similar argument.)  Let $G'$ be the graph
obtained by deleting the $\leq |X_1|-2$ edges incident to $x$ or
$y$, and adding back as many edges between $x$ and $X_1 \setminus
\{x,y\}$.

Consider any $[X_1 \setminus \{x,y\}, X_2 \setminus \{x,y\}]$-regular
partial coloring of $V \setminus \{x,y\}$.  If it is monochromatic in
$X_1$, which happens at least half the time by assumption, then in $G'$ it has
exactly 2 extensions to $x$, followed by 3 further extensions to the
newly-isolated vertex $y$.  The rest of the time, the partial coloring
is monochromatic in $X_2$ and uses at most 2 colors in $X_1$.  Then,
in $G'$ it has at least 1 extension to $x$, followed by 3 further
extensions to $y$.

On the other hand, since $x$ and $y$ both have degree at least 1 and do not form an isolated edge, one of them, 
say $x$, has a  neighbor in the rest of the graph.  Therefore, in $G$
the same partial coloring has at most $2$ extensions to the vertex $x$,
and then at most $2$ further extensions to the non-isolated vertex $y$. 
Yet by Claim 2, almost all colorings of $G$
arise in this way, so the ratio of $G'$-colorings to $G$-colorings is
at least $\frac{1}{2} \big(\frac{2 \cdot 3}{2 \cdot 2} + \frac{1 \cdot
  3}{2 \cdot 2}\big) - o(1) = \frac{9}{8} - o(1) > 1$, contradiction.
\hfill $\Box$

\vspace{3mm}

\noindent \textbf{Proof of Claim 3.}\, If there is an isolated edge
$w_1 w_2$, then Lemma \ref{lem:exact:q=3:sum-degrees} implies that any
other vertex $x$ has $d(x) + 1 = d(x) + d(w_1) \geq \min\{|X_1|,
|X_2|\} - 1 = \Theta(n)$, giving exception (i).  Otherwise, the same
lemma implies there is at most one vertex $v_0$ of degree $\leq 2
\epsilon_3 n$, giving exception (ii).  The rest of this claim, that
each $|V_i|$ is within $O(\epsilon_2 n)$ of $u_i$, follows by the same
argument as in Claim 3 of Section \ref{sec:exact:sparse:details}, but
with $\sqrt{m}$ replaced by $n$ throughout.  \hfill $\Box$

\vspace{3mm}

\noindent \textbf{Proof of Claim 4.}\, Note that a coloring is
$[V_1^*, V_2^*]$-regular if and only if it does not use any color on
both $V_i^*$.  So, we bound the colorings that share a color on both
$V_i^*$, but \textbf{(i)} use only one color on $X_1$ and a subset of the other two on
$X_2$, or \textbf{(ii)} one on $X_2$ and a subset of the other two on $X_1$.
Since almost all colorings are $[X_1, X_2]$-regular, it suffices to
show that these two types of colorings constitute $o(1)$-fraction of
all colorings.  The same calculation as in Claim 4 of Section
\ref{sec:exact:sparse:details}, with $q$ replaced by 3 and $\sqrt{m}$
replaced by $n$, shows that the number of type-(i) colorings is at
most:
%  The key observation is
%  that every $v \in V_2^*$ has a neighbor in $X_1$, because $d(v) \geq 2
%  \epsilon_3 n$ is too large to entirely fit within $|V_1 \setminus X_1|
%  \leq O(\epsilon_2 n)$.
%
%  To bound the number of colorings of type (i) above, first choose a
%  color $c_1$ for all $X_1$.  By the observation, $c_1$ cannot appear on
%  $V_2^*$, leaving 2 choices for the shared color $c_2$ which appears on
%  some pair of vertices $x \in V_1^* \setminus X_1$ and $y \in V_2^*$.
%  The $\geq \epsilon_3 n$ neighbors of $x$ in $V_2^*$ must avoid $c_2$
%  as well as $c_1$, so they each have $\leq 1$ color choice.  Every
%  other vertex of $V_2^*$ must still avoid $c_1$, so we use the bound of
%  $\leq 2$ choices there.  Using the trivial bound $\leq 3$ for all
%  other vertices, and the fact that $|X_i|$ and $|V_i^*|$ are within
%  $O(\epsilon_2 n)$ of $u_i = \Theta(n)$, we find that the number of
%  type-(i) colorings is at most:
\begin{eqnarray*}
  \Sigma_2 &:=& 3 \cdot 2 \cdot |V_1^* \setminus X_1| |V_2^*| \cdot
  1^{\epsilon_3 n} \cdot 2^{|V_2^*| - \epsilon_3 n} \cdot
  3^{n-|X_1|-|V_2^*|-1} \\
  &\leq& e^{O(\epsilon_2 n)} \cdot O(n^2) \cdot 2^{-\epsilon_3 n} \cdot 2^{u_2} \cdot 3^{n-u_1-u_2}.
\end{eqnarray*}
On the other hand, we showed at the end of the proof of Claim 2 that
$G$ had at least $\Sigma_0 = e^{-\epsilon_1 n} \cdot 2^{u_2} \cdot
3^{n-u_1-u_2}$ colorings.  Since $\epsilon_1 \ll \epsilon_2 \ll
\epsilon_3$, we have $\Sigma_2 / \Sigma_0 \leq e^{-\Theta(\epsilon_3
  n)} = o(1)$, as desired.  The analogous result for type-(ii)
colorings follows by a similar argument.  \hfill $\Box$

\vspace{3mm}

\noindent \textbf{Proof of Claim 5.}\, We first show that $v_0$ cannot
have neighbors in both $V_i^*$.  Suppose for contradiction that this
is not the case.  Almost all colorings are $[V_1^*, V_2^*]$-regular by
Claim 4, so there is $I \in \{1,2\}$ such that $V_I^*$ is
monochromatic in at least $\big(\frac{1}{2} - o(1)\big)$-fraction of
all colorings.  Let $G'$ be obtained by deleting the $\leq 2
\epsilon_3 n$ edges incident to $v_0$, and replacing them with edges
to $|V_I^*| = \Theta(n)$ only.  Consider any partial $[V_1^*,
V_2^*]$-regular coloring of $V \setminus \{v_0\}$.  If it uses only
one color on $V_I^*$ (which happens at least half the time by assumption), in $G'$
it has exactly 2 extensions to $v_0$.  The rest of the time, it still
uses at most 2 colors on $V_I^*$, so there is at least 1 extension.
On the other hand, in $G$ the same partial coloring always has at most
1 extension to $v_0$, because $v_0$'s neighbors in $V_1^*$ are colored
differently from its neighbors in $V_2^*$.  By Claim 2, almost all
colorings of $G$ arise in this way, so the ratio of number of
colorings of $G'$ to $G$ is at least $\frac{1}{2} \cdot
\big(\frac{2}{1} + \frac{1}{1}\big) - o(1) = \frac{3}{2} - o(1)$,
contradiction.  Therefore, $v_0$ cannot have neighbors in both
$V_i^*$, as claimed.

It remains to show that both $G[V_i^*]$ are empty.  Suppose for
contradiction that some $x \in V_2^*$ has neighbors within $V_2^*$.
(The analogous result for $V_1^*$ follows by a similar argument.)
Almost every coloring is $[V_1^*, V_2^*]$-regular, but $V_2^*$ can
never be monochromatic because it contains edges.  So, almost all
colorings are in fact $(V_1^*, V_2^*)$-regular.\footnote{Recall that
  round brackets denote ``ordered'' regularity, where $V_1^*$ is
  monochromatic, and $V_2^*$ has the other two colors.}  Therefore,
the same argument as in Claim 5 of Section
\ref{sec:exact:sparse:details}, with $q$ replaced by 3 and $\sqrt{m}$
replaced by $n$, shows that $x$ has at most $\epsilon_3 n$ neighbors
within $V_2^*$.

\vspace{2mm}

\textbf{Case 1: there is some $z_0 \in \boldsymbol{V_0}$.}\, Let $G'$
be obtained by deleting the $\leq \epsilon_3 n$ edges between $x$ and
$V_2^*$ and the $\leq 2 \epsilon_3 n$ edges incident to anything in
$V_0$, and adding back as many edges between $z_0$ and $|V_1^*| =
\Theta(n)$.  Every $(V_1^*, V_2^* \setminus \{x\})$-regular partial
coloring of $V \setminus (V_0 \cup \{x\})$ has exactly $2 \cdot 2
\cdot 3^{|V_0|-1}$ extensions to all of $G'$, because $x$ and $z_0$
only need to avoid the single color which appears on $V_1^*$, and the
rest of $V_0$ (if any) is now isolated.  On the other hand, in $G$ the
same partial coloring has at most 1 extension to $x$ because $x$ must
avoid the color of $V_1^*$ as well as some (different) color which
appears on its neighbor in $V_2^*$.  Then, it has at most
$3^{|V_0|-1}$ further extensions to $V_0 \setminus \{z_0\}$ by the
trivial bound, and at most 2 further extensions to the non-isolated
vertex $z_0$.  Note that all $(V_1^*, V_2^*)$-regular colorings of $G$
arise in this way, which is almost all of the total by our remark
before we split into cases.  Hence for sufficiently large $m$, $G$ has
fewer colorings than $G'$, contradiction.

\vspace{2mm}

\textbf{Case 2: $\boldsymbol{V_0} = \emptyset$, but there is some
  isolated vertex $\boldsymbol{z}$.}\, Define $G'$ by deleting the
$\leq \epsilon_3 n$ edges between $x$ and $V_2^*$, and adding back as
many edges between $z$ and $|V_1^*| = \Theta(n)$.  By the same
arguments as in Case 1, all $(V_1^*, V_2^* \setminus \{x\})$-regular
partial colorings of $V \setminus \{x,z\}$ have exactly $2 \cdot 2$
extensions to $G'$, but in $G$ they have at most 1 extension to $x$,
followed by 3 further extensions to the isolated $z$.  This produces
almost all colorings of $G$, so $G'$ has more colorings for large $m$,
contradiction.

\vspace{2mm}

\textbf{Case 3: $\boldsymbol{V_1^* \cup V_2^* = V}$.}\, We observed
that the edges in $V_2^*$ force almost all colorings to use only one
color for $V_1^*$ and the other two on $V_2^*$ (hence $G[V_2^*]$ is
bipartite).  There are 3 color choices for $V_1^*$, so the number of
colorings of $G$ is $(3+o(1)) \cdot \#\{\text{2-colorings of
  $V_2^*$}\}$.  Recall that the number of 2-colorings of any bipartite
graph $F$ is precisely $2^r$, where $r$ is its number of connected
components.

We claim that the bipartite $G[V_2^*]$ has at most $|V_2^*| -
2\sqrt{t} + 1$ components, where $t$ is the number of edges in
$G[V_2^*]$.  Indeed, for fixed $t$, the optimal configuration is to
have all isolated vertices except for a single nontrivial (bipartite)
component $C$.  The sizes $a,b$ of the sides of that bipartite $C$
should minimize $a+b$ subject to the constraint $ab \geq t$, so by the
inequality of the arithmetic and geometric means, we have $a+b \geq
2\sqrt{t}$, as desired.  Therefore, $G$ has at most $(3+o(1)) \cdot
2^{|V_2^*| - 2\sqrt{t} + 1}$ colorings.

Let $G'$ be the complete bipartite graph with sides $s$ and $n-s$,
such that $s$ is as large as possible subject to $s(n-s) \geq m$.
Note that $|V_1^*| \cdot |V_2^*| \geq m-t$ because all but $t$ of
$G$'s $m$ edges cross between the $V_i^*$, so Inequality
\ref{ineq:exact:q=3:claim5} routinely shows that $s \geq |V_2^*| -
\lceil \sqrt{t} \rceil$.  Since $G'$ is complete bipartite, it has
exactly $3 \cdot 2^s + 3 \cdot 2^{n-s} - 6$ colorings, and thus our
bound on $s$ implies that $G'$ has strictly more than $3 \cdot 2^s
\geq 3 \cdot 2^{|V_2^*| - \lceil \sqrt{t} \rceil}$ colorings.  Yet for
$t \geq 3$, one may check that $-\lceil \sqrt{t} \rceil \geq
(-2\sqrt{t} + 1) + 0.4$, giving $G'$ more colorings than $G$, which is
impossible.

We are left with the cases $t \in \{1,2\}$, but for these values there
is always a vertex $y \in V_2^*$ with exactly 1 neighbor $z$ in
$G[V_2^*]$.  This forces all edges to be present between the $V_i^*$,
because otherwise we could increase the number of $(V_1^*,
V_2^*)$-regular colorings by a factor of 2 by deleting the edge $yz$
and adding one of the missing edges between the $V_i^*$.  The presence
of the complete bipartite graph forces \emph{every}\/ coloring of $G$
to use exactly two colors on $V_2^*$, and the other on $V_1^*$.
Together with the observation that the maximum number of connected
components of $G[V_2^*]$ is $|V_2^*| - t$ when $t \in \{1,2\}$, we
find that $G$ has \emph{exactly}\/ $3 \cdot 2^r \leq 3 \cdot
2^{|V_2^*| - t}$ colorings.  On the other hand, we showed above that
$G'$ had more than $3 \cdot 2^{|V_2^*| - \lceil \sqrt{t} \rceil}$
colorings.  Since $t = \lceil \sqrt{t} \rceil$ for $t \in \{1,2\}$,
$G'$ has more colorings than $G$, contradiction.  \hfill $\Box$

\vspace{3mm}

\noindent \textbf{Proof of Claim 6.}\, Let $G_0 = G[V_1 \cup V_2]$ be
the graph formed by the non-isolated vertices of $G$, and let $n_0 =
|V_1 \cup V_2|$.  Since the number of colorings of $G$ is precisely
$3^{n - n_0}$ times the number of colorings of $G_0$, the optimality
of $G$ implies that $G_0$ must also be optimal among $n_0$-vertex
graphs with $m$ edges.  Furthermore, Claim 4 also implies that almost
all colorings of $G_0$ are $[V_1^*, V_2^*]$-regular.

\vspace{2mm}

\textbf{Case 1: $\boldsymbol{V_0}$ is empty.}\, Let $\{a,b\}$ be the
sizes of the $V_i^*$, with $a \leq b$.  If there are less than $a$
missing edges between the $V_i^*$, then Lemma \ref{lem:semi-complete}
shows that $G_0$ is semi-complete, so we are done.  On the other hand,
if there are at least $a$ missing edges, then $K_{a, b - 1}$ plus one
isolated vertex has $n_0$ vertices and at least $m$ edges, but also
exactly $(3 \cdot 2^a + 3 \cdot 2^{b-1} - 6) \cdot 3$ colorings.  Yet
$G_0$ has no vertices outside $V_1^* \cup V_2^*$, and almost all
colorings are $[V_1^*, V_2^*]$-regular, so $G_0$ has at most $(1+o(1))
\cdot (3 \cdot 2^a + 3 \cdot 2^b )$ colorings, which is smaller,
contradiction. \hfill $\Box$

\vspace{2mm}

\textbf{Case 2: $\boldsymbol{V_0}$ is the single edge $\boldsymbol{w_1
    w_2}$.}\, We show that this is impossible.  Let $\{a,b\}$ be the
sizes of the $V_i^*,$ with $a \leq b$.  Since there are always exactly
6 ways to color the endpoints $\{w_1, w_2\}$ of the isolated edge
independently of the rest of $V$, and almost all colorings are
$[V_1^*, V_2^*]$-regular, $G_0$ has $(6+o(1)) \cdot ( 3 \cdot 2^a + 3
\cdot 2^b )$ colorings.  Let $G'$ be the complete bipartite graph
$K_{a-1, b+3}$, and let $G''$ be the complete bipartite graph $K_{a-1,
  b+2}$ plus one isolated vertex.  Both graphs have the same number of
vertices as $G_0$, so it suffices to show that at least one of them
has more edges and more colorings than $G_0$.

Claim 3 gives $\frac{a}{b} \geq \frac{u_1}{u_2} - O(\epsilon_2)$, and
Proposition \ref{prop:solve-opt} implies that $\frac{u_1}{u_2} \geq
\frac{\log 3/2}{\log 3} \approx 0.37$.  So for small $\epsilon_2$ and
large $n$, we have that $ab + 3a - b - 3 > ab+1$, hence $G'$ has more
edges than $G_0$.  Also, $G'$ has $3 \cdot 2^{b+3} = 24 \cdot 2^b$
colorings that use only one color on the $(a-1)$-side and the other
two on the $(b+3)$-side.  We claim that this already exceeds the
number of colorings of $G_0$ whenever $b \geq a+2$.  Indeed, then $2^a
\leq \frac{1}{4} \cdot 2^b$, so the number of colorings of $G_0$ is at
most:
\begin{displaymath}
  (6+o(1)) \cdot ( 3 \cdot 2^a + 3 \cdot 2^b )
  \ \ \leq \ \ 
  (6+o(1)) \cdot \frac{5}{4} \cdot 3 \cdot 2^b
  \ \ = \ \
  (22.5 + o(1)) \cdot 2^b,
\end{displaymath}
which is indeed less than the number of colorings of $G'$.

It remains to consider $a \leq b \leq a+1$.  Here, $G''$ has $ab + 2a
- b - 2 > ab+1$ edges, and exactly $(3 \cdot 2^{a-1} + 3 \cdot 2^{b+2}
- 6)\cdot 3$ colorings.  Using $a \geq b-1$, this is at least
$(1-o(1)) \cdot \frac{17}{16} \cdot 3 \cdot 2^{b+2} \cdot 3 = (38.25 -
o(1)) \cdot 2^b$.  On the other hand, using $a \leq b$, the number of
colorings of $G_0$ is at most $(36 + o(1)) \cdot 2^b$, which is smaller.
Therefore, $G''$ is superior on this range, and we are done.  \hfill
$\Box$

\vspace{2mm}

\textbf{Case 3: $\boldsymbol{V_0}$ is the single vertex
  $\boldsymbol{v_0}$.}\, Let $I$ be the index (unique by Claim 5) such
that $V_I^*$ contains neighbors of $v_0$.  Let $J = 3-I$ be the other
index, and let $a = |V_I^*|$, $b = |V_J^*|$.  Note that $G_0$ is
bipartite with partition $(V_I^*, V_J^* \cup \{v_0\})$.  If at least
$d(v_0)$ edges are missing between $V_I^*$ and $V_J^*$, then we can
isolate $v_0$ while only adding edges between $V_I^*$ and $V_J^*$.
This increases the number of $[V_I^*, V_J^*]$-regular colorings by a
factor of $\frac{3}{2} + o(1)$, which is impossible.  So, less than $d(v_0)$
edges are missing between $V_I^*$ and $V_J^*$, which implies that less
than $a$ edges are missing between $V_I^*$ and $V_J^* \cup \{v_0\}$.
Hence $G_0$ is a subgraph of $K_{a, b+1}$ with less than $a$ missing
edges.

When $a \leq b+1$, Lemma \ref{lem:semi-complete} shows that $G_0$ is
semi-complete, as desired.  It remains to consider $a > b+1$.  Some
vertex of the set $V_I^*$ of size $a$ is complete to $V_J^* \cup \{v_0\}$, because
less than $a$ edges are missing between $V_I^*$ and $V_J^* \cup
\{v_0\}$.  But we also showed that less than $d(v_0) \leq 2 \epsilon_3
n \ll |V_J^*|$ edges are missing between $V_I^*$ and $V_J^*$, so some
vertex of $V_J^*$ must be complete to $V_I^*$.  Thus, Lemma
\ref{lem:subgraph-bipartite:q=3} implies that since $G_0$ is an
optimal graph, the missing edges $E(K_{a,b+1}) \setminus E(G_0)$ form
a star, which must have center $v_0$ because $d(v_0) \leq 2 \epsilon_3
n \ll \min\{a,b\}$.  In particular, the number of missing edges is
then exactly $a-d$, where $d = d(v_0)$, and then the same lemma shows
that $G_0$ has exactly $3 \cdot 2^a + 3 \cdot 2^{b+1} + 6 \cdot
(2^{a-d} - 2)$ colorings.

Consider the graph $G'$ obtained by removing a $(b-d)$-edge star from
the complete bipartite graph $K_{a+1,b}$.  This has as many vertices
and edges as $G_0$, and $3 \cdot 2^{a+1} + 3 \cdot 2^b + 6 \cdot
(2^{b-d} - 2)$ colorings by Lemma \ref{lem:subgraph-bipartite:q=3}.
The difference between the numbers of colorings of $G'$ and $G_0$ is
\begin{displaymath}
  3 \cdot 2^a - 3 \cdot 2^b + 6 \cdot (2^{b-d} - 2^{a-d}) 
  \ \ = \ \ \left(3 - \frac{6}{2^d}\right) \cdot (2^a - 2^b),
\end{displaymath}
which exceeds zero for $d \geq 2$ because we are in the case $a >
b+1$.  Optimality of $G_0$ thus forces $d(v_0) = 1$.

We showed there were less than $d(v_0)$ edges missing between the
$V_i^*$, so now we know that the non-isolated vertices of $G$ form a
complete bipartite subgraph $(V_1^*, V_2^*)$ plus a pendant edge to
$v_0$.  Finally, observe that $G$ cannot have any isolated vertex $z$,
or else we could replace the pendant edge with the (isolated) edge
$v_0 z$, and this would not change the number of colorings because
every partial coloring of $V \setminus \{v_0\}$ would still have
exactly 2 extensions to the degree-1 vertex $v_0$.  But the resulting
graph is not optimal by the same argument as in Case 2 of this claim.
Therefore, $G$ is only a complete bipartite subgraph plus a pendant
edge, with no isolated vertices.  This completes the final case of our
final claim, and our entire proof.  \hfill $\Box$

\section{Exact result for Tur\'an graphs}
\label{sec:exact:turan}

We now study the extremality of Tur\'an graphs.  As we mentioned in
the introduction, Lazebnik conjectured that Tur\'an graphs $T_r(n)$
were the unique graphs that maximized the number of $q$-colorings
whenever $r \leq q$.  Note that Theorem \ref{thm:main:q=3} implies
this result for $q=3$ and $r=2$ when $n$ is large, because it shows
that all optimal graphs are bipartite, and no other bipartite graph
has as many edges as $T_2(n)$.  In this section, we prove Theorem
\ref{thm:exact:turan}, which confirms (for large $n$) Lazebnik's
conjecture when $r = q-1$, for all remaining $q$.  Our proof relies on
the following special case of a result of Simonovits
\cite{Simonovits}.  Let $t_r(n)$ denote the number of edges of the
$r$-partite Tur\'an graph $T_r(n)$ with $n$ vertices.

\begin{fact}
  \label{fact:simonovits}
  Let $F$ be a graph with chromatic number $r+1$.  Suppose there is an
  edge whose deletion makes $F$ $r$-colorable.  Then for all
  sufficiently large $n$, the Tur\'an graph $T_r(n)$ is the unique
  $n$-vertex graph with at least $t_r(n)$ edges that does not contain
  a subgraph isomorphic to $F$.
\end{fact}

We use this fact to prove the following lemma, which we will need
later.

\begin{lemma}
  \label{lem:exact:turan}
  Let $q \geq 4$ be fixed.  The following holds for all sufficiently
  large $n$.  Let $G \neq T_{q-1}(n)$ have $n$ vertices, and at least
  as many edges and $q$-colorings as $T_{q-1}(n)$.  Let $\Delta$ be
  the difference between the number of edges of $G$ and $T_{q-1}(n)$,
  and let $n' = n-(q-1)$.  Then there is an $n'$-vertex graph $H$ with
  at least $\Delta + 1$ more edges than $T_{q-1}(n')$, and at least
  half as many $q$-colorings as $G$ has.
\end{lemma}

\noindent \textbf{Proof.}\, We begin with a convenient technical
adjustment.  If $G$ has $k \geq 2$ connectivity components $C_i$ that
are not isolated vertices, then choose vertices $v_i \in C_i$ and glue
the components together by merging all of the $v_i$ into a single
vertex $v$.  Add $k-1$ isolated vertices $w_1, \ldots, w_{k-1}$ to
restore the vertex count, and let $G'$ be the resulting graph.
Clearly, $G'$ has as many edges as $G$, and it also is not
$T_{q-1}(n)$ because $G'$ has a vertex whose deletion increases the
number of components while $T_{q-1}(n)$ does not.  Furthermore, we
claim that $G$ and $G'$ have the same number of colorings.  Indeed, by
symmetry, for an arbitrary color $c$, the total number of colorings of
$G$ is precisely $q^k$ times the number of colorings of $G$ which use
$c$ for every $v_i$.  The obvious correspondence gives a bijection
between these colorings and partial colorings of $G' \setminus \{w_1,
\ldots, w_{k-1}\}$ which use $c$ on the merged vertex $v$.  Yet the
$w_i$ are isolated, so each of these partial colorings has exactly
$q^{k-1}$ extensions to all of $G'$.  Again by symmetry, the total
number of colorings of $G'$ is precisely $q$ times the number that use
$c$ on $v$.  Putting everything together, we find that $G$ and $G'$
indeed have the same number of colorings.  Therefore, by replacing $G$
with $G'$, we may assume without loss of generality that $G$ has only
one nontrivial connectivity component.

Fact \ref{fact:simonovits} implies that for large $n$, $G$ has a
subgraph $F$ which is the complete $(q-1)$-partite graph on $V(F) =
X_1 \cup \ldots \cup X_{q-1}$ with each part $X_i = \{u_i, w_i\}$
consisting of two vertices, plus an extra edge $u_1 w_1$.  Let 
$U = \{u_1, \ldots, u_{q-1}\}$ and $W = \{w_1, \ldots, w_{q-1}\}$, and let
$A = U \cup \{w_1\}$.

Let $\delta$ be the difference between the number of edges of
$T_{q-1}(n)$ and $T_{q-1}(n')$.  We claim that if there is a set $Y$
of $q-1$ vertices of $A$ such that the sum of their degrees is at most
$\delta + {q-1 \choose 2} - 1$, then $H = G - Y$ satisfies the lemma's
assertion.  Clearly, $H$ has the correct number of vertices, and it
has the correct number of edges because $Y \subset A$ induces a
complete graph $K_{q-1}$, so the number of deleted edges is at most
$\delta - 1$.  We now show that every $q$-coloring of $H$ extends to
at most two $q$-colorings of $G$.

If $Y = U$, since $\{u_1\} \cup W$ induces a $K_q$-subgraph in $G$,
every coloring of $H \supset W$ has at most 1 extension to $u_1$.
Then, every other $u_i$ has at most 1 choice because $\{u_1, u_i\}
\cup (W \setminus \{w_i\})$ induces a $K_q$-subgraph in which $u_i$ is
the only uncolored vertex.  Thus when $Y = U$, every coloring of $H$
colors $W$ and hence has at most 1 extension to $G$.  On the other
hand, up to a symmetry of $F$, the only other case is when $Y =
\{w_1\} \cup (U \setminus \{u_{q-1}\})$.  As before, $\{u_1\} \cup W$
induces a $K_q$-subgraph in $G$, but this time $H$ contains neither
$u_1$ nor $w_1$ (although it contains the rest).  Any partial coloring
of $q-2$ vertices of $K_q$ has only 2 completions, so there are at
most 2 ways to extend any coloring of $H$ to include $u_1$ and $w_1$.
But then every other $u_i$ has at most 1 choice because $\{u_1, u_i\}
\cup (W \setminus \{w_i\})$ induces a $K_q$-subgraph in which $u_i$ is
the only uncolored vertex.  Therefore, every coloring of $H$ has at
most 2 extensions to $G$, as claimed.

It remains to consider the case when every set of $q-1$ vertices of
$A$ has degrees summing to at least $\delta + {q-1 \choose 2}$.  We
will show that then $G$ has fewer colorings than $T_{q-1}(n)$, which
is impossible.  Let $B = V(G) \setminus A$.  By an averaging argument,
the sum of degrees of $A$ is at least $\frac{q}{q-1} \big[ \delta +
{q-1 \choose 2} \big]$. Since $|A|=q$, the number of edges between $A$ and $B$ is
at least $\frac{q}{q-1} \big[ \delta + {q-1 \choose 2} \big] - 2 {q
  \choose 2}$.

Let $B_0$ be the set of isolated vertices of $G$, and for $2 \leq i
\leq q-1$, let $B_i$ be the set of vertices of $B$ that send $i$ edges
to $A$.  Note that no vertex can send $q = |A|$ edges to $A$ because
that would create a $K_{q+1}$-subgraph, making $G$ not $q$-colorable.
So, if we let $B_1 = B \setminus (B_0 \cup B_2 \cup \cdots \cup
B_{q-1})$, then every vertex of $B_1$ either sends exactly 1 edge to
$A$, or it is a non-isolated vertex that sends no edges to $A$.  Let
$b_i = |B_i|$.  By counting the number of edges between $A$ and $B$,
we obtain:
\begin{equation}
  \label{ineq:sum-ibi}
  \sum_{i=1}^{q-1} i b_i 
  \ \ \geq \ \ 
  \frac{q}{q-1} \left[ \delta + {q-1 \choose 2} \right] 
  - 2 {q \choose 2}.
\end{equation}

We now bound the number of $q$-colorings of $G$ in terms of the $b_i$.
There are exactly $q!$ ways to color $A$ because it induces $K_q$.
Then, there are exactly $q^{b_0}$ ways to extend this partial coloring
to $B_0$ because each isolated vertex has a free choice of the $q$
colors.  Next, for every $i \in \{2, \ldots, q-1\}$, each vertex in
$B_i$ has at most $q-i$ color choices left because it is adjacent to
$i$ vertices in $A$, all of which received different colors since
$G[A] = K_q$.  Finally, we color the vertices of $B_1$ by considering
them in an order such that whenever we color a vertex, it always has a
neighbor that we already colored.  This is possible because our
initial technical adjustment allows us to assume that $G$ has only one
nontrivial connectivity component.  Hence each vertex in $B_1$ will
have at most $q-1$ choices.  Putting this all together, we find that
the number of $q$-colorings of $G$ is at most
\begin{displaymath}
  q! \cdot \prod_{i=0}^{q-1} (q-i)^{b_i}
  \ \ \leq \ \ 
  q! \cdot \prod_{i=0}^{q-1} 2^{(q-i-1)b_i} 
  \ \ \leq \ \ 
  q! \cdot 2^{(q-1)(n-q)} \cdot 2^{
    -\frac{q}{q-1} \left[ \delta + {q-1 \choose 2} \right] 
    + 2 {q \choose 2}},
\end{displaymath}
where we used the inequality $x+1 \leq 2^x$ for $x \in \mathbb{Z}$,
the identity $\sum b_i = n-q$ (since $\cup B_i = V(G) \setminus A$),
and the bound for $\sum i b_i$ from inequality \eqref{ineq:sum-ibi}.
Inequality \ref{ineq:exact:turan:routine} routinely verifies that this
final bound is always strictly less than the number of colorings of
$T_{q-1}(n)$, contradicting our assumption that $G$ had at least that
many colorings.  \hfill $\Box$

\vspace{3mm}

\noindent \textbf{Proof of Theorem \ref{thm:exact:turan}.}\, Let $q
\geq 4$ be fixed, and let $N$ be the corresponding minimum number of
vertices for which Lemma \ref{lem:exact:turan} holds (it is valid only
for sufficiently large $n$).  We will show that Theorem
\ref{thm:exact:turan} holds for all $n \geq q {N \choose 2}$.  So,
suppose for contradiction that $G \neq T_{q-1}(n)$ is an $n$-vertex
graph with at least as many edges and $q$-colorings as $T_{q-1}(n)$.

Define a sequence of graphs as follows.  Start with $G_0 = G$.  If
$G_i$ is the current graph, stop if $G_i$ has fewer colorings than the
$(q-1)$-partite Tur\'an graph with $n-(q-1)i$ vertices.  Otherwise,
let $G_{i+1}$ be the graph $H$ obtained by applying Lemma
\ref{lem:exact:turan} to $G_i$.  We claim that this process terminates
before the graph $G_i$ has fewer than $N$ vertices, so we will always
be able to apply the lemma.  Indeed, each $G_i$ has exactly $n-(q-1)i$
vertices, so it will take more than ${N \choose 2}$ iterations before
$G_i$ has fewer than $N$ vertices.  Yet if $\Delta \geq 0$ is the
difference between the number of edges of $G$ and $T_{q-1}(n)$, then
each $G_i$ has at least $\Delta+i$ more edges than the $(q-1)$-partite
Tur\'an graph with $n-(q-1)i$ vertices.  So, after ${N \choose 2}$
iterations, $G_i$ would certainly have more than the maximum number of
edges of an $N$-vertex graph, and we indeed can never reach a graph
with fewer than $N$ vertices.

Therefore, we stop at some $G_t$, which has $n' = n-(q-1)t$ vertices
and fewer colorings than $T_{q-1}(n')$, but at least $2^{-t}$ times as
many colorings as $G$.  Divide $n$ by $q-1$, so that $n = s(q-1) + r$
with $0 \leq r < q-1$, and note that $n' = (s-t)(q-1) + r$.  Lemma
\ref{lem:turan-number-colorings} calculates that $T_{q-1}(n')$ has
exactly $q! \cdot \big[ (q-1+r)2^{s-t-1} - q + 2 \big]$ colorings, so
$G$ has at most $2^t$ times that many, hence fewer than $q! \cdot
\big[ (q-1+r)2^{s-1} - q + 2 \big]$.  Yet by the same lemma, that
final bound equals the number of colorings of $T_{q-1}(n)$.  Thus $G$
has fewer colorings than $T_{q-1}(n)$, contradiction.  \hfill $\Box$

\section{Concluding remarks}

\begin{itemize}
\item We have developed an approach that we hope future researchers
  can use to determine the graphs that maximize the number of
  $q$-colorings.  Theorems \ref{thm:asymp-number} and
  \ref{thm:asymp-stability} reduce any instance of this problem to a
  quadratically-constrained linear program, which can be solved for
  any case of interest.  Thus, thanks to modern computer algebra
  packages, these theorems imply that for any fixed $q$, approximately
  determining the extremal graphs amounts to a finite symbolic
  computation.

  The remaining challenge is to find analytic arguments which solve
  the optimization problem for general $q$, and then refine the
  approximate structure into precise results.  We accomplished this
  for low densities $m/n^2$, and the natural next step would be to
  extend the result to the range $\frac{m}{n^2} \leq \frac{1}{4}$.  In
  this range, and for all $q$, we expect the solution to the
  optimization problem to correspond to a bipartite graph plus
  isolated vertices.  This common form gives hope that perhaps one can
  find a solution which works across all $q$.

\item For $q=3$, we also know the approximate form of the extremal
  graphs when $\frac{m}{n^2} > \frac{1}{4}$, since Proposition
  \ref{prop:solve-opt} solved the entire $q=3$ case of the
  optimization problem.  However, we did not pursue the precise
  structure of the optimal graphs because it appears that their
  description is substantially more involved, and this paper was
  already quite long.

%  [But perhaps we should not conjecture that the semi-complete graphs
%  are best.  They might not be the best, because, say 100 colors,
%  should split 50-50, and if we lose, say, $r^2$ edges, then maybe it
%  is better to split the loss by taking away $K_{r,r}$ from each side.
%  That will reduce us by $r$ many 50's on the left, and $r$ many 50's
%  on the right.  But if we take away a $K_{1,r^2}$, then the left
%  loses only one 50, but the right loses $r^2$ many 50's.]

\item Our methods in Section \ref{sec:reduction-to-opt} can easily be
  adapted to maximize the number of graph homomorphisms to an
  arbitrary $H$ (not just $K_q$).  The analogues of Theorems
  \ref{thm:asymp-number} and \ref{thm:asymp-stability} show that for
  any fixed $H$, the asymptotic maximum number of homomorphisms from
  an $n$-vertex, $m$-edge graph to $H$ can be determined by solving a
  certain quadratically-constrained linear program.  Although this can
  in principle be done, it appears that the computations become rather
  messy even for graphs $H$ of small order.

  However, in the interesting case when $H$ is the two-vertex graph
  consisting of a single edge plus a loop, one can easily determine
  the extremal graphs via a direct argument.  As we mentioned in the
  introduction, this corresponds to maximizing the number of
  independent sets.  By considering the complement of the graph, this
  is equivalent to maximizing the number of cliques.

  We claim that for any $n,m$, the same graph that Linial found to
  minimize the number of colorings also happens to maximize the number
  of cliques.  This graph $G^*$ was a clique $K_k$ with an additional
  vertex adjacent to $l$ vertices of the $K_k$, plus $n-k-1$ isolated
  vertices, where $k,l$ are the unique integers satisfying $m = {k
    \choose 2} + l$ with $k > l \geq 0$.  We will show that for any
  $t$, every $n$-vertex graph $G$ with $m$ edges has at most as many
  $t$-cliques as $G^*$.  The only nontrivial values of $t$ to check
  are $2 \leq t \leq k$.

  If $l+2 \leq t \leq k$, then $G^*$ has exactly ${k \choose t}$
  cliques of size $t$.  Suppose for contradiction that $G$ has more
  $t$-cliques.  Construct a $t$-uniform hypergraph with at least ${k
    \choose t} + 1 = {k \choose t} + {t-1 \choose t-1}$ hyperedges by
  defining a hyperedge for each $t$-clique.  By the Kruskal-Katona
  theorem (see, e.g., the book \cite{Kruskal-Katona}), the number of
  2-sets that are contained in some hyperedge is at least ${k \choose
    2} + {t-1 \choose 1} \geq {k \choose 2} + (l+1)$, which exceeds
  the number of edges of $G$.  This contradicts the definition of the
  hyperedges, because each of these 2-sets must be an edge of $G$.

  On the other hand, if $2 \leq t \leq l+1$, $G^*$ has exactly ${k
    \choose t} + {l \choose t-1}$ cliques of size $t$.  A similar
  argument shows that if $G$ has at least ${k \choose t} + {l \choose
    t-1} + 1 = {k \choose t} + {l \choose t-1} + {t-2 \choose t-2}$
  cliques of size $t$, then $G$ must have at least ${k \choose 2} + {l
    \choose 1} + {t-2 \choose 0} \geq {k \choose 2} + l + 1$ edges,
  contradiction.

  Therefore, $G^*$ indeed maximizes the number of cliques.
  Furthermore, we can classify all extremal graphs, because our
  argument shows that any other graph $G$ with as many cliques as
  $G^*$ must also have exactly the same number of $t$-cliques for all
  integers $t$.  In particular, using $t=k$, we see that $G$ must also
  contain a $K_k$.  If $l \neq 1$, we can use $t=l+1$ to conclude that
  the remaining edges form a star with all endpoints in the $K_k$.
  Therefore, the maximizer is unique unless $l = 1$, in which case the
  extremal graphs are $K_k$ plus an arbitrary edge (not necessarily
  incident to the $K_k$).

%   Therefore, it is impossible to do better than $G^*$.  Our arguments
%   also show that if a graph $G$ has the same number of cliques as
%   $G^*$, then it must also have the same number of $t$-cliques for all
%   $t$.  In particular, $G$ must also have a copy of $K_k$.
%   \textbf{When $l \geq 2$.}  The remaining $l$ edges must form a star
%   with all endpoints in the $K_k$ because $G$ must have exactly ${k
%     \choose l+1} + {l \choose l}$ cliques of size $l+1$.

\end{itemize}

\appendix

\section{Routine verifications for Optimization Problem 2}

In this section, we present the postponed proofs of the results stated
in Section \ref{sec:solve-opt-sparse-q>=9}.  We begin by disposing of
Lemma \ref{lem:sparse:Fq}, which states some analytical facts about
the function $F_q(x) =\log \frac{q}{q-x} \cdot \log \frac{q}{x}$.

\vspace{3mm}

\begin{figure}[htbp]
  \centering
  \begin{minipage}[h]{0.4\linewidth}
    \begin{center}
    \includegraphics[scale=0.6]{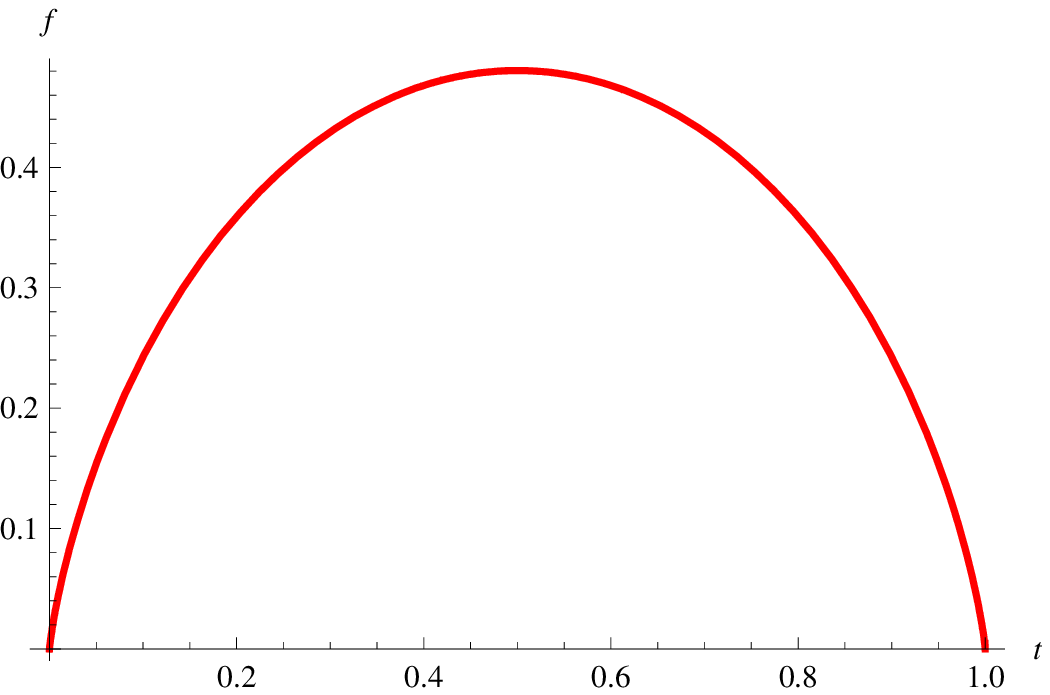}

    \noindent \textbf{Plot (i)}
    \end{center}
  \end{minipage}
  \hspace{0.05\linewidth}
  \begin{minipage}[h]{0.4\linewidth}
    \begin{center}
    \includegraphics[scale=0.6]{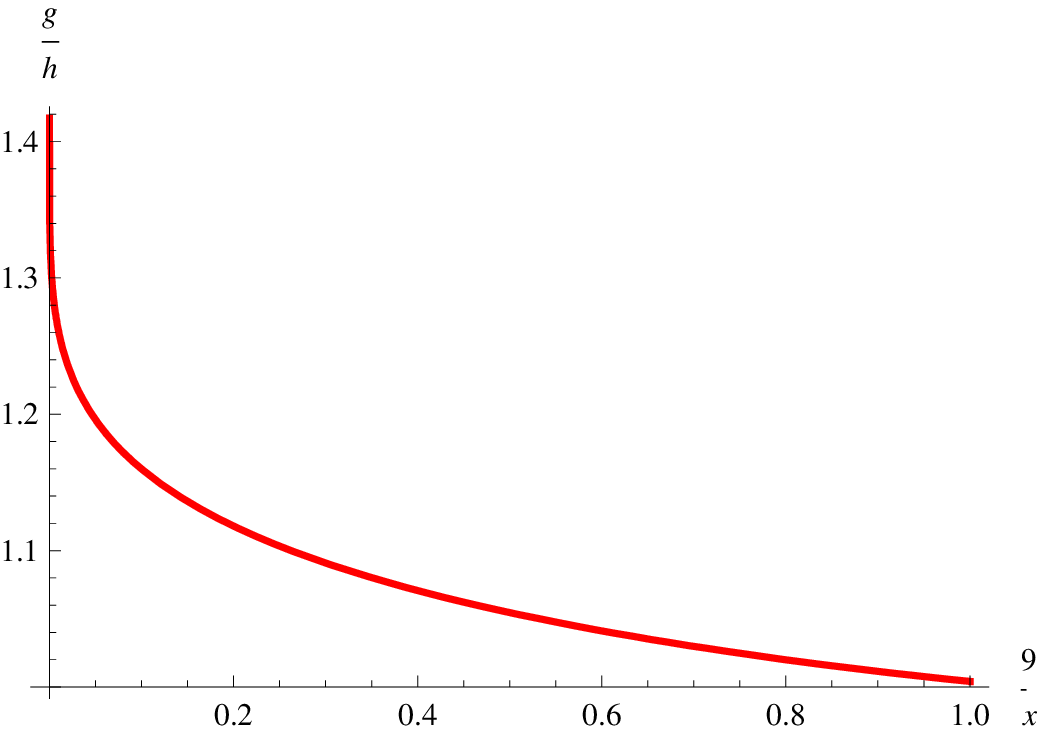}

    \noindent \textbf{Plot (ii)}
    \end{center}
  \end{minipage}
  \caption{\footnotesize Plot (i) displays the function $f(t) = \log
    \frac{1}{1-t} \log \frac{1}{t}$.  Plot (ii) displays the ratio
    $g(x)/h(x)$, where $g$ and $h$ are as defined above, and the
    horizontal axis is parameterized by $9/x$.}
  \label{fig:sparse:Fq}
\end{figure}

\noindent \textbf{Proof of Lemma \ref{lem:sparse:Fq}.}\, For part (i),
observe that if we reparameterize with $t = x/q$, then we need to show
that the function $f(t) = \log \frac{1}{1-t} \log \frac{1}{t}$ is
strictly increasing on $0 < t < 1/2$ and strictly decreasing on $1/2 <
t < 1$.  Instead of presenting a tedious analytic proof (which is
routine and not very interesting), we refer the reader to
\emph{Mathematica}'s plot of $f(t)$ in Figure \ref{fig:sparse:Fq}(i).

For part (ii), define the functions $g(x) = F_x(3) = \log
\frac{x}{x-3} \log \frac{x}{3}$ and $h(x) = 2F_x(1) \cdot
\frac{x-3}{x-2} = 2 \cdot \log \frac{x}{x-1} \log x \cdot
\frac{x-3}{x-2}$.  We need to show that $g(x) > h(x)$ for all $x \geq
9$.  Direct substitution yields $g(9) \approx 0.4454$ and $h(9)
\approx 0.4437$, so it is true at $x=9$.

Also, a quick estimate shows that asymptotically, as $x \rightarrow
\infty$, $g(x) = \log \big(1 + \frac{3}{x-3}\big)\cdot \log \frac{x}{3} =
(1+o(1)) \frac{3}{x} \cdot \log x$ and $h(x) = 2 \cdot \log \big(1 +
\frac{1}{x-1}\big)\cdot \log x \cdot \frac{x-3}{x-2} = (2+o(1)) 
\frac{1}{x} \cdot \log x$.  Therefore, the ratio $g(x) / h(x)$ tends
to 1.5, which is indeed greater than 1.

Again, instead of writing a routine analytic proof to fill in the gap
between 9 and infinity, we refer the reader to Figure
\ref{fig:sparse:Fq}(ii), which shows that the ratio $g/h$ steadily
increases as $x$ grows from 9.  Thus, $g(x) > h(x)$ for all $x \geq
9$, as required. \hfill $\Box$

\vspace{3mm}

The monotonicity of $F_q(x)$ on $0 < x < q/2$, which we just
established, is useful for our next proof.  This is Lemma
\ref{lem:sparse:2sets=>done}, which stated that if $\blpha$ solves
$\optS$ and is supported by a partition of $[q]$ consisting of exactly
two sets, then $\blpha$ must have the same form as $\blpha^*$, the
claimed optimal vector in Proposition \ref{prop:solve-opt-2}.

\vspace{3mm}

\noindent \textbf{Proof of Lemma \ref{lem:sparse:2sets=>done}.}\, Let
$A$ and $B$ denote the two sets in the support, with $|A| \leq |B|$.
Write $a = |A|$.  Flipping the fractions to make the logarithms
positive, we have $\objS(\blpha) = -\alpha_A \log \frac{q}{a} -
\alpha_B \log \frac{q}{q-a} \leq -2 \sqrt{ \alpha_A \log \frac{q}{a}
  \cdot \alpha_B \log \frac{q}{q-a} }$ by the inequality of arithmetic
and geometric means.  Yet $\alpha_A \alpha_B = \ec(\blpha) \geq 1$
since $\blpha$ is in the feasible set $\feasS$, so $\objS(\blpha) \leq
-2 \sqrt{ \log \frac{q}{a} \cdot \log \frac{q}{q-a}} =
-2\sqrt{F_q(a)}$.  Here, $F_q$ is the function which Lemma
\ref{lem:sparse:Fq}(i) claimed was strictly increasing between 0 and
$q/2$.  In particular, since $1 \leq a \leq q/2$, the final bound is
at most $-2 \sqrt{F_q(1)}$, which we recognize as $\objS(\blpha^*)$,
where $\blpha^*$ is the claimed unique optimal vector in Proposition
\ref{prop:solve-opt-2}.

Since $\blpha$ was assumed to be maximal, we must have equality in all
of the above inequalities.  Checking the equality conditions, we find
that $\blpha$ must indeed have the unique form claimed in Proposition
\ref{prop:solve-opt-2}. \hfill $\Box$

\vspace{3mm}

The remaining lemma from Section \ref{sec:solve-opt-sparse-q>=9} ruled
out a handful of partitions as possible supports for optimal vectors.
It turns out that each of those excluded partitions is a special case
of the following result.

\begin{lemma}
  \label{lem:sparse:almost-singletons}
  Fix any integer $q \geq 3$, and let $\blpha$ be a vector which
  solves $\optS$, whose support is a partition of $[q]$.  Then that
  partition cannot be $\{1, \ldots, t\} \cup \{t+1\} \cup \{t+2\} \cup
  \ldots \cup \{q\}$, where $1 \leq t \leq q-2$.
\end{lemma}

\noindent \textbf{Proof.}\, Assume for the sake of contradiction that
$\blpha$ is supported by the above partition.  Let $x =
\alpha_{\{t+1\}} = \cdots = \alpha_{\{q\}}$, which are all equal by
Lemma \ref{lem:sparse:IJ}(ii).  We assumed that $\blpha$ was maximal,
so in particular $\objS(\blpha) \geq \objS(\blpha^*) = -2\sqrt{\log
  \frac{q}{q-1} \log q}$, where $\blpha^*$ is the feasible vector
constructed in Proposition \ref{prop:solve-opt-2}.  Therefore,
\begin{displaymath}
  (q-t)x \log \frac{1}{q} 
  \ \ > \ \ \alpha_{\{1, \ldots, t\}} \log \frac{t}{q} + (q-t)x \log \frac{1}{q}
  \ \ = \ \ \objS(\blpha) \ \ \geq \ \ -2 \sqrt{\log \frac{q}{q-1} \log q},
\end{displaymath}
and we conclude that $(q-t)x < 2\sqrt{\log \frac{q}{q-1} / \log q}$.
On the other hand, we also know by Lemma \ref{lem:sparse:IJ}(i) for
the set $A = \{1, \ldots, t\}$ that $(q-t)x = I_A/\alpha_A = 2J_A/\alpha_A =
\big(2\log \frac{t}{q}\big) / \objS(\blpha)$.  Using the final bound
for $(q-t)x$ above, this gives
\begin{displaymath}
  \objS(\blpha) = \left(2 \log \frac{t}{q}\right) / ((q-t)x)
  \ \ < \ \ \log \frac{t}{q} \cdot \sqrt{(\log q) / \log \frac{q}{q-1}}.
\end{displaymath}
(The inequality reversed because $\log \frac{t}{q}$ is negative.)

To get our contradiction, it remains to show that this is less than
$\objS(\blpha^*) = -2\sqrt{\log \frac{q}{q-1} \log q}$.  Cancelling
the common factor of $\sqrt{\log q}$ and rearranging terms, this
reduces to showing that $\log \frac{q}{t} > 2 \log \frac{q}{q-1}$.

Since $t \leq q-2$ by definition, it suffices to show that $\log
\frac{q}{q-2} > 2 \log \frac{q}{q-1}$.  Removing the logarithms
reduces us to showing that $\frac{q}{q-2} > \frac{q^2}{(q-1)^2}$.
This is equivalent to $(q-1)^2 > q(q-2)$, which is easily seen to be
true by multiplying out each side. \hfill $\Box$

\vspace{3mm}

\noindent \textbf{Proof of Lemma \ref{lem:sparse:extreme}}.\, Part
(i), the partition of all singletons, is precisely the case of the
previous lemma when $t=1$.  Similarly, part (ii), the partition of all
singletons except for a 2-set, corresponds to the $t=2$ case.  For
part (iii), which concerns partitions that include a $(q-2)$-set,
first note that if the partition is a $(q-2)$-set plus two singletons,
then it is precisely the $t=q-2$ case of the previous lemma.  The only
other possibility is that the partition is a $(q-2)$-set plus a 2-set,
and this is excluded by Lemma \ref{lem:sparse:2sets=>done}.  \hfill
$\Box$

\section{Routine verifications for exact results}

% In this section...

\begin{proposition}
  \label{prop:pendant-necessary}
  Let $r$ be a sufficiently large positive integer.  Then the complete
  bipartite graph $K_{r,2r}$ plus one pendant edge achieves the
  maximum number of colorings among all $(3r+1)$-vertex graphs with
  $2r^2 + 1$ edges.
\end{proposition}

\noindent \textbf{Proof.}\, Every 3-coloring of $K_{r,2r}$ has exactly
2 extensions to the pendant vertex, so Lemma
\ref{lem:subgraph-bipartite:q=3} shows that the above graph has
exactly $\big(3 \cdot 2^r + 3 \cdot 2^{2r} - 6\big) \cdot 2 = (1+o(1))
\cdot 3 \cdot 2^{2r+1}$ colorings.  Plugging $n = 3r+1$ and $m = 2r^2
+ 1$ into the dense case of Theorem \ref{thm:main:q=3}, we see that
the only other graphs we need to consider are semi-complete subgraphs
of some $K_{a,b}$ with $a = (1+o(1))r$ and $b = (2+o(1))r$, plus
isolated vertices.  Note that we must have $a \geq r$, because when $a
\leq r-1$ and $a+b \leq 3r+1$, convexity implies that $ab \leq
(r-1)(2r+2) = 2r^2 - 2 < 2r^2 + 1$, and there would not be enough
edges.

Let $G'$ be one of the above graphs with $a = r+t$ for some $t \geq
0$.  We must have $b \geq 2r-2t+1$, because $(r+t)(2r-2t) = 2r^2 -
2t^2 < 2r^2 + 1$, so any smaller $b$ would not produce enough edges.
This leaves $n - a - b \leq t$ isolated vertices.  Observe that when
$t=0$, this forces $G'$ to be a semi-complete subgraph of $K_{r,
  2r+1}$ with exactly $r-1$ missing edges.  Lemma
\ref{lem:subgraph-bipartite:q=3} then shows that the number of
colorings of $G'$ is $3 \cdot 2^r + 3 \cdot 2^{2r+1} + 6 \cdot
\big(2^{r-1} - 2\big)$, which is exactly the same as $G$.

It remains to consider $t > 0$.  By definition, any semi-complete
subgraph of $K_{a,b}$ is missing at most $a-1$ edges, so Lemma
\ref{lem:subgraph-bipartite:q=3} implies that the number of
3-colorings of $G'$ is at most $3^{n-a-b} \cdot \big(3 \cdot 2^a + 3
\cdot 2^b + 6 \cdot \big(2^{a-1} - 2\big) \big)$.  This expression is
largest when $b$ is as small as possible, so using $b \geq 2r-2t+1$
and $n = 3r+1$, we find that $G'$ has at most $3^t \cdot \big(3 \cdot
2^a + 3 \cdot 2^{2r-2t+1} + 6 \cdot \big(2^{a-1} - 2\big) \big)$
colorings.  Since $a = (1+o(1))r$, this is at most
$\big(\big(\frac{3}{4}\big)^t + o(1)\big) \cdot 3 \cdot 2^{2r+1}$,
which is indeed less than the number of colorings of $G$ when $r$ is
large.  \hfill $\Box$

\vspace{3mm}

\noindent \textbf{Remark.}\, A similar argument shows that for any $c
\in \{0, \pm 1, \pm 2\}$ and large $r$, $K_{r, 2r+c}$ plus a pendant
edge is optimal among graphs with $3r+c+1$ vertices and $r(2r+c) + 1$
edges.  Interestingly enough, it can also be shown that these values
of $n,m$ are the only ones which produce optimal graphs that are not
semi-complete plus isolated vertices, when $n,m$ are large.
% WHY: suppose we have $a \leq b$ and take $G$ to be $K_{a,b}$ plus a
% pendant edge attached to the $b$-side.  By Lemma
% \ref{lem:subgraph-bipartite:q=3}, $G$ has exactly $2(3 \cdot 2^a + 3
% \cdot 2^b - 6) = 3 \cdot 2^{a+1} + 3 \cdot 2^{b+1} - 12$ colorings.
%
% Suppose $b \leq 2a-3$.  Then take the complete bipartite graph
% $K_{a-1, b+2}$, which has $ab+2a-b-2 \geq ab+1$ edges, good.  Yet it
% has at least $3 \cdot 2^{a-1} + 3 \cdot 2^{b+2} - 6$ colorings,
% where the $-6$ comes from the double-counted colorings that only use
% 2 colors at all.  This exceeds the number of colorings of $G$.
%
% On the other hand, suppose $b \geq 2a+3$.  Then take the complete
% bipartite graph $K_{a+1, b-2}$ plus 2 isolated vertices.  Number of
% edges is $ab-2a+b-2 \geq ab+1$, which is enough.  Yet the number of
% colorings is $3^2 \cdot (3 \cdot 2^{a+1} + 3 \cdot 2^{b-2} - 6)$.
% The first term clearly beats the $3 \cdot 2^{a+1}$ of the number of
% colorings of $G$, and the second term is $27 \cdot 2^{b-2}$ compared
% to the $3 \cdot 2^{b+1} = 24 \cdot 2^{b-2}$ of the corresponding
% term in $G$.  So there is a gain as long as $b \geq 6$, which
% certainly happens for sufficiently large $m,n$ because our
% approximate result shows that we will never have one side fixed
% below a constant (6) as $m,n$ grow.

\vspace{3mm}

\begin{inequality}
  \label{ineq:partition-colors}
  Let $a,b,t$ be positive integers, with $t \geq 3$ and $\frac{b}{a}
  \geq \log t / \log \frac{t-1}{t-2}$.  Then:
  \begin{description}
  \item[(i)] The product $i^a (t-i)^b$ falls by a factor of at least
    $1.5^a$ when $i$ increases by 1, for all $i \in \{1, \ldots,
    t-2\}$.
  \item[(ii)] If we further assume that $a$ is sufficiently large
    (depending only on $t$), then $\sum_{i=1}^{t-1} {t \choose i} i^a
    (t-i)^b \leq 1.1 \cdot t (t-1)^b$, i.e., the first summand
    dominates.
  \end{description}
\end{inequality}

\noindent \textbf{Proof.}\, When $i \in \{1, \ldots, t-2\}$ increases
by 1, $i$ grows by a factor of at most 2, but $t-i$ falls by at least
$\frac{t-1}{t-2}$.  Thus, the product $i^a (t-i)^b$ falls by a factor
of at least $\big(\frac{1}{2}\big)^a \big(\frac{t-1}{t-2}\big)^b =
\big(\frac{1}{2} \cdot \big(\frac{t-1}{t-2}\big)^{b/a}\big)^a \geq
\big(\frac{1}{2} \cdot t\big)^a$.  Since $t \geq 3$, this gives (i).

For part (ii), when $i$ increases by 1, the term ${t \choose i}$ in
the summand grows by a factor of at most $t$, but by (i) the rest of
the summand falls by a factor of at least $1.5^a$.  Thus for
sufficiently large $a$, each successive term of the sum falls by a
factor of at least $1.4^a > 20$.  The result follows by bounding the
sum by a geometric series, since $1 + \frac{1}{20} + \frac{1}{20^2} +
\cdots < 1.1$.  \hfill $\Box$

\vspace{3mm}

\begin{inequality}
  \label{ineq:exact:q=3:claim5}
  Let $m$, $n$, $t$, and $v_1$ be positive integers, with $m \leq
  n^2/4$ and $v_1 (n-v_1) \geq m-t$.  Let $s$ be the largest integer
  that satisfies $s(n-s) \geq m$.  Then $s \geq v_1 - \sqrt{t}$.
\end{inequality}

\noindent \textbf{Proof.}\, The inequality for $s$ rearranges to $s^2
- ns + m \leq 0$, so the quadratic formula implies that $s$ is
precisely $\big\lfloor \frac{n + \sqrt{n^2 - 4m}}{2} \big\rfloor$.
Similarly, the inequality for $v_1$ rearranges to $v_1^2 - n v_1 +
(m-t) \leq 0$, so the quadratic formula implies that $v_1 \leq
\big\lfloor \frac{n + \sqrt{n^2 - 4m + 4t}}{2} \big\rfloor$.
Therefore,
\begin{eqnarray*}
  v_1 - s &\leq& \left\lfloor \frac{n + \sqrt{n^2 - 4m + 4t}}{2} \right\rfloor
  - \left\lfloor \frac{n + \sqrt{n^2 - 4m}}{2} \right\rfloor \\
  &\leq& \left\lceil \frac{n + \sqrt{n^2 - 4m + 4t}}{2} 
    - \frac{n + \sqrt{n^2 - 4m}}{2} \right\rceil 
  \ \ = \ \ \left\lceil \frac{\sqrt{(n^2 - 4m) + 4t}
    - \sqrt{n^2 - 4m}}{2} \right\rceil.
\end{eqnarray*}
Since the function $\sqrt{x}$ is concave and we assumed $n^2 - 4m \geq
0$, this final bound is largest when $n^2 - 4m = 0$.  Therefore, $v_1
- s \leq \lceil \sqrt{t} \rceil$, which gives the claimed result.
\hfill $\Box$

\vspace{3mm}

%Now we prove a supporting result for the other flavor of exact result
%which we proved in Section \ref{sec:exact:turan}.  

%The following lemma
%calculates the number of colorings of $T_{q-1}(n)$.

\begin{lemma}
  \label{lem:turan-number-colorings}
  The number of $q$-colorings of the Tur\'an graph $T_{q-1}(n)$ is
  exactly $q! \cdot \big[ (q-1+r)2^{s-1} - q + 2 \big]$, where $s$ and
  $r$ are defined by $n = s(q-1) + r$ with $0 \leq r < q-1$.
\end{lemma}

\noindent \textbf{Proof.}\, The complete $(q-1)$-partite graph
$T_{q-1}(n)$ has $r$ parts of size $s+1$ and $q-1-r$ parts of size
$s$, and any $q$-coloring must use different colors on each part.  The
number of $q$-colorings that use exactly one color on each part is
exactly $q \cdot (q-1) \cdots 2 = q!$.  All other colorings use 2
colors on one part, and one color on each of the other parts.  There
are ${q \choose 2}$ ways to choose which two colors are paired.  If
the pair of colors is used on one of the $r$ parts of size $s+1$, then
there are $2^{s+1} - 2$ ways to color that part with exactly 2 colors,
followed by $(q-2)!$ ways to choose which color goes to each of the
remaining parts.  Otherwise, if the pair of colors appears on one of
the $q-1-r$ parts of size $s$, then there $(2^s - 2) (q-2)!$ colorings
of this form.  Therefore, the number of $q$-colorings of $T_{q-1}(n)$
is exactly
\begin{displaymath}
  q! + {q \choose 2} \cdot \left[r \cdot (2^{s+1} - 2) (q-2)! +
    (q-1-r) \cdot (2^s - 2) (q-2)! \right] 
  \ \ = \ \  q! \cdot \left[ (q-1+r)2^{s-1} - q + 2 \right],
\end{displaymath}
as claimed.  \hfill $\Box$

\vspace{3mm}

\begin{inequality}
  \label{ineq:exact:turan:routine}
  Fix any $q \geq 4$.  For all sufficiently large $n$, the number of
  $q$-colorings of the Tur\'an graph $T_{q-1}(n)$ is strictly greater
  than
  \begin{equation}
    \label{exp:exact:turan}
    q! \cdot 2^{(q-1)(n-q)} \cdot 2^{
    -\frac{q}{q-1}  \left[ \delta + {q-1 \choose 2} \right] 
    + 2 {q \choose 2}},
  \end{equation}
  where $\delta$ is the difference between the number of edges of
  $T_{q-1}(n)$ and $T_{q-1}(n-q+1)$.
\end{inequality}

\noindent \textbf{Proof.}\, Divide $n$ by $q-1$, so that $n = s(q-1) +
r$ with $0 \leq r < q-1$.  Then $T_{q-1}(n)$ has exactly $r$ parts of
size $s+1$ and $q-1-r$ parts of size $s$, and $T_{q-1}(n-q+1)$ is
obtained by deleting one vertex per part.  Each deleted vertex in a
part of size $s+1$ had degree $n-s-1$, while each deleted vertex in a
part of size $s$ had degree $n-s$.  Thus, the number of deleted edges
is $\delta = r(n-s-1) + (q-1-r)(n-s) - {q-1 \choose 2}$, where we had
to subtract the double-counted edges of the $K_{q-1}$ induced by the
set of deleted vertices.  Substituting this into
\eqref{exp:exact:turan} and using $n = s(q-1)+r$ to simplify the
expression, we obtain:
\begin{eqnarray*}
  q! \cdot 2^{(q-1)(n-q)} \cdot 2^{
    -\frac{q}{q-1} \left[ \delta + {q-1 \choose 2} \right] 
    + 2 {q \choose 2}}
  &=& q! \cdot 2^{(q-1)(n-q)} \cdot 2^{
    -\frac{q}{q-1} \left[ r(n-s-1) + (q-1-r)(n-s) \right] 
    + 2 {q \choose 2}} \\
  &=& %q! \cdot 2^{s + \frac{r}{q-1}}  
%  \ \ = \ \ 
  q! \cdot 2^s \cdot 2^{\frac{r}{q-1}}.
\end{eqnarray*}
It remains to show that this is strictly less than the number of
colorings of $T_{q-1}(n)$, which Lemma
\ref{lem:turan-number-colorings} calculated to be $q! \cdot \big[
(q-1+r)2^{s-1} - q + 2 \big] = (1-o(1)) \cdot q! \cdot 2^s \cdot
\frac{q-1+r}{2}$.  Here, the $o(1)$ term tends to zero as $n$ grows
(and $s = \big\lfloor \frac{n}{q-1} \big\rfloor$ grows).  Recall that
$0 \leq r < q-1$, so when $r \geq 1$ and $q \geq 4$ we always have
$2^{\frac{r}{q-1}} < 2^1 \leq \frac{q-1+r}{2}$, giving the desired
result.  On the other hand, when $r = 0$, the result follows from
$2^{\frac{r}{q-1}} = 2^0 < \frac{3}{2} \leq \frac{q-1+r}{2}$.  \hfill
$\Box$

\section{\emph{Mathematica} computations for Optimization Problem 2}
\label{sec:mathematica}

The next 9 pages contain the complete \emph{Mathematica}\/ program
(and output), solving Optimization Problem 2 for $q < 9$.

\newpage
\pagestyle{empty}

\topmargin 0in
\headheight -0.5in
\headsep 0in

\centerline{\includegraphics{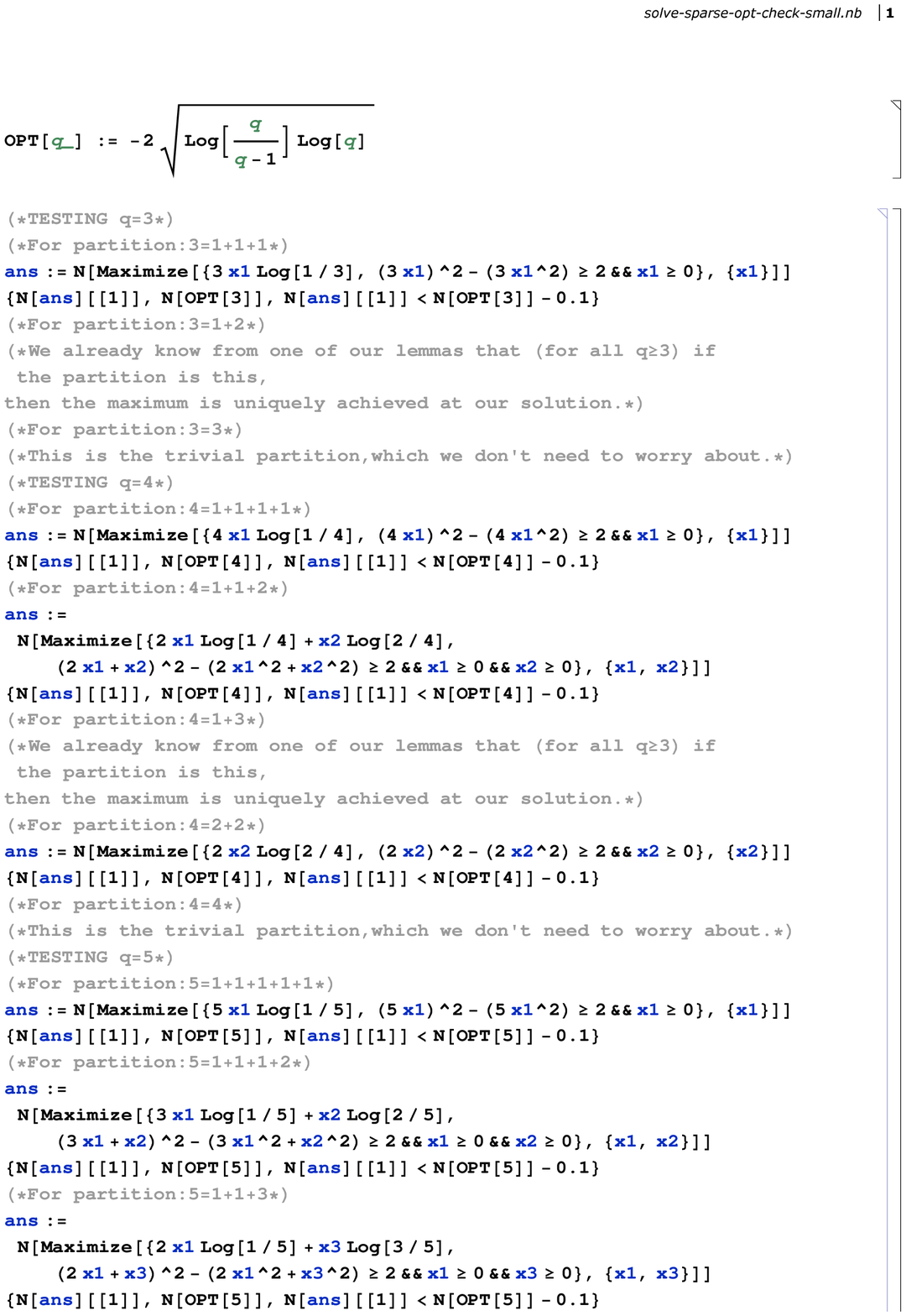}}

\newpage

\centerline{\includegraphics{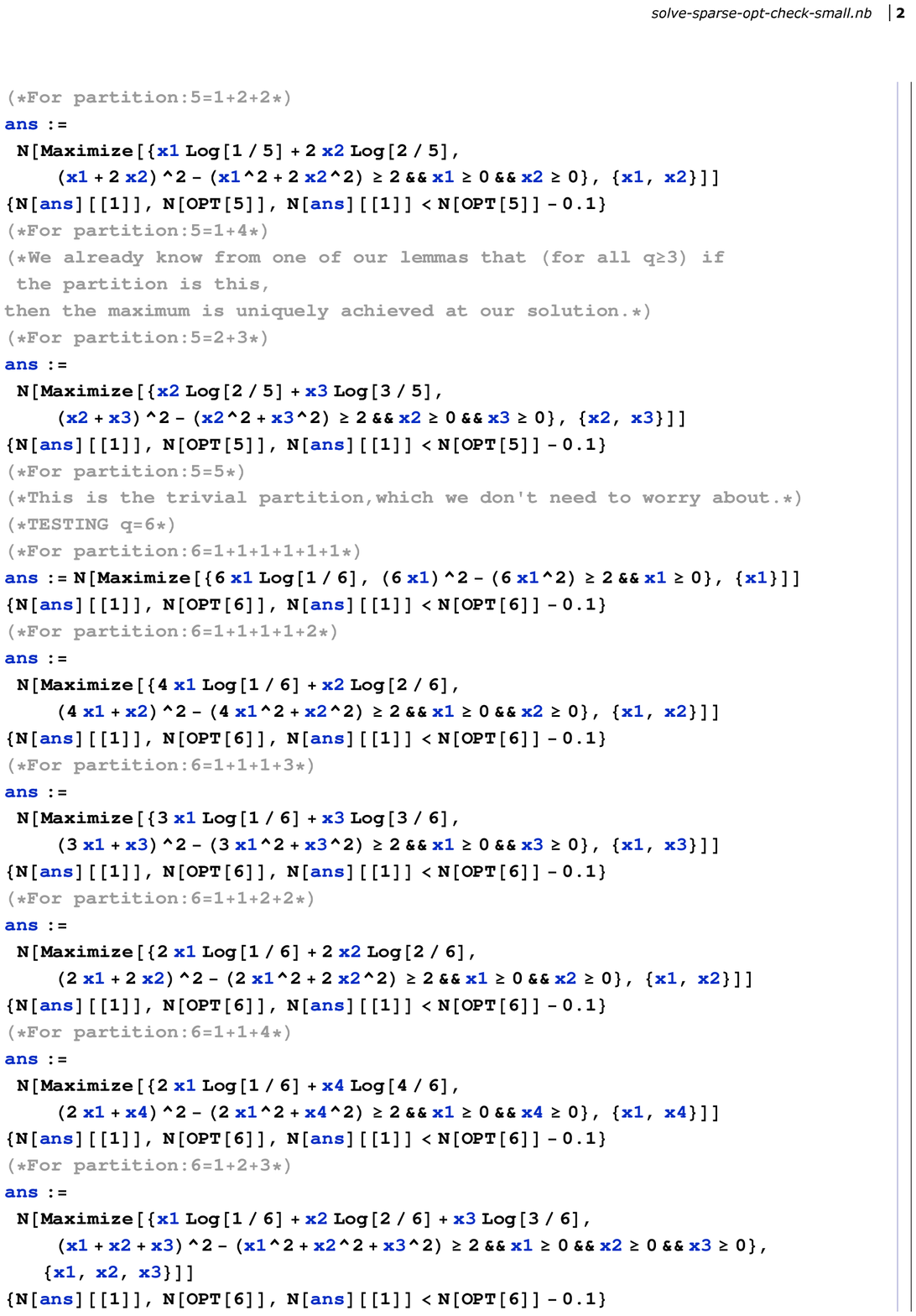}}

\newpage

\centerline{\includegraphics{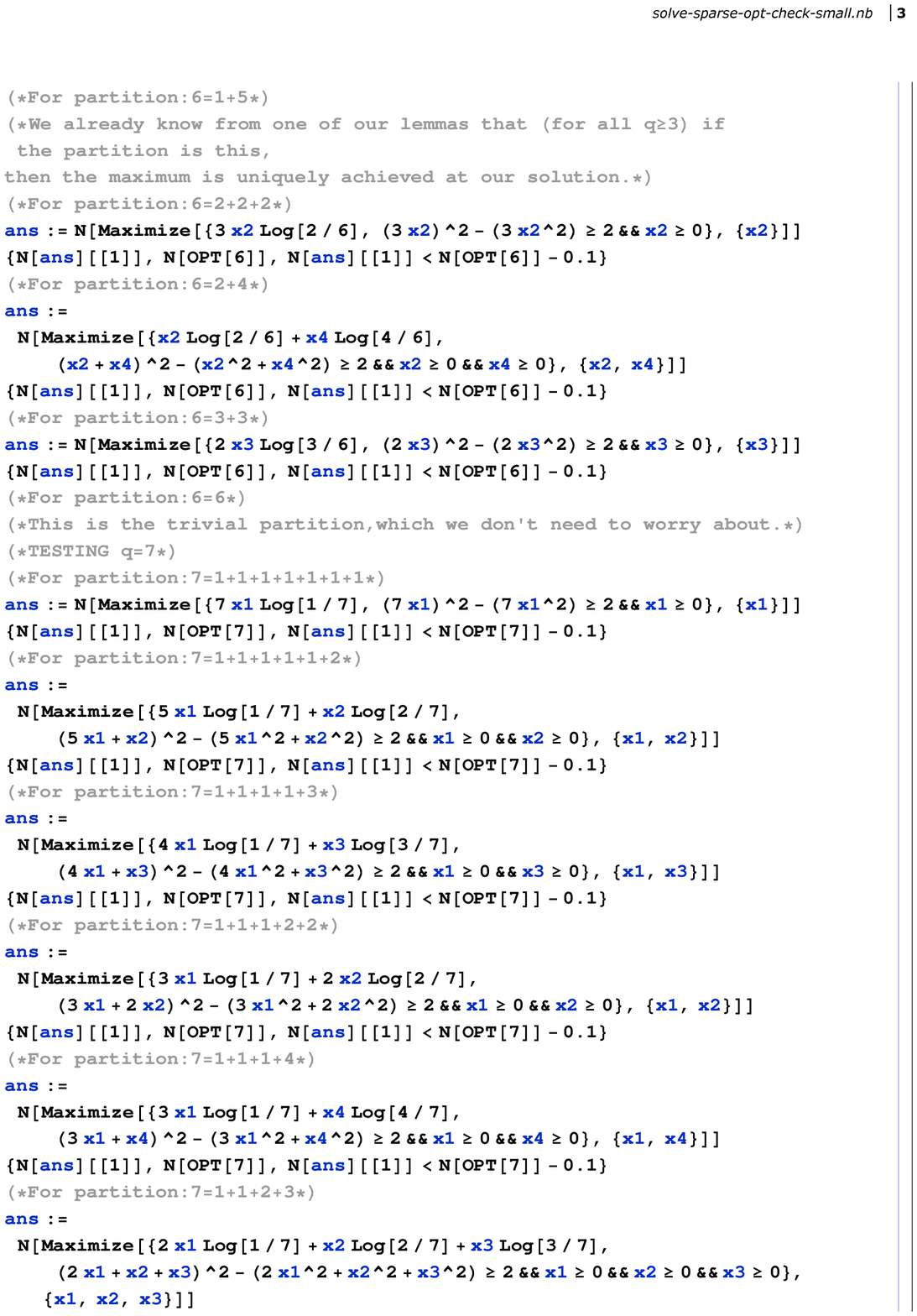}}

\newpage

\centerline{\includegraphics{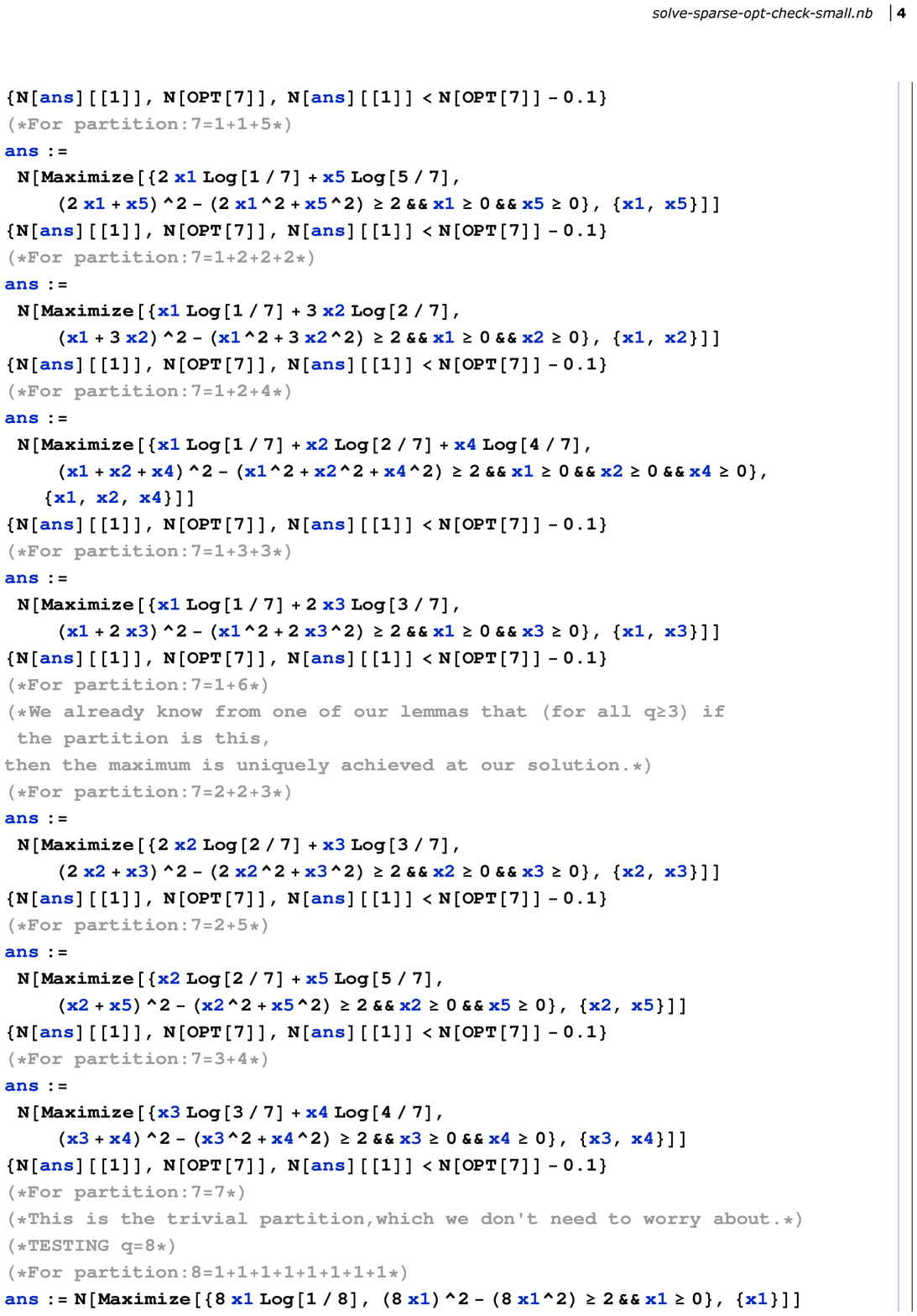}}

\newpage

\centerline{\includegraphics{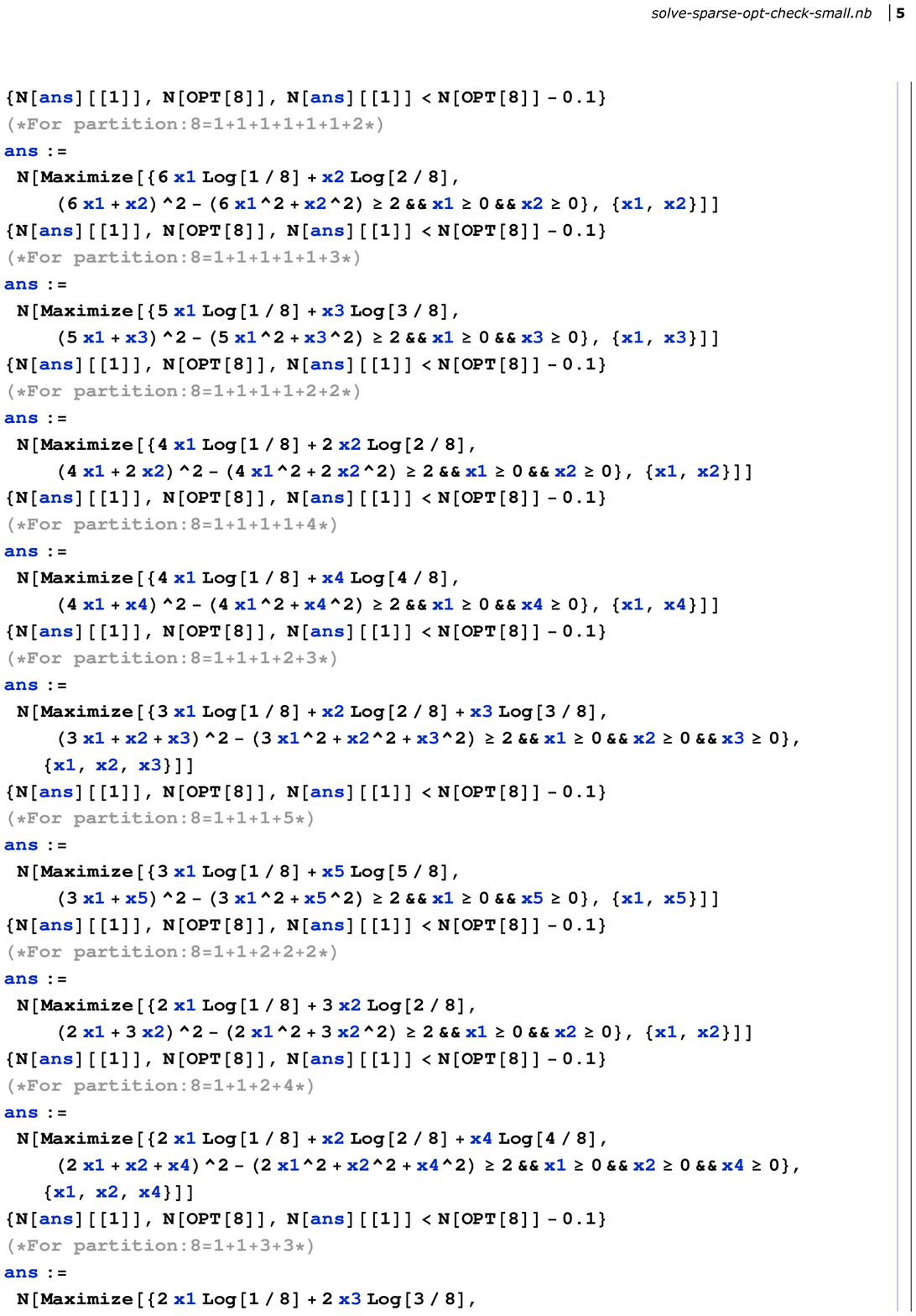}}

\newpage

\centerline{\includegraphics{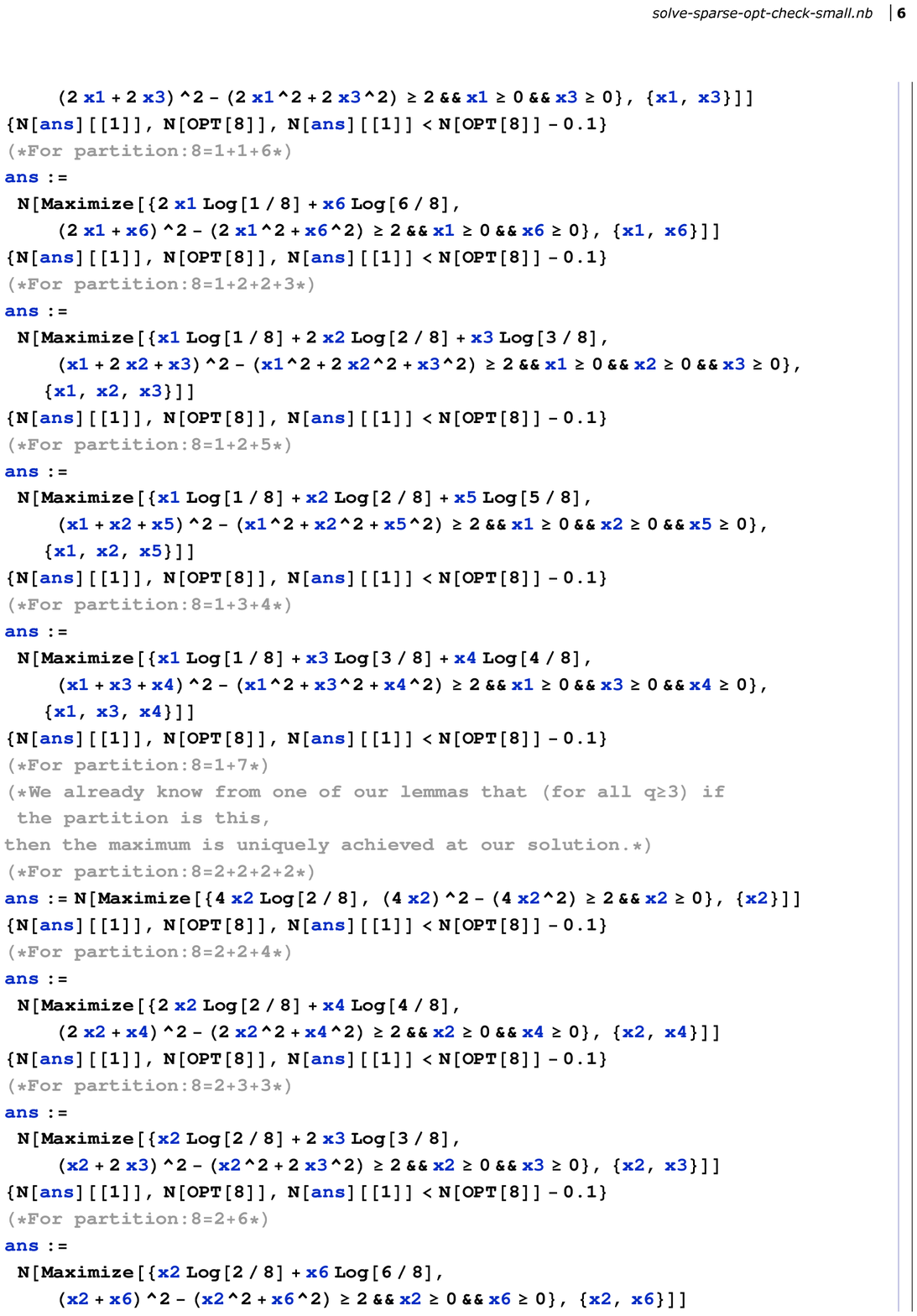}}

\newpage

\centerline{\includegraphics{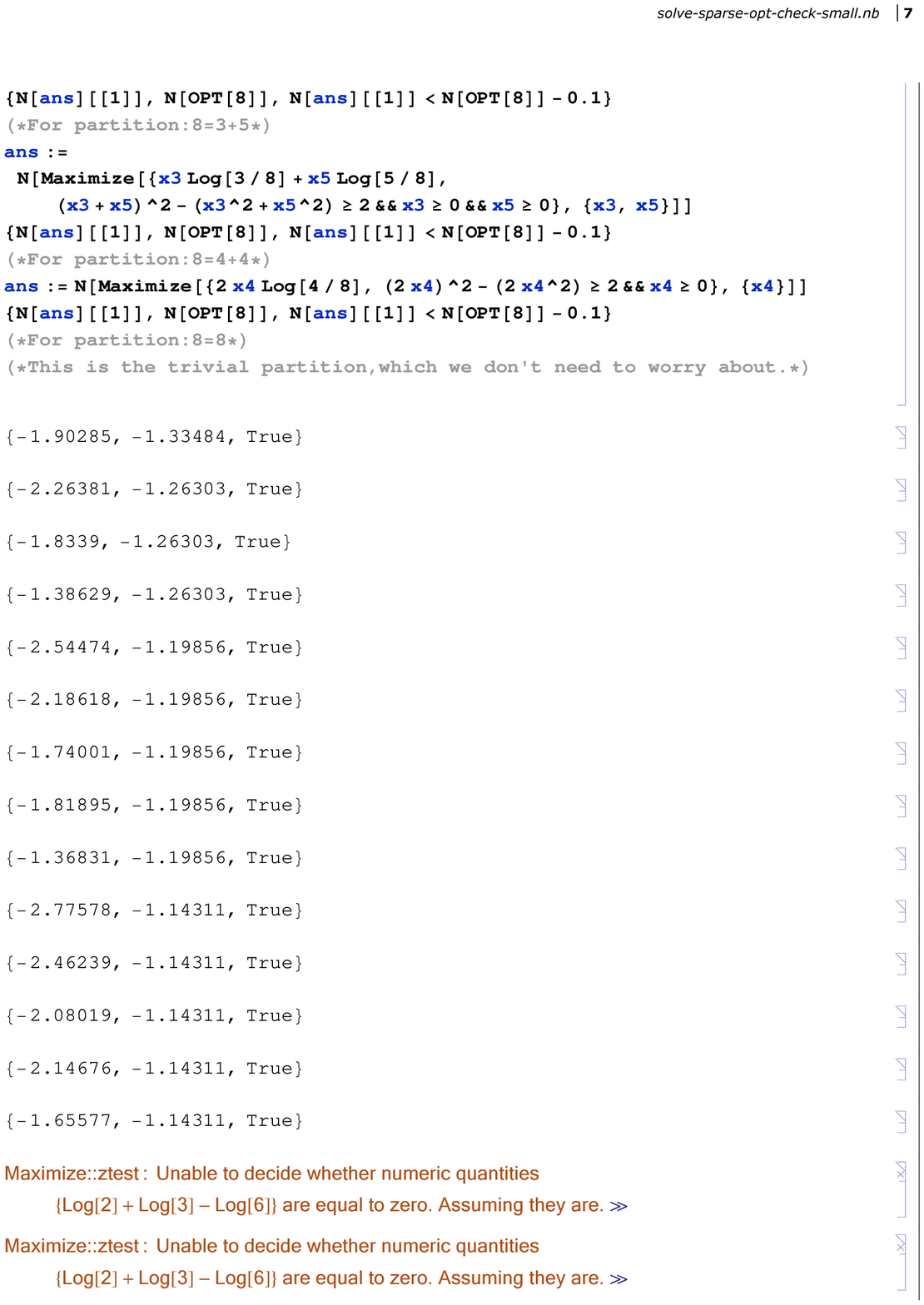}}

\newpage

\centerline{\includegraphics{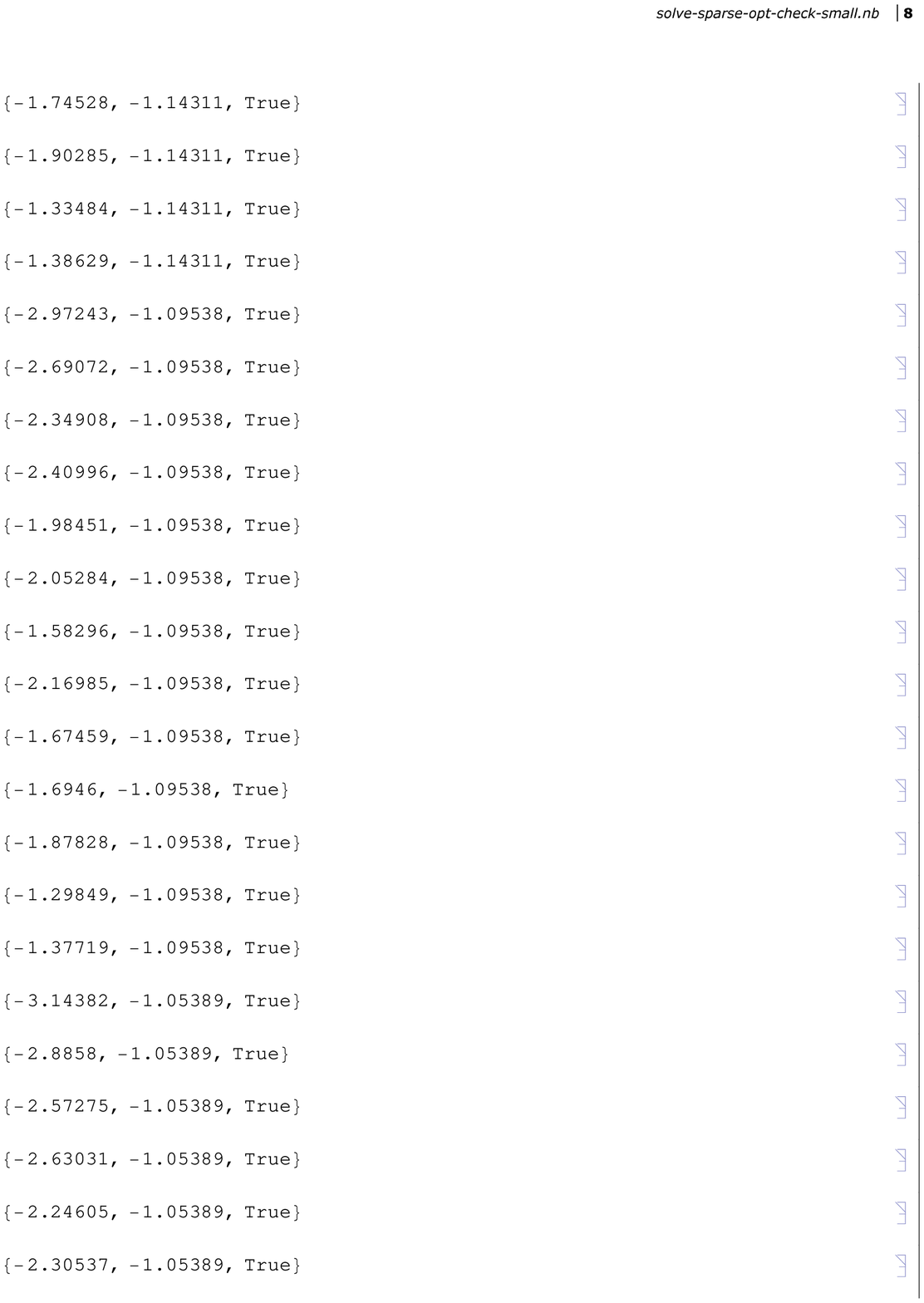}}

\newpage

\centerline{\includegraphics{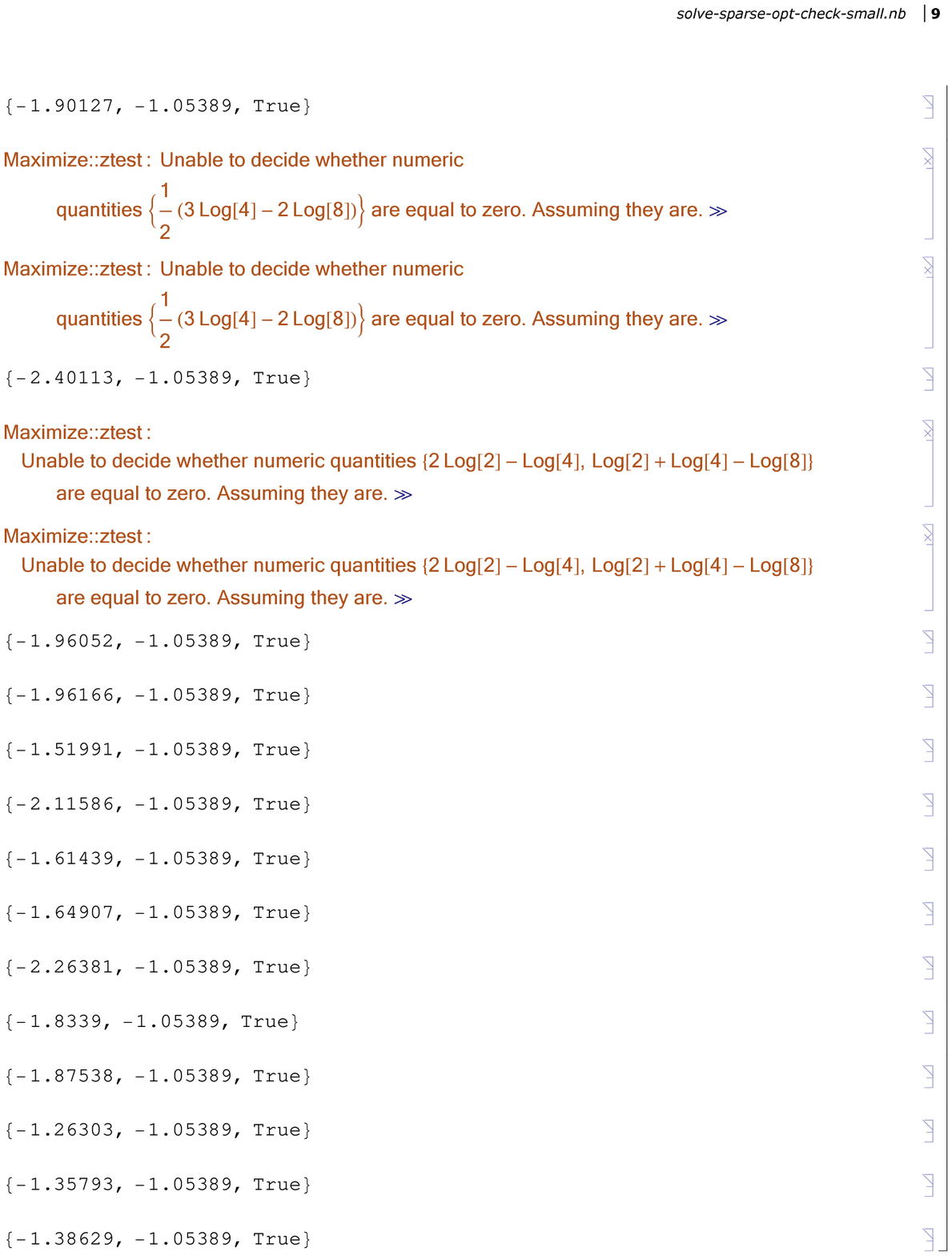}}

\end{document}